\documentclass[11pt,a4paper]{article}
\usepackage{graphicx}
\usepackage{mathrsfs}
\usepackage{color}
\usepackage[utf8]{inputenc}
\usepackage[francais, english]{babel}

\usepackage{libertine}
\usepackage[T1]{fontenc}
\usepackage{lipsum}
\usepackage{titlesec}
\usepackage[colorlinks=true,linkcolor=magenta,citecolor=blue]{hyperref}

\usepackage[cmex10]{amsmath}
\usepackage{stmaryrd}
\usepackage{color}
\usepackage{cmll}
\usepackage{amsmath,amssymb}
\usepackage[sloped]{fourier}
\usepackage{amsthm}
\usepackage{esint}
\usepackage{psfrag}

\author{S\'ebastien Alvarez}
\date{}
\title{Gibbs measures for foliated bundles with negatively curved leaves}

\begin{document}

\newtheorem{maintheorem}{Theorem}
\newtheorem{maincoro}[maintheorem]{Corollary}
\renewcommand{\themaintheorem}{\Alph{maintheorem}}
\newcounter{theorem}[section]
\newtheorem{exemple}{\bf Exemple \rm}
\newtheorem{exercice}{\bf Exercice \rm}
\newtheorem{conj}[theorem]{\bf Conjecture}
\newtheorem{defi}[theorem]{\bf Definition}
\newtheorem{lemma}[theorem]{\bf Lemma}
\newtheorem{proposition}[theorem]{\bf Proposition}
\newtheorem{coro}[theorem]{\bf Corollary}
\newtheorem{theorem}[theorem]{\bf Theorem}
\newtheorem{rem}[theorem]{\bf Remark}
\newtheorem{ques}[theorem]{\bf Question}
\newtheorem{propr}[theorem]{\bf Property}
\newtheorem{question}{\bf Question}
\def\bp{\noindent{\it Proof. }}
\def\ep{\noindent{\hfill $\fbox{\,}$}\medskip\newline}
\renewcommand{\theequation}{\arabic{section}.\arabic{equation}}
\renewcommand{\thetheorem}{\arabic{section}.\arabic{theorem}}
\newcommand{\eps}{\varepsilon}
\newcommand{\disp}[1]{\displaystyle{\mathstrut#1}}
\newcommand{\fra}[2]{\displaystyle\frac{\mathstrut#1}{\mathstrut#2}}
\newcommand{\dif}{{\rm Diff}}
\newcommand{\homeo}{{\rm Homeo}}
\newcommand{\Per}{{\rm Per}}
\newcommand{\Fix}{{\rm Fix}}
\newcommand{\A}{\mathcal A}
\newcommand{\Z}{\mathbb Z}
\newcommand{\Q}{\mathbb Q}
\newcommand{\R}{\mathbb R}
\newcommand{\C}{\mathbb C}
\newcommand{\N}{\mathbb N}
\newcommand{\T}{\mathbb T}
\newcommand{\U}{\mathbb U}
\newcommand{\D}{\mathbb D}
\newcommand{\PP}{\mathbb P}
\newcommand{\Sp}{\mathbb S}
\newcommand{\K}{\mathbb K}
\newcommand{\car}{\mathbf 1}
\newcommand{\g}{\mathfrak g}
\newcommand{\gs}{\mathfrak s}
\newcommand{\h}{\mathfrak h}
\newcommand{\rr}{\mathfrak r}
\newcommand{\s}{\sigma}
\newcommand{\fhi}{\varphi}
\newcommand{\ffhi}{\tilde{\varphi}}
\newcommand{\moins}{\setminus}
\newcommand{\ds}{\subset}
\newcommand{\W}{\mathcal W}
\newcommand{\WW}{\widetilde{W}}
\newcommand{\F}{\mathcal F}
\newcommand{\G}{\mathcal G}
\newcommand{\CC}{\mathcal C}
\newcommand{\RR}{\mathcal R}
\newcommand{\DD}{\mathcal D}
\newcommand{\M}{\mathcal M}
\newcommand{\B}{\mathcal B}
\newcommand{\cS}{\mathcal S}
\newcommand{\HH}{\mathcal H}
\newcommand{\Hyp}{\mathbb H}
\newcommand{\UU}{\mathcal U}
\newcommand{\Pp}{\mathcal P}
\newcommand{\QQ}{\mathcal Q}
\newcommand{\E}{\mathcal E}
\newcommand{\GG}{\Gamma}
\newcommand{\LL}{\mathcal L}
\newcommand{\KK}{\mathcal K}
\newcommand{\TT}{\mathcal T}
\newcommand{\X}{\mathcal X}
\newcommand{\Y}{\mathcal Y}
\newcommand{\ZZ}{\mathcal Z}
\newcommand{\tcW}{\widetilde{\mathcal W}}
\newcommand{\bE}{\widehat{E}}
\newcommand{\bF}{\widehat{F}}
\newcommand{\wF}{\widetilde{F}}
\newcommand{\hcF}{\widehat{\mathcal F}}
\newcommand{\bW}{\widehat{W}}
\newcommand{\bcW}{\widehat{\mathcal W}}
\newcommand{\tL}{\widetilde{L}}
\newcommand{\hB}{\widehat{B}}
\newcommand{\hN}{\widehat{N}}
\newcommand{\hM}{\widehat{M}}
\newcommand{\hL}{\widehat{L}}
\newcommand{\diam}{{\rm diam}}
\newcommand{\diag}{{\rm diag}}
\newcommand{\Jac}{{\rm Jac}}
\newcommand{\Plong}{{\rm Plong}}
\newcommand{\Tr}{{\rm Tr}}
\newcommand{\Conv}{{\rm Conv}}
\newcommand{\Ext}{{\rm Ext}}
\newcommand{\Spec}{{\rm Sp}}
\newcommand{\Isom}{{\rm Isom}\,}
\newcommand{\Supp}{{\rm Supp}\,}
\newcommand{\Grass}{{\rm Grass}}
\newcommand{\Hold}{{\rm H\ddot{o}ld}}
\newcommand{\Ad}{{\rm Ad}}
\newcommand{\ad}{{\rm ad}}
\newcommand{\Aut}{{\rm Aut}}
\newcommand{\End}{{\rm End}}
\newcommand{\Leb}{{\rm Leb}}
\newcommand{\Int}{{\rm Int}}
\newcommand{\cc}{{\rm cc}}
\newcommand{\grad}{{\rm grad}}
\newcommand{\proj}{{\rm proj}}
\newcommand{\mass}{{\rm mass}}
\newcommand{\dive}{{\rm div}}
\newcommand{\dist}{{\rm dist}}
\newcommand{\PSL}{{\rm PSL}}
\newcommand{\im}{{\rm Im}}
\newcommand{\re}{{\rm Re}}
\newcommand{\codim}{{\rm codim}}
\newcommand{\Map}{\longmapsto}
\newcommand{\vide}{\emptyset}
\newcommand{\tr}{\pitchfork}
\newcommand{\ssl}{\mathfrak{sl}}

\newenvironment{demo}{\noindent{\textbf{Proof.}}}{\quad \hfill $\square$}
\newenvironment{pdemo}{\noindent{\textbf{Proof of the proposition.}}}{\quad \hfill $\square$}
\newenvironment{IDdemo}{\noindent{\textbf{Idea of proof.}}}{\quad \hfill $\square$}

\def\to{\mathop{\rightarrow}}
\def\act{\mathop{\curvearrowright}}
\def\To{\mathop{\longrightarrow}}
\def\Sup{\mathop{\rm Sup}}
\def\Max{\mathop{\rm Max}}
\def\Inf{\mathop{\rm Inf}}
\def\Min{\mathop{\rm Min}}
\def\lims{\mathop{\overline{\rm lim}}}
\def\limi{\mathop{\underline{\rm lim}}}
\def\egal{\mathop{=}}
\def\dans{\mathop{\subset}}
\def\surj{\mathop{\twoheadrightarrow}}

\def\h#1#2#3{ {\mathfrak{h}}^{#1}_{#2\rightarrow#3}}

\maketitle

\begin{abstract}
In this paper we develop a notion of Gibbs measure for the geodesic flow tangent to a foliated bundle over a compact negatively curved base. We also develop a notion of $F$-harmonic measure and prove that there exists a natural bijective correspondence between these two concepts.

For projective foliated bundles with $\C\PP^1$-fibers without transverse invariant measure, we show the uniqueness of these measures for any Hölder potential on the base. In that case we also prove that $F$-harmonic measures are realized as weighted limits of large balls tangent to the leaves and that their conditional measures on the fibers are limits of weighted averages on the orbits of the holonomy group.
\end{abstract}

\section{Introduction}

Ergodic theory of flows and diffeomorphisms is the statistical study of the long term behaviour of their trajectories. By Birkhoff's theorem it is the study of their invariant measures. In the present work, we will study the statistical, or ergodic, behaviour of the leaves of a foliation. The first difficulty we face is that invariant measures for foliations exist in fact very rarely. In order to overcome this difficulty, we have to search for natural candidates describing where a leaf spends most of its time. 

An important step in the development of the theory was Garnett's introduction of \emph{harmonic measures} in 1983 (see \cite{Gar}). These measures describe the behaviour of Brownian paths tangent to the leaves of a foliation in the sense that they are invariant by a \emph{leafwise heat diffusion} operator. But leaves of a foliation of dimension $\geq 2$ may have a very complicated topology. Thus there could be other natural ways to visit the leaves than following Brownian paths. We will restrict our discussion to the context where the leaves are negatively curved and transverse to a fiber bundle. Before stating our goals and results, let us introduce the main characters of the present paper.

\subsection{Foliated bundles with negatively curved base and the foliated geodesic flow}

\paragraph{Suspension.} Let $B$ and $X$ be two closed (i.e. $C^{\infty}$, compact and boundaryless) manifolds. Consider a representation $\rho:\pi_1(B)\to\dif(X)$. Then the group $\pi_1(B)$ acts on the universal cover $N=\widetilde{B}$ by deck transformations as well as on $X$ by $\rho$. Hence it acts diagonally on the product $N\times X$.

It is proven in \cite{CL} that the quotient $M$ is a smooth manifold, which is endowed with a structure of \emph{foliated bundle}. The partition $(\{z\}\times X)_{z\in N}$ descends to a fiber bundle $\Pi:M\to B$ with base $B$ and fiber $X$, and the partition $(N\times\{x\})_{x\in X}$ descends to a foliation $\F$ whose leaves are covering spaces of the base and are transverse to the fibers. This foliation is called the \emph{suspension} of the representation $\rho$. The latter is called the \emph{holonomy representation} of the foliated bundle. The image group $\rho(\pi_1(B))$ is called the \emph{holonomy group}.

The foliated bundle is determined by the data $(\Pi,M,B,X,\F)$. We can prove that two foliated bundles with conjugated holonomy representations are leaf equivalent: they are conjugated by a diffeomorphism of the ambient spaces which sends leaf to leaf. Every foliation transverse to a fiber bundle $\Pi:M\to B$ is conjugated to a suspension.

In particular, the dynamics of the foliation is entirely determined by $\rho$. For example, one can prove (as in \cite[Lemma 1.4]{BM}) that the existence of a family of transverse holonomy invariant measures for $\F$ is equivalent to the existence of a measure in $X$ which is $\rho(\pi_1(B))$-invariant.

When $X=\C\PP^1$ and $\rho$ takes values in the group $\PSL_2(\C)$ of projective transformations, we say that the foliated bundle is \emph{projective}.

\paragraph{Parametrization.} If $(\Pi,M,B,X,\F)$ is a foliated bundle, we can parametrize the leaves of $\F$ by a Riemannian structure on $B$. When $g$ is a Riemannian metric on $B$, we lift it to the leaves via the fibration, so as to get a leafwise Riemannian metric on $M$. When restricted to the leaves, $\Pi$ becomes a Riemannian submersion: locally it is an isometry. We say that the bundle is \emph{parametrized by the Riemannian structure of} $B$.

\paragraph{The unit tangent bundle.} The \emph{unit tangent bundle} of $B$, i.e. the set of unit vectors tangent to $B$, is denoted by $T^1B$. The unit tangent bundle of $\F$, i.e. the set of unit vectors tangent to leaves of $\F$, is denoted by $T^1\F$.

The differential of $\Pi$ induces a fiber bundle $\Pi_{\ast}:T^1\F\to T^1B$. There is a foliation $\widehat{\F}$ transverse to the fibers of $\Pi_\ast$ whose leaves are precisely the unit tangent bundles $T^1L$ of the leaves $L$ of $\F$. It is easy to see that the holonomy representation $\widehat{\rho}:\pi_1(T^1B)\to\dif(X)$ factorizes as $\rho:\pi_1(B)\to\dif(X)$: there is no holonomy along the sphere factor. 

In particular, there is a transverse measure invariant by the holonomy of $\F$ if and only if there is one for that of $\widehat{\F}$.

\paragraph{The foliated geodesic flow.} The \emph{geodesic flow} of $T^1B$ is defined by flowing unit vectors along the geodesic they direct at unit speed. For every leaf $T^1L$ of $\hcF$, $\Pi_{\ast|T^1L}:T^1L\to T^1B$ is a local isometry. Hence the lift via $\Pi_{\ast}$ of $g_t$ to the leaves of $\hcF$ is the \emph{foliated geodesic flow} $G_t:T^1\F\to T^1\F$ i.e. the leaf-preserving flow of $T^1\F$ which induces on every leaf $T^1L$ its geodesic flow.

We emphasize that we will see this flow as a flow of the \emph{compact} manifold $T^1\F$ on which the leaves $T^1L$ are typically immersed in a complicated way. We are not interested in the intrinsic dynamics \emph{inside $T^1L$}, which is a priori non-compact. This will be our principal tool to understand the dynamical behaviour of the foliation. A situation where this tool appears to be useful is when leaves of $\F$ are parametrized by a negatively curved metric on the base. In that case the foliated geodesic flow must enjoy some hyperbolic properties.

\paragraph{Foliated hyperbolicity.} Suppose the metric $g$ on $B$ is negatively curved. In that case it is known that the geodesic flow is  \emph{Anosov}. This means that it preserves two continuous foliations of $T^1B$ the first one being exponentially contracted, the other one, exponentially expanded. These foliations are respectively called the \emph{stable} and \emph{unstable} foliations and denoted by $\W^s$ and $\W^u$.

These foliations can be lifted to the leaves of $\hcF$ via $\Pi_{\ast}$, which is a local isometry in restriction to every $T^1L$. As a consequence we obtain two continuous and $G_t$-invariant subfoliations of $\hcF$, denoted by $\bcW^s$ and $\bcW^u$. The first one, the \emph{stable one} is uniformly contracted, and the second one, the \emph{unstable one}, unformly expanded.

This weak form of hyperbolicity has been formalized by Bonatti, G\'omez-Mont and Mart\'inez in a recent preprint \cite{BGM}. They called it \emph{foliated hyperbolicity}. It has the flavour of partial hyperbolicity. Here, there is uniform hyperbolicity \emph{inside the leaves}, and there is a direction transverse to $\hcF$, similar to the central direction of partially hyperbolic flows. We don't know much about the way $G_t$ acts on this transverse direction: this is prescribed by the holonomy of the underlying foliation. Actually this is precisely understanding this action that interests us. In particular the major difference between these two weak forms of hyperbolicity is that we don't know a priori that contraction and expansion in the transverse direction are \emph{dominated} by the contraction in the stable direction, and the expansion in the unstable one.

\subsection{Unique ergodicities}
We now describe three natural ways to visit the leaves of a foliation.

\paragraph{Harmonic measures.} As mentionned earlier, we can study the distribution of Brownian paths along the leaves. Measures describing this distribution are Garnett's \emph{harmonic measures}. One possible way to define them is by their disintegration in local charts of $\F$. In a foliated chart $P\times T$ of $\F$ a harmonic measure $m$ reads as
$$dm_{|P\times T}=\left(h(x,t)d\,\Leb_{P\times\{t\}}\right)\,d\nu(t),$$
where $\nu$ is a finite Borel measure on $T$, $\Leb_{P\times \{t\}}$ denotes the Lebesgue measure of the plaque $P\times\{t\}$ and $h:P\times T\to [0,\infty)$ 
 is a measurable function such that for $\nu$-almost every $t\in T$ $h(.,t)$ is positive, smooth and harmonic for the Laplace operator on $P\times\{t\}$. Garnett proved in \cite{Gar} the existence of such measures for foliations of compact manifolds.

\paragraph{SRB measures.} It is also natural to wonder where are going ``almost all'' geodesics in the leaves, in the sense of Lebesgue. A way to formalize this is by looking for measures in the unit tangent bundle of the foliation $T^1\F$ which are \emph{physical}, or \emph{SRB measures} for the foliated geodesic flow. Such an object is a $G_t$-invariant measure $\mu$ on $T^1\F$ whose \emph{basin} (the set of $v\in T^1\F$ such that the averages of the Dirac masses along the orbit of $v$ converges to $\mu$ in the weak$^{\ast}$ sense) has positive volume. These measures are named after by Sinai, Ruelle and Bowen who introduced them for uniformly hyperbolic dynamics in \cite{BR,Si}.

\paragraph{Large balls.} Finally we can ask what is the distribution of large balls tangent to the leaves of $\F$. More precisely,  let $N$ be the universal cover of a leaf $L$. Consider the \emph{uniform distribution} (i.e. the normalized Lebesgue measure) on a large ball in $N$ and project it down somewhere onto the leaf. These are multidimensional analogues of Birkhoff averages (we refer to \S \ref{slargeballs} for more details). What can be said about accumulation points of such measures?

\paragraph{Triple unique ergodicity.} Let us now state a consequence of the main results of this paper.

\begin{maintheorem}
\label{tripleuniqueergodicity}
Let $(\Pi,M,B,\C\PP^1,\F)$ be a projective foliated bundle parametrized by a closed and negatively curved Riemannian manifold $B$. Assume that the holonomy group has no invariant probability measure on $\C\PP^1$. Then the following assertions hold true:
\begin{itemize}
\item there exists a unique harmonic measure on $M$ for $\F$;
\item there exists a unique SRB measure for the foliated geodesic flow $G_t$ on $T^1\F$;
\item there exists a unique measure on $M$ which is accumulation point of normalized volumes of large balls tangent to the leaves.
\end{itemize}

Moreover, suppose $B$ is a closed surface with variable negative curvature and the holonomy representation $\rho:\pi_1(B)\to \PSL_2(\C)$ is quasifuchsian. Then the harmonic measure, the projection on $M$ of the SRB measure for $G_t$ and the accumulation point of normalized area of large discs are pairwise mutually singular.
\end{maintheorem}

In the early 2000's, in a series of three papers \cite{BG,BGV,BGVil}, Bonatti and G\'omez-Mont, together with Vila and Viana, showed the dichotomy when the base is a \emph{hyperbolic surface}. It turns out that in this case the three measures coincide. This fact is a bit magical and one of our goals is to shed some light about the case where the base is of variable curvature, and of arbitrary dimension.

Note that the uniqueness of the harmonic measure in this case has also been obtained by the author in \cite{Al1} by completely different methods. Note also that in \cite{DK}, Deroin and Kleptsyn proved it in the much more general context of transversally conformal foliations.

The uniquenes of the SRB measures has also been proven in a more general context: that of transversally conformal foliations with negatively curved leaves without holonomy invariant measure. First Bonatti, G\'omez-Mont and Mart\'inez treated the case where the leaves are hyperbolic in \cite{BGM}. They used Deroin-Kleptsyn's result as well as the bijective correspondence between harmonic measures and some special $G_t$-invariant measures on the unit tangent bundle, that we will be led to introduce later on, proven in \cite{Al2,BMar,Ma}. Then together with Yang, we proved in \cite{AY} the general case using Pesin's theory as well as a criterion for the existence of holonomy-invariant measures proved in \cite{Al3}.

The measures obtained in Theorem \ref{tripleuniqueergodicity} may be singular: it reveals that the behaviour of the leaves at infinity has different meanings according to the chosen point of view. It would be nice to prove a rigidity result, since the example we give is too simple.
\begin{ques}
\label{rigidity}
For a projective foliated bundle over a negatively curved closed surface $\Sigma$ with fiber $\C\PP^1$, and without holonomy invariant measure, is it true that if two of the three measures obtained in Theorem \ref{tripleuniqueergodicity} (harmonic, projection of the SRB measures or accumulation of large discs tangent to the leaves) are equivalent, the curvature is constant?
\end{ques}

\subsection{Dictionary between Gibbs and harmonic measures}
A priori the three measures obtained in Theorem \ref{tripleuniqueergodicity} are of different nature and live in different spaces. Harmonic measures live on $M$ and are defined by specifying the density of conditional measures in the plaques of $\F$. SRB measures live on $T^1\F$ and will be obtained as a $G_t$-invariant measures satisfying a special property that we will call a \emph{Gibbs property}. Finally, limits of large balls live on $M$  are obtained by an equidistribution process.

The main goal of this work is to provide a dictionary which allows to unify these points of view and to give a unified proof of our theorem. The key lies in the theory of Gibbs states, and in the theory of Hölder cocycles developped by Ledrappier \cite{L3} and carried out especially by Roblin \cite{Rob} and Paulin-Pollicott-Schapira \cite{PPS}.

We will associate two objects to every Hölder continuous potential $F:T^1B\to\R$. One is a $G_t$-invariant measure on the unit tangent bundle $T^1\F$ with a special property called \emph{Gibbs property}. The other one is a measure on the manifold $M$ having a local description analogous to the one of harmonic measures.

\paragraph{Gibbs measures.} The theory of Gibbs measures is well developped in uniformly hyperbolic dynamics because these dynamics have nice symbolic representations (see \S \ref{ssmarkov} and \ref{sgibbsstates}). We refer to \cite{BL,Bo1,Bo2,Ra,Si} for the classical theory, and to \cite{PPS} for the many characterizations and features of Gibbs states, as well as for a theory of Gibbs measures of the geodesic flow on noncompact manifolds with negative curvature.

Although there exist lots of analogies with these works, we deal with a weaker form of hyperbolicity, which resembles partial hyperbolicity. For such systems Pesin and Sinai defined the Gibbs property by specifying the class of the disintegration of invariant measures in unstable manifolds: see \cite{PS}. If one requires that the class is Lebesgue, we get the definition of a \emph{Gibbs u-state}. This can be easily adapted to our context. Gibbs $u$-states have succesfully been used to construct SRB measures for some partially hyperbolic diffeomorphisms in \cite{ABV,BV,Car,Y} and for some Hénon maps in \cite{BY}. We refer to \cite[Chapter 11]{BDV} for an account of the theory.

Let us explain the strategy. Associated to the potential $F:T^1B\to\R$ there is a unique Gibbs state $\mu_F$ for the geodesic flow $g_t:T^1B\to T^1B$ characterized by a local product structure (see \S\ref{lpsgibbsstates}). This provides a well defined \emph{measure class} in the leaves of $\W^u$ which is invariant by the geodesic flow (see Theorem \ref{familiesofmeasures} for details). We can lift it to the leaves of $\bcW^u$ via $\Pi_{\ast}$. This will allow us to define a \emph{Gibbs measure} for $G_t$ associated to the lifted potential $\bF=F\circ\Pi_{\ast}:T^1\F\to\R$ as a $G_t$-invariant probability measure whose conditional measures in plaques of $\bcW^{u}$ are in the measure class we just defined.

\paragraph{$F$-harmonic measures.} Let $N$ denote the universal cover of $B$ and $N(\infty)$ the sphere at infinity (see \S \ref{sgeometry}). It is proven in \cite{AS} that a positive harmonic function $h:N\to(0,\infty)$ has an integral representation
\begin{equation}
\label{Poissonrep}
h(z)=\int_{N(\infty)} k(o,z;\xi)d\eta(\xi),
\end{equation}
where $o\in N$ is some base point, $k:N\times N\times N(\infty)\to[0,\infty)$ is the \emph{Poisson kernel} and $\eta$ is a finite Radon measure on $N(\infty)$. This simple fact enabled us in \cite{Al2} to ``unroll'' harmonic measures and to lift them to the unit tangent bundle when the foliation has negatively curved leaves (see \S \ref{unrolling} for an outline of the argument).

Our strategy is to use, instead of the Poisson kernel, a kernel provided by Ledrappier's theory of Hölder cocycles (see \cite{L3}) that we shall define and call \emph{Gibbs kernel} in \S \ref{sgibbsstatesgeo}. We can define $F$-harmonic functions as the ones with an integral representation as in \eqref{Poissonrep} but with Gibbs kernel instead of Poisson kernel. Finally, we can define $F$-harmonic measures for $\F$ as we define harmonic measures. Locally, we require the conditional measures in the plaques to have an $F$-harmonic density with respect to Lebesgue. Note that the existence of these measures is far from being obvious, because a priori they are not invariant by an operator semigroup.

\paragraph{Special cases.} The three measures obtained in Theorem \ref{tripleuniqueergodicity}, as we shall see, all come from some potentials.

\begin{itemize}
\item \emph{The harmonic measure} corresponds to the potential
$$H(v)=d/dt|_{t=0}\log k(c_v(0),c_v(t);c_v(-\infty)),$$
where $k$ is the Poisson kernel and $c_v$ denotes the geodesic directed by $v$.
\item \emph{The SRB measure} corresponds to the potential
$$\phi^u(v)=-d/dt|_{t=0}\log\,\Jac^u g_t(v),$$
where $\Jac^u=\det Dg_t{|E^u}$ denotes the unstable jacobian of the geodesic flow.
\item \emph{The limit of large discs} corresponds to constant potentials.
\end{itemize}

The uniqueness part in Theorem \ref{tripleuniqueergodicity} is then a consequence of the following theorem: the precise definitions will be given in Section \ref{sgibbsharmonicbundles}.

\begin{maintheorem}
\label{tripleuniquegibbs}
Let $(\Pi,M,B,X,\F)$ be a foliated bundle parametrized by a closed and negatively curved Riemannian manifold $B$ and whose fiber $X$ is a closed manifold. Let $F:T^1B\to\R$ be a Hölder continuous potential and $\bF:T^1\F\to\R$ be its lift.
\begin{enumerate}
\item There is a canonical bijective correspondence between $F$-harmonic measures and Gibbs measures associated to $\bF$.
\item There exist Gibbs measures associated to $\bF$.
\item If $X=\C\PP^1$ and if $\rho(\pi_1(B))\leq\PSL_2(\C)$ does not preserve any measure on $\C\PP^1$, then the Gibbs measure associated to $\bF:T^1\F\to\R$, and the $F$-harmonic measure for $\F$ are unique.
\end{enumerate}

\end{maintheorem}

\subsection{Overview of the paper} In Section \ref{sgeneralities}, we give some tools that will be useful in the whole paper. In Section \ref{sgibbsharmonicbundles}, we define Gibbs measures for the foliated geodesic flow as well as the notion of $F$-harmonic measures for foliated bundles with negatively curved leaves. There we describe a canonical bijective correspondence between the two objects. In Section \ref{suniqueness} we discuss the uniqueness of these objects in the context of $\C\PP^1$-bundles and we prove in particular Theorem \ref{tripleuniquegibbs}. We also prove Theorem \ref{disintegrationfharmonic} which gives the desciption of conditional measures of the unique $F$-harmonic measure in the fibers of $\Pi:M\to B$. In Sections \ref{slargeballs} and \ref{sequidistribution}, we give some properties of $F$-harmonic measures in the context of foliated $\C\PP^1$-bundles with no transverse invariant measures. We prove, Theorem \ref{limitlargeballs}, that they are obtained as limits of averaging sequences and, Theorem \ref{countingmeasures}, that their conditional measures on the fibers are limits of weighted counting measures. In Section \ref{sfuchsianquasifuchsian}, we treat the examples of fuchsian and quasifuchsian representations and we finish the proof of Theorem \ref{tripleuniqueergodicity}.

\section{Some generalities}
\label{sgeneralities}

\subsection{Holonomy of foliated bundles}
\label{holonomyfolbdle}

Let $(\Pi,M,B,X,\F)$ be a foliated bundle where $B$ and $X$ are closed manifolds. Let $(V_i,\phi_i)_{i\in I}$ be an atlas in $B$ consisting of trivializing charts. By this we mean that for every $i\in I$, $\phi_i:U_i\to V_i\times X$ is a fiber preserving diffeomorphism (where $U_i=\Pi^{-1}(V_i)$), such that $\phi_i(\F_{|U_i})$ is the partition $(V_i\times\{x\})_{x\in X}$.

\paragraph{Composition of transition maps.} Assume that $(V_i,\phi_i)_{i\in I}$ is a \emph{good atlas} in the sense that when $V_i\cap V_j\neq\vide$, then the intersection is connected and the closure of $V_i\cup V_j$ is contained in a trivializing chart. \emph{Transition maps} are the diffeomorphisms $\tau_{ij}:X\to X$ defined when $V_i\cap V_j\neq\vide$ in such a way that $\phi_j^{-1}\circ\phi_i=(Id,\tau_{ij})$ in restriction to $(V_i\cap V_j)\times X$.

Let $c$ be a path on $B$ and $(V_{i_0},...,V_{i_n})$ be a chain of charts covering $c$. The \emph{holonomy over $c$} is by definition the diffeomorphism of $X$ defined as $\tau_c=\tau_{i_{n-1}i_n}\circ...\circ\tau_{i_0i_1}$. It does not depend on the covering nor on the path $c$ but only on the homotopy class $[c]$.

\paragraph{Holonomy via lifting.} We give another, more geometric yet equivalent, point of view. We will make use of both of them. If $L$ is a leaf of $\F$, $\Pi_{|L}:L\to B$ is a covering map (see \cite{CL}). Hence by choosing $p\in B$, $c$ a path of $B$ starting at $p$ and $x\in X_p$ there exists a unique lift of $c$ to $L$ starting at $x$. The ending point of this path is $\tau_c(x)\in X_q$, $q$ denoting the ending point of $c$.

In particular if $c$ is a loop based at $p$ and $\gamma\in\pi_1(B)$ is its homotopy class then $\tau_c=\rho(\gamma)^{-1}$.

\subsection{Disintegration and foliations}
\label{sdisintegration}

\paragraph{Absolutely continuous/singular disintegration on the leaves of a foliation.} A priori, disintegration on the leaves of a foliation is impossible. For example, the linear foliation of a torus by lines of irrational slope is not a measurable partition in the sense of Rokhlin \cite{Ro}.

What is possible to do is to disintegrate a finite measure on the \emph{plaques} of any foliated chart and to compare these measures (which a priori depend on the choice of the chart) with a family of measures on the leaves. The memoir \cite{KL} is a good introduction for the notion of tangential and transverse measures for foliations.

\begin{defi}
\label{disintegrationleaves}
Let $(\F,\lambda^{\F}_x)$ be a foliation of a compact manifold $M$, together with a family of Radon measures on the leaves satisfying the two following assumptions. Firstly, for $x,y$ in the same leaf, $\lambda^{\F}_x=\lambda^{\F}_y$. Secondly, $x\mapsto\lambda^{\F}_x$ varies measurably transversally in local charts of $\F$.

We will say that a finite Radon measiure $\mu$ has an absolutely continuous (resp. singular) disintegration with respect to $(\F,\lambda^{\F}_x)$ if its conditional measure in almost every plaque  
$P_x$ is absolutely continuous (resp. singular) with respect to $\lambda^{\F}_x$.

In the special case where $\lambda^{\F}_x$ is the Lebesgue measure, we will say that $\mu$ has Lebesgue disintegration.
\end{defi}

\subsection{Negatively curved manifolds and the geodesic flow}
\label{sgeometry}
In this paragraph, $N$ denotes the universal cover of a closed Riemannian manifold $B$ with negative sectional curvature. Its sectional curvature is pinched between two negative constants $-b^2\leq-a^2<0$.

\paragraph{Sphere at infinity.} We refer to \cite{AS} for the proofs of the assertions below. Let $N(\infty)$ denote the \emph{sphere at infinity} of $N$, i.e. the set of equivalence classes of geodesic rays for the relation ``stay at bounded distance''.

We consider $\pi_o:T^1_oN\to N(\infty)$ the natural projection which associates to $v$ the class of the geodesic ray it determines. All the maps $\pi_{o'}^{-1}\circ\pi_o:T^1_oN\to T^1_{o'}N$ are Hölder continuous: these projections determine a natural Hölder structure on the sphere at infinity.

The union $N\cup N(\infty)$ is endowed with the \emph{cone topology}: the neighbourhoods of infinity are given by the \emph{truncated cones} $T_o(v,\theta,R)=C_o(v,\theta)\moins B(o,R)$, where:
$$C_o(v,\theta)=\{x\in N|\angle_o(v,v_{ox})\leq\theta\},$$
with $v\in T^1_oN$, $\theta\in\R$ and $v_{ox}$, the unit vector tangent to $o$ which points to $x$. Moreover, one says that a sequence $(x_i)_{i\in\N}$ \emph{converges nontangentially} to a point at infinity $\xi$ if it converges in the cone topology while staying at bounded distance from a geodesic ray.

\paragraph{Busemann cocycle and horospheres.}We define the \emph{Busemann cocycle} as
$$\beta_{\xi}(y,z)=\lim_{t\to\infty}\dist(c(t),z)-\dist(c(t),y),$$
where $\xi\in N(\infty)$, $y,z\in N$ and $c$ is any geodesic ray parametrized by arc length and pointing to $\xi$. It satisfies the following cocycle relation
\begin{equation}
\label{Eq:cocyclebusemann}
\beta_{\xi}(x,y)+\beta_{\xi}(y,z)=\beta_{\xi}(x,z).
\end{equation}

\paragraph{Remark.} The sign convention is different from those of \cite{K1,PPS}. It is used so that the classical formula $k(y,z;\xi)=e^{-(d-1)\beta_{\xi}(y,z)}$ holds in the hyperbolic space $\mathbb{H}^d$ where $k$ is the Poisson kernel (see \cite{AS,L1,L2}).\\

The \emph{horospheres} are the level sets of this cocycle: two points $y,z$ are said to be on the same horosphere centered at $\xi$ if $\beta_{\xi}(y,z)=0$. It is possible (see next paragraph) to see that horospheres are $C^{\infty}$ manifolds.

\paragraph{The geodesic flow.} 
A vector $v\in T^1B$ directs a unique geodesic. The \emph{geodesic flow} is defined by flowing $v$ along this geodesic at unit speed. We denote it by $g_t:T^1B\to  T^1B$. It is well known that this flow is an \emph{Anosov flow} \cite{An}: there is a (Hölder) continuous splitting $TT^1B=E^s\oplus E^u\oplus\R X$, where $X$ is the generator of the flow and $E^s$, $E^u$ are $D g_t$-invariant subbundles which are respectively uniformly contracted and dilated by the flow

$$\begin{cases}\displaystyle{\frac{a}{b}e^{-bt}||v_s||\leq ||Dg_t(v_s)||\leq \frac{b}{a}e^{-at}||v_s||} & \text{if $v_s\in E^s$ and $t>0$} \\
                                                                                                                      \\
               \displaystyle{\frac{a}{b}e^{-bt}||v_u||\leq ||Dg_{-t}(v_u)||\leq \frac{b}{a}e^{-at}||v_u||} & \text{if $v_u\in E^u$ and $t>0$}
   .\end{cases}$$
The bundles $E^s$ and $E^u$ are respectively called \emph{stable} and \emph{unstable} bundles. They are uniquely integrable: we denote by $\W^s$, $\W^u$ the integral foliations called respectively \emph{stable} and \emph{unstable} foliations. Moreover, by the stable manifold theorem, the stable and unstable manifolds are as smooth as the geodesic flow that is $C^{\infty}$.

If we lift this flow to $T^1N$ via the differential of the universal covering map, we obtain the geodesic flow of $T^1N$, denoted by $G_t:T^1N\to T^1N$. The lifts of invariant foliations shall be denoted by $\tcW^{\star}$ and their leaves $\widetilde{W}^{\star}(v)$, $\star=s,u,cs,cu$. Stable (resp. unstable) leaves are identified with stable (resp. unstable) horospheres, i.e. horospheres endowed with the normal inward (resp. outward) vector field. This proves that horospheres are smooth, as are stable and unstable manifolds.

We have the identification $T^1N\simeq N\times N(\infty)$ obtained by sending $v$ on the couple $(c_v(0),c_v(-\infty))$ where $c_v$ is the directed geodesic determined by $v$. This identification conjugates the actions of the group of direct isometries $\Isom^+(N)$ on $T^1N$ by differentials and on $N\times N(\infty)$ by diagonal maps. Moreover, it also trivializes the unstable foliation: a slice $N\times\{\xi\}$ has to be thought as filled with unstable horospheres centered at $\xi$.

\paragraph{Notation.} In the particular case of holonomy maps along (center) stable and (center) unstable leaves, we will adopt the following notation
$$\h {\star}{T_1}{T_2}$$
for $\star=s,u,cs,cu$ and $T_1$, $T_2$ two small transversals to $\W^{\star}$. This map is defined on a small relatively compact open set $S_1\dans T_1$. When $T_1$ and $T_2$ are pieces of small center stable (in the case $\star=u$), or center unstable (in the case $\star=s$) manifolds, we will just denote the corresponding holonomy map
$$\h {\star}{v_1}{v_2}.$$
where $v_1$ and $v_2$ belong to the same leaf of $\W^{\star}$.

\subsection{Markov partitions}
\label{ssmarkov}
\paragraph{Rectangles.}
We will need to consider Markov partitions for $g_t$. We will follow the steps of Bowen and Ratner who managed to prove their existence in the seventies (see \cite{Bo1} and \cite{Ra}).

For any small $\eps$, the local invariant manifold of size $\eps$ of a point $v$, denoted by $W_{\eps}^{\star}(v)$ for $\star=s,u,cs$ or $cu$, is defined as the connected component of $B(v,\eps)\cap W^{\star}(v)$ which contains $v$, where $B(v,\eps)$ stands for the Riemannian ball centered at $v$ and of radius $\eps$.

\begin{proposition}
There exist numbers $0<\delta<\gamma<2\delta$ such that if $v,w\in T^1B$ are at distance at most $\gamma$, then the local unstable manifold of $v$ and the local center stable manifold of $w$ of size $\delta$ intersect transversally in a unique point denoted by $[v,w]$
$$W_{\delta}^{cs}(v)\tr W_{\delta}^{u}(w)=\{[v,w]\}.$$
\end{proposition}

If $X\dans W^u_{\delta}(v)$ and $Y\dans W^{cs}_{\delta}(v)$, $[X,Y]$ will  denote the image of $X\times Y$ by the bracket $[.\,,\,.]$.

\begin{defi}
\label{rectangles}
A rectangle $R$ of $T^1B$ is a set $R=[A^u,A^{s}]$ where $A^{\star}\dans W^{\star}_{\delta}(v)$, for some point $v$ and $\star=u$ or $s$, has a nonempty interior and is the adherence of its interior.
\end{defi}

\paragraph{Remark.} A rectangle $R$ is \emph{topologically transverse} to the flow in the sense that there exists an $\alpha>0$ such that $R\cap g_{(0,\alpha]}(R)=\vide$.

\begin{defi}
A finite family of rectangles  $\RR=(R_i)_{i\in I}$ will be called proper of size $\alpha$ if the union of the $R_i$ intersects each positive orbit and if for any $i\neq j$, $g_{[0,\alpha]}(R_i)\cap R_j=\vide$.
\end{defi}

We introduce the following notation. If $R_i=[A_i^u,A_i^{s}]$ is a rectangle with $A_i^{\star}\dans W^{\star}_{\delta}(v)$ ($\star=u$ or $s$) and $w_1,w_2\in R_i$, we denote $A^u_i(w_1)=[A^u_i,w_1]$ and $A^s_i(w_2)=[w_2,A_i^s]$. Note that we always have $A^u_i(w_1)\dans W^u(w_1)$, whereas $A^s_i(w_2)$ is not included in $W^s(w_2)$ since $\W^s$ and $\W^u$ are not jointly integrable.

For a proper family of rectangles $\RR=(R_i)_{i\in I}$, $R_i=[A^u_i,A^s_i]$, one can (by compactness of $T^1B$) associate to $v\in T^1B$ the least $t>0$ for which $g_t(v)\in\bigcup R_i$. It defines the \emph{first return time} $r(v)$, as well as the \emph{first return map} $T(v)$. These functions are continuous on
\begin{equation}
\label{Eq:residual}
\CC=\{v\in T^1B|T^k(v)\in\bigcup \Int R_i\,\,\,\textrm{for\,\,all}\,\,k\in\Z\},
\end{equation}
which is a residual set and hence is dense by Baire's theorem (see \cite{Bo1}).

Associated to such a proper family of rectangles, there is always a decomposition of $T^1B$ as a union of \emph{cubes}
\begin{equation}
\label{cube}
C_i=\bigcup_{v\in R_i}g_{[0,r(v)]}(v).
\end{equation}
Moreover, we can consider the \emph{open cubes}
\begin{equation}
\label{Eq:opencube}
C_i^{\ast}=\bigcup_{v\in \Int R_i}g_{[0,r(v)]}(v).
\end{equation}
Observe that if $\CC'$ is the union of all the orbits of the $C_i\moins C_i^{\ast}$, we have $\CC=T^1B\moins\CC'$ (these two sets are invariant by the flow).
\paragraph{Markov property.} For a proper family of rectangles $\RR$ and two indices $i$ and $j$ we denote $R_{ij}=\Int\,R_i\cap T^{-1}(\Int\,R_j)$.

\begin{defi}
\label{Markovproperty}
A proper family of rectangles $\RR$ is said to have the Markov property, if for any $i$ and $j$ such that $R_{ij}\neq\vide$ and $v\in R_{ij}$, we have
$$A^s_i(v)\dans \overline{R_{ij}},$$
$$A^u_j(T(v))\dans\overline{T(R_{ij})}.$$
\end{defi}

\begin{defi}
A Markov partition is a family of cubes $(C_i)_{i\in I}$ as defined in \eqref{cube} which is associated to a proper family of rectangles $(R_i)_{i\in I}$ satisfying the Markov property.

The sets $A_i^u(v)$, for $v\in C_i^{\ast}$, $i\in I$ will be called the Markovian unstable plaques.
\end{defi}

Bowen and Ratner showed in \cite{Bo1} and \cite{Ra} the existence of proper families of rectangles with the Markov property and arbitrarily small size. We will need a version of this existence theorem given in \cite{BGV}.

\begin{theorem}[Existence of Markov partitions]
\label{Markovpartition}
Let $\UU=(U_j)_{j\in J}$ be any open cover of $T^1B$. Then there exists an adapted Markov partition, that is a proper family of rectangles $\RR=(R_i)_{i\in I}$ with the Markov property such that there exist two surjective functions $\alpha$, $\beta:I\to J$ satisfying, for every $i\in I$
\begin{enumerate}
\item $C_i\dans U_{\alpha(i)}$;
\item for every $v\in C_i$, $g_{r(v)}(v)\in U_{\beta(i)}$.
\end{enumerate}
Here the cube $C_i$ is defined as in \eqref{cube}. Moreover the set $\CC$ defined by \eqref{Eq:residual} is a dense $G_{\delta}$ of full measure for any ergodic invariant measure for $g_t$, which is positive on all open subsets.
\end{theorem}

\subsection{Gibbs states}
\label{sgibbsstates}

\paragraph{Cocycles on stable and unstable foliations.} Associated to every Hölder potential $F:T^1B\to\R$, there are two cocycles $k_F^s$ and $k_F^u$ defined by

\begin{equation}
\label{Eq:cocs}
k_F^s(v,w)=\exp\left[\int_0^{\infty}(F\circ g_t(w)-F\circ g_t(v)) dt\right]
\end{equation}
for $v,w$ in the same stable leaf and
\begin{equation}
\label{Eq:cocu}
k_F^u(v,w)=\exp\left[\int_0^{\infty}(F\circ g_{-t}(w)-F\circ g_{-t}(v)) dt\right]
\end{equation}
for $v,w$ in the same unstable leaf.

The existence of these cocycles is an immediate consequence of the Hölder continuity of $F$ and of the usual distortion controls.

\paragraph{Local product structure.} The following result can be found in \cite[Lemma 2.1]{BL}.

\begin{theorem}
\label{familiesofmeasures}
Let $B$ be a closed Riemannian manifold with negative sectional curvature and $F:T^1B\to\R$ be a Hölder continuous potential. Then:
\begin{enumerate}

\item there exists a family $(\lambda^{cu}_{F,v})_{v\in T^1B}$ of measures defined on the leaves of $\W^{cu}$ which satisfy $\lambda^{cu}_{F,v}=\lambda^{cu}_{F,w}$ when $w\in W^{cu}(v)$, as well as the cocycle relation
\begin{equation}
\label{Eq:cocycleweakunstable}
\frac{d\left[\h sv{v'}\,_{\ast}\lambda^{cu}_{F,v}\right]}{d\lambda^{cu}_{F,v'}}(w)=k_F^s(w,\h s{v'}v(w)),
\end{equation}
for $v,v'$ in the same stable leaf and $w$ in the domain of the holonomy map $\h s{v'}v$.

Moreover, this family are unique up to multiplication by a positive constant;

\item there exist a number $P(F)$, as well as a family $(\lambda^{u}_{F,v})_{v\in T^1B}$ of measures defined on the leaves of $\W^{u}$ which satisfy $\lambda^u_{F,v}=\lambda^u_{F,w}$ when $w\in W^u(v)$ and which are quasi-invariant by the flow $g_t$ with the cocycle relation for for $v\in T^1B$
\begin{equation}
\label{Eq:cocyclestrongunstable}
\frac{d\left[g_T\,_{\ast}\lambda^{u}_{F,g_{-T}(v)}\right]}{d\lambda^{u}_{F,v}}(v)=\exp\left[\int_0^T (F\circ g_{-t}(v)-P(F))dt\right],
\end{equation}

Moreover, this family is unique up to multiplication by a positive constant and this relation uniquely determines the number $P(F)$;
\item we have with obvious abusive notation $\lambda^{cu}=\lambda^u\times dt$ (see \cite{L4}).  In particular, the families of measures $\lambda^{c\star}_{F,v}$ are quasi invariant by the action of the flow, with the same cocycle relations \eqref{Eq:cocyclestrongunstable}.

\item We have also the corresponding families of measures $(\lambda^{cs}_{F,v})_{v\in T^1B}$ and $(\lambda^s_{F,v})_{v\in T^1B}$ that are uniquely defined up to a multiplicative constant, that are constant respectively on center stable and stable leaves and satisfy

\begin{equation}
\label{Eq:cocycleweakstable}
\frac{d\left[\h uv{v'}\,_{\ast}\lambda^{cs}_{F,v}\right]}{d\lambda^{cs}_{F,v'}}(w)=k_F^u(w,\h u{v'}v(w))
\end{equation}

\begin{equation}
\label{Eq:cocyclestrongstable}
\frac{d\left[g_T\,_{\ast}\lambda^{s}_{F,g_{-T}(v)}\right]}{d\lambda^{s}_{F,v}}(v)=\exp\left[-\int_0^T (F\circ g_{-t}(v)-P(F))dt\right].
\end{equation}

and $\lambda^{cs}=\lambda^s\times dt.$

\end{enumerate}
\end{theorem}

Now we can state a result concerning local product structure of Gibbs states, which can be found in \cite[Lemma 2.1]{BL}.

\begin{theorem}
\label{lpsgibbsstates}
Let $B$ be a closed Riemannian manifold with negative sectional curvature and $F:T^1B\to\R$ be a Hölder continuous potential. Then there exists a unique $g_t$-invariant probability measure $\mu_F$, called the Gibbs state associated to $F$, with the following local product structure.
\begin{enumerate}
\item Let $U=[W^u_{loc}(v_0),W^{cs}_{loc}(v_0)]$. When restricted to $U$, the Gibbs state $\mu_F$ disintegrates as
$$d\mu_F\,_{|U}=\left(k^u_F(v,w)\,d\lambda^u_{F,v}(w)\right)\,d\lambda^{cs}_{F,v_0}(v).$$
\item When restricted to $U$, $\mu_F$ also disintegrates as
$$d\mu_F\,_{|U}=\left(k_F^s(v,w)\,d\lambda^s_{F,v}(w)\right)\,d\lambda^{cu}_{F,v_0}(v).$$
\end{enumerate}
\end{theorem}

The following theorem states that it is possible to reconstruct Gibbs states by forward iteration of a unstable small set with positive $\lambda^u_{F,v}$ measure. It can be found in \cite{Ba} in a much more general context.

\begin{theorem}
\label{reconstructthegibbsstate}
Let $F:T^1B\to\R$ be a Hölder continuous potential. Let $v\in T^1B$ and $D\dans W^u_{loc}(v)$ with positive $\lambda^u_{F,v}$ measure. Then the measure
\begin{equation}
\label{Eq:trueconvergence}
\frac{g_T\,_{\ast}(\lambda^u_{F,v})_{|D}}{\lambda^u_{F,v}(D)}
\end{equation}
converges to $\mu_F$ as $T$ tends to infinity.

\end{theorem}

The following result is a consequence of the uniqueness of Gibbs states (see Theorem \ref{lpsgibbsstates}) and of Theorem \ref{reconstructthegibbsstate} and can be seen as an adaptation of \cite[Corollary 11.14.]{BDV} for Gibbs $u$-states.

\begin{coro}
Let $F:T^1B\to\R$ be a Hölder continuous potential. Then the Gibbs state $\mu_F$ is the only measure on $T^1B$ which:
\begin{itemize}
\item is $g_t$-invariant;
\item has an absolutely continuous disintegration with respect to $(\W^u,\lambda^u_{F,v})$.
\end{itemize}
\end{coro}

\subsection{Gibbs kernel and Ledrappier's measures}
\label{sgibbsstatesgeo}

The Gibbs states for the geodesic flow of $T^1B$ have a nice description in the universal cover. Let $F$ be a Hölder continuous function on $T^1B$ and $\wF$ be its lift to $T^1N$. When $z_1,z_2\in N$, we denote by $\int_{z_1}^{z_2}\wF$ the integral of $\wF$ on the \emph{directed} geodesic going from $z_1$ to $z_2$ and parametrized by arc length.

\paragraph{The Gibbs kernel.} Ledrappier defines in \cite{L3} the following cocycle on $N\times N\times N(\infty)$
\begin{equation}
\label{gibbskernel}
k^F(y,z;\xi)=\exp\left[\int_{\xi}^{z}\wF-\int_{\xi}^{y}\wF\right]\exp[-P(F)\beta_{\xi}(y,z)],
\end{equation}
where $\beta_{\xi}$ is the Busemann cocycle at $\xi$. Here we use an abusive notation: the difference of these two integrals has the following meaning. Consider some geodesic ray $c$ which is asymptotic to $\xi$, then the following limit exists, doesn't depend on the choice of $c$ and we set by definition
$$\int_{\xi}^{z}\wF-\int_{\xi}^{y}\wF=\lim_{T\to\infty}\left(\int_{c(T)}^{z}\wF-\int_{c(T)}^{y}\wF\right).$$

\paragraph{Remark.} This cocycle appears also in \cite{PPS} with a different sign conddition. Our sign condition is coherent with the usual Poisson kernel (see \cite{L1,L2}).\\

This cocycle will play the role of the Poisson kernel, and we call it the \emph{Gibbs kernel}. We clearly have the following cocycle relation
\begin{equation}
\label{Eq:Fcocyclerelation}
k^F(x,y;\xi)k^F(y,z;\xi)=k^F(x,z;\xi)
\end{equation}
when $x,y,z\in N$ and $\xi\in N(\infty)$.

Let us describe the Gibbs kernel $k^F$ in three basic examples.
\begin{itemize}
\item The kernel associated to the null function is $k^0(y,z;\xi)=e^{-h\beta_{\xi}(y,z)}$, where $h$ denotes the topological entropy of the geodesic flow $g_t$.
\item The kernel associated to $\phi^u(v)=-d/dt|_{t=0}\log\,\Jac^u g_t(v)$ is:
$$k^u(y,z;\xi)=\lim_{t\to\infty}\frac{\Jac^u G_{-t-\beta_{\xi}(y,z)}(v_{\xi,z})}{\Jac^u G_{-t}(v_{\xi,y})}.$$
\item The Poisson kernel $k$ defined in \cite{AS} is the kernel $k^H$ associated to 
 $$H(v)=d/dt|_{t=0}\log k(c_v(0),c_v(t);c_v(-\infty)).$$
\end{itemize}
Note that the last two potentials have zero pressure (see \cite{BR} for the first one and \cite{L1,L4} for the second one).

\paragraph{Ledrappier's measures.} The Gibbs kernel is normalized in the sense of Ledrappier \cite{L3}. Using a Patterson-Sullivan argument, Ledrappier was able to show that this cocycle thus realizes as a Radon-Nikodym cocycle for a unique (up to multiples) family of measures on $N(\infty)$ (see also \cite{K1}).

\begin{theorem}[Ledrappier]
\label{ledrappiermeasures}
Let $B$ be a closed Riemannian manifold with negative sectional curvature, whose universal Riemannian cover is denoted by $N$. Then, there exists a family $(\nu^F_z)_{z\in N}$ of finite measures on $N(\infty)$ satisfying:
\begin{enumerate}
\item the equivariance property $\gamma_\ast\nu^F_z=\nu^F_{\gamma z}$ for $\gamma\in\pi_1(B)$ and $z\in N$;
\item the cocycle property:
\begin{equation}
\frac{d\nu_{z}^F}{d\nu_{y}^F}(\xi)=k^F(y,z;\xi).
\end{equation}
\end{enumerate}

Moreover, this family is unique up to multiplication by a constant.
\end{theorem}

\paragraph{Remark.} The family of Ledrappier's measures is of course closely related to the Gibbs states. With our convention for the kernel $k^F$ the typical trajectories of the geodesic flow for $\mu_F$ have the following description. Their lifts to $T^1N$ have:
\begin{itemize}
\item their past extremities which are distributed according to $\nu_o^F$;
\item their future extremities which are distributed according to $\nu_o^{\check{F}}$, where $\check{F}=F\circ\iota$, $\iota:T^1N\to T^1N$ being the involution $v\mapsto -v$.
\end{itemize}
Of course one can say more. Lifts of Gibbs states to $T^1N$ have a very nice description in terms of Ledrappier's measures in the so called \emph{Hopf's coordinates}. We won't need this description here and we refer to \cite{PPS} for more details. This is the reason why we called the cocycle $k^F$ the Gibbs cocycle.

\section{Gibbs measures and \emph{F}-harmonic measures for foliated bundles}
\label{sgibbsharmonicbundles}
The goal of this section is to define Gibbs measures, $F$-harmonic measures. We will prove their existence and give a bijective correspondence between the two.

In all the following, we consider a foliated bundle $(\Pi,M,B,X,\F)$  parametrized by a closed and negatively curved base $B$. As we saw in the previous section, the differential of $\Pi$ induces a foliated bundle $(\Pi_{\ast},T^1\F,T^1B,X,\hcF)$ and a foliated geodesic flow $G_t:T^1\F\to T^1\F$.

\subsection{Leafwise hyperbolic flow}
\label{lfwsehyperbolic}

It is possible to lift to the leaves of $\hcF$

\begin{itemize}
\item the distributions $E^{\star}$, $\star=s,u,cs$, or $cu$; the lifted distributions are denoted by $\bE^{\star}$ and are invariant by the flow $G_t$: this is the \emph{foliated hyperbolicity} \cite{BGM,Al3};
\item the invariant foliations $\W^{\star}$, $\star=s,u,cs$, or $cu$; the lifted foliations are denoted by $\bcW^{\star}$ and form subfoliations of $\hcF$ invariant by the flow $G_t$;
\item any Hölder continuous potential $F:T^1B\to\R$: the lifted potential $\bF=F\circ \Pi_{\ast}$ is continuous and remains uniformly Hölder continuous inside the leaves;
\item the families of measures $(\lambda^{\star}_{F,v})_{v\in T^1B}$ with $\star=s,u,cs$, or $cu$: we define the families $(\widehat{\lambda}^{\star}_{F,w})_{w\in T^1\F}$ as follows, if $w\in T^1\F$, $v=\Pi_{\ast}(w)\in T^1B$ and $D\dans \bW^{\star}(w)$, one has $\widehat{\lambda}^{\star}_{F,w}(D)=\lambda^{\star}_{F,v}(\Pi_{\ast}(D))$.
\end{itemize}

We can also lift the cocycles defined by Formulas \eqref{Eq:cocs} and \eqref{Eq:cocu}
\begin{equation}
\label{Eq:cocslift}
\widehat{k}_F^s(v,w)=\exp\left[\int_0^{\infty}(\bF\circ G_t(v)-\bF\circ G_t(w)) dt\right]=k_F^s(\Pi_{\ast}(v),\Pi_{\ast}(w));
\end{equation}
for $v,w$ in the same stable leaf $\bW^s$;
\begin{equation}
\label{Eq:coculift}
\widehat{k}_F^u(v,w)=\exp\left[\int_0^{\infty}(\bF\circ G_{-t}(w)-\bF\circ G_{-t}(v)) dt\right]=k_F^s(\Pi_{\ast}(v),\Pi_{\ast}(w)).
\end{equation}
for $v,w$ in the same unstable leaf $\bW^u$.
\subsection{Gibbs measures for the foliated geodesic}
\label{sleafwisegibbsmeasures}

\paragraph{Gibbs measures.} The family of measures $(\widehat{\lambda}^{u}_{F,v})_{v\in T^1\F}$ on the unstable manifolds is quasi-invariant by the foliated hyperbolic flow with cocycle relations given by

\begin{equation}
\label{Eq:cocyclestrongunstablelift}
\frac{d\left[G_T\,_\ast\widehat{\lambda}^{u}_{F,G_{-T}(v)}\right]}{d\widehat{\lambda}^{u}_{F,v}}(w)=\exp\left[\int_0^T (\bF\circ G_{-t}(v)-P(F))dt\right],
\end{equation}
for $T\in\R$, $v\in T^1\F$ and $w\in \bW^u(v)$. By analogy with the partially hyperbolic context (see the classical definition of Gibbs $u$-states defined in \cite{BV,BDV,PS}), we shall define the Gibbs property by specifying the class of the conditional measures in unstable manifolds.

\begin{defi}
\label{Defigibbs}
Let $F:T^1B\to\R$ be a Hölder continuous potential. A Gibbs measure of the foliated geodesic flow $G_t$ associated to the potential $\bF$ is a probability measure $\mu$ on $T^1\F$ such that
\begin{itemize}
\item $\mu$ is invariant by $G_t$;
\item $\mu$ has an absolutely continuous disintegration with respect to $(\bcW^u,\widehat{\lambda}^u_{F,v})$.
\end{itemize} 
\end{defi}

\paragraph{Existence.} The existence of Gibbs measures for $G_t$ follows from the following theorem. 

\begin{theorem}
\label{characgibbsmeasures}
Let $F:T^1B\to\R$ be a Hölder continuous potential. Then
\begin{enumerate}
\item for any small disc $D\dans\bW^u_{loc}(v)$, the accumulation points of
$$\frac{1}{T}\int_0^T\frac{G_t\,_\ast(\widehat{\lambda}^u_{F,v})_{|D}}{\widehat{\lambda}^u_{F,v}(D)} dt$$
are Gibbs measures for $G_t$ associated to the potential $\bF$;
\item ergodic components of Gibbs measures for $\bF$ are Gibbs measures for the same potential;
\item the densities $\widehat{\psi}_{F,v}^u$ of the conditional measure of any Gibbs measure $\mu$ on the unstable plaques $\bW^u_{loc}(v)$ with respect to $\widehat{\lambda}^u_{F,v}$ are uniformly log-bounded (i.e. bounded away from $0$ and $\infty$) and satisfy

\begin{equation}
\label{Eq:characgibbsmeasures1}
\frac{\widehat{\psi}_{F,v_0}^u(w)}{\widehat{\psi}_{F,v_0}^u(v)}=\widehat{k}^u_F(v,w),
\end{equation}
for $v_0\in T^1\F$, $v,w\in\bW^u_{loc}(v_0)$ and $T\in\R$;
\item the projection to $T^1B$ of a Gibbs measure for $G_t$ associated to $\bF$ is $\mu_F$.
\end{enumerate}
\end{theorem}

\begin{proof}
The proof of this theorem is a simple adaptation of Section 11.2.2. of \cite{BDV}, where the authors study the Gibbs $u$-states of partially hyperbolic systems. Here, we know that $\bcW^u$ is a foliation of a \emph{compact} manifold $T^1\F$ which is \emph{uniformly} expanded by 
$G_t$. By using the Hölder continuity of the lifted potential $\bF$ in the leaves of $\hcF$, we see that these results can be generalized without difficulty (the precise verifications can be found in the author's thesis \cite[Section VI.1.3]{Al5}).
\end{proof}

The next result shows that any Gibbs measure for $G_t$ yields a family of measures transverse to $\bcW^u$ with an absolute continuity property which is analogous to Property \eqref{Eq:cocycleweakstable}. Here, a transversal to $\bcW^{u}$ is of higher dimension than $\bcW^{cs}$. Typically, we can think of a transversal as being the product of a center stable manifold by a little open subset of the fiber $X$, which is transverse to $\hcF$.

\begin{proposition}
\label{transverseGibbs}
Let $F:T^1B\to\R$ be a Hölder continuous potential and $\mu$ be a Gibbs measure for $G_t$ associated to the potential $\bF$. Let $\A=(U_i,\phi_i)_{i\in I}$ be a foliated atlas for $\bcW^u$. Call $(T_i)_{i\in I}$ the associated complete system of transversals. Then, there is a family of measures $(\nu_i)_{i\in I}$ on these transversals such that
\begin{enumerate}
\item $\nu_i$ is equivalent to the projection on $T_i$ of the restriction $\mu_{|U_i}$ with a Radon-Nikodym derivative which is uniformly log-bounded;
\item it satisfies the following quasi-invariance relation
$$\frac{d[\h u{T_i}{T_j}\,_\ast\nu_i]}{d\nu_j}(v)=\widehat{k}_F^u(v,\h u{T_j}{T_i}(v)),$$
for $v$ in the domain of $\h u{T_j}{T_i}$.
\end{enumerate}
\end{proposition}

\begin{proof}
By Theorem \ref{characgibbsmeasures}, the conditional measures of $\mu_{|U_i}$ in the unstable plaques $\bW^u_{loc}(v)$ have the form $\widehat{\psi}^u_{F,v}\,\widehat{\lambda}^u_{F,v}$, $v\in T_i$.

The densities $\widehat{\psi}^u_{F,v}$ are uniformly log-bounded, so we can divide them by $\widehat{\psi}^u_{F,v}(v)$. Hence we are free to assume that $\widehat{\psi}^u_{F,v}$ is equal to $1$ in each $T_i$. By uniqueness of the disintegration, it determines the transverse measure $\nu_i$.

The Radon-Nikodym derivative with respect to the projection on $T_i$ of the restriction $\mu_{|U_i}$ is given by $v\mapsto\widehat{\psi}^u_{F,v}$ which is uniformly log-bounded.

Since the densities satisfy Relation \eqref{Eq:characgibbsmeasures1}, we obtain the desired cocycle relation by evaluating $\mu$ on an intersection $U_i\cap U_j$.
\end{proof}

\subsection{\emph{F}-harmonic measures}
\label{sfharmonic}

Recall that the Gibbs cocycle $k^F$ defined by Relation \eqref{gibbskernel} plays the role of the Poisson kernel. We are interested in the class of continuous functions on $N$, that will be called \emph{F-harmonic} and are defined below.

\begin{defi}
\label{Fharmonicfunction}
Let $F:T^1B\to\R$ be a Hölder continuous potential and $o\in N$. A positive function $h:N\to\R$ is said to be $F$-harmonic if there exists a finite measure $\eta_o$ on $N(\infty)$ such that for every $z\in N$
$$h(z)=\int_{N(\infty)}k^F(o,z;\xi)\,d\eta_o(\xi).$$
If $\Gamma\leq\Isom^+(N)$ acts on $N$ properly discontinuously and $h$ is a $\Gamma$-invariant $F$-harmonic function, the projection of $h$ to $N/\Gamma$ will still be called $F$-harmonic.
\end{defi}
\paragraph{Remark 1.} This definition doesn't depend on the point $o$, because of cocycle relation (\ref{Eq:Fcocyclerelation}): the definition will work with another point $o'\in N$ if we state $\eta_{o'}(\xi)=k_F(o,o';\xi)\eta_o(\xi)$.

\paragraph{Remark 2.} Recall the family $(\nu^F_z)_{z\in N}$ of \emph{Ledrappier's measures} defined in Theorem \ref{ledrappiermeasures}, which satisfy
\begin{itemize}
\item for all $\gamma\in\pi_1(B)$, we have $\gamma_\ast\nu^F_z=\nu^F_{\gamma z}$;
\item for all $y,z\in N,\xi\in N(\infty)$ $k^F(y,z;\xi)=d\nu_{z}^F/d\nu_{y}^F(\xi)$.
\end{itemize}
Then for any $z\in N$, $\mass(\nu^F_z)=\int_{N(\infty)}k^F(o,z;\xi)d\nu^F_o(\xi)$. Hence the function $z\mapsto\mass(\nu_z^F)$ is $F$-harmonic and $\pi_1(B)$-invariant (by equivariance of $(\nu^F_z)_{z\in N})$. Consequently, it induces a $F$-harmonic map on the quotient. This fonction on $B$ will be denoted by $h_0^F$. We will prove later that, up to a multiplicative constant, it is the unique $F$-harmonic function on $B$: this is Theorem \ref{uniquenessharmonicfunction}.

\paragraph{Remark 3.} It follows from the main theorem of \cite{Al4} that $F$-harmonic functions satisfy an analogue of Fatou's theorem of nontangential convergence. We deduce there that \emph{the integral representation of} $F$-\emph{harmonic functions is unique}.\\

We can now define the notion of $F$-harmonic measures for foliated bundles. Since leaves are quotients of $N$ (they are Riemannian covers of the base), it makes sense to talk about $F$-harmonic functions on them.
\begin{defi}
\label{Fharmonicmeasure}
Let $F:T^1B\to\R$ be a Hölder continuous function. A probability measure $m$ on $M$ will be called $F$-harmonic if it has Lebesgue disintegration in the leaves of $\F$ and if the densities in the plaques are $F$-harmonic.
\end{defi}
The existence of these measures comes from Theorem \ref{bijectionFgibbsFharmonic}, which is proven in the next paragraph.

\subsection{Bijective correspondence between Gibbs and \emph{F}-harmonic measures}
\label{unrolling}
\paragraph{Induced measures.} We fix a good atlas of $B$, $\A=(U_i,\phi_i)_{i\in I}$ which trivializes the bundle $\Pi:M\to B$, the unit tangent bundle $T^1B\to B$, as well as the covering map $N\to B$. We can lift this atlas to the unit tangent bundle $T^1B$ so as to obtain a new one denoted by $\widehat{\A}=(\widehat{U}_i,\widehat{\phi_i})_{i\in I}$, where the charts $\widehat{U}_i$ are trivially foliated by the tangent spheres.

Consider the fiber $X_{p_i}$, $p_i\in U_i$. The measure induced on the fiber $X_{p_i}$ by an $F$-harmonic measure $m$ on $M$ is by definition the projection of $m_{|U_i\times X}$ on $X_{p_i}$ along the plaques of $\F$.

The measure induced on $X_{p_i}$ by an $F$-Gibbs measure $\mu$ on $T^1\F$ is by definition the projection of $\mu_{|\widehat{U}_i\times{X}}$ along the plaques of $\hcF$.

\begin{theorem}
\label{bijectionFgibbsFharmonic}
Let $(\Pi,M,B,X,\F)$ be a projective foliated bundle parametrized by a closed and negatively curved Riemannian manifold $B$, whose fiber $X$ is a compact manifold. Let $F:T^1B\to\R$ be a Hölder continuous potential, and $\bF$ denote its lift to $T^1\F$.
\begin{enumerate}
\item For every Gibbs measure for $G_t$ associated to $\bF$, there is a unique $F$-harmonic measure for $\F$ that induces the same measures on the fibers $X_p$.
\item Reciprocally, for every $F$-harmonic measure for $\F$, there is a unique Gibbs measure for $G_t$ associated to $\bF$ that induces the same measures on the fibers $X_p$.
\end{enumerate}
\end{theorem}

\paragraph{Proof of Theorem \ref{bijectionFgibbsFharmonic}.} The proof of this is a copy of the proof of the Main Theorem of \cite{Al2}, where it is done for a special potential. Instead of copying verbatim the proof, let us show how we use the integral representation to lift measures in a very special case.

Consider a harmonic measure on the hyperbolic plane $\Hyp^2$ of the form $dm(z)=h(z)d\Leb(z)$ where $h(z)=\int_{\R\PP^1}k(o,z;\xi)d\eta(\xi)$ is a harmonic function and $k(o,z;\xi)$ is the Poisson kernel and $\eta$ is a finite Radon measure on $N(\infty)$.

The space $\Hyp^2\times\{\xi\}$ can be viewed as a center unstable manifold and carries the canonical measure $m_{\xi}=k(o,z;\xi)\,d\Leb_{\Hyp^2\times\{\xi\}}$. Note that $m_{\xi}$ is, up to a multiplicative constant, the unique measure on $\Hyp^2\times\{\xi\}$ invariant by the joint action of the geodesic and unstable horocyclic flow (see \cite[Proposition 1.1]{Al2}).

We obtain a measure $\mu$ on $T^1\Hyp$ by integration of measures $m_{\xi}$ against $\eta$. This measure projects down to $m$ and is invaraint by the joint action of the  geodesic and unstable horocyclic flows. Moreover the correspondence $m\mapsto\mu$ is a bijection between harmonic measures of $\Hyp^2$ and measures on $T^1\Hyp^2$ invariant by the joint action of the geodesic and horocyclic flow (see \cite[Proposition 1.1]{Al2} or \cite[Section IV.1]{Al5}for more details).

How to perform this ``unrolling argument'' for leaves of $\F$ is the purpose of \cite[Proposition 3.6]{Al2}. This is the first step of the proof where we associate to every $F$-harmonic measure $m$ a measure $\widehat{m}$ on $T^1\F$ which induces the same measure as $m$ in the fibers and is not a priori invariant by the foliated geodesic flow. The second step is a reparametrization argument given in \cite[Proposition 4.8]{Al2}. We associate to $\widehat{m}$ a Gibbs measure $\mu$ for $\bF$ inducing the same measure in the fibers. An argument ``à la Hopf'' given in \cite[Proposition 4.6]{Al2} shows that two different Gibbs measures induce different measures in the fibers which proves that this correspondence is bijective.

\quad \hfill $\square$

\paragraph{$F$-harmonic functions on a compact manifold.} By Remark 2 above, the function ``mass of Ledrappier measure'' on $N$ induces on the quotient $B$ a canonical $F$-harmonic function on $B$ denoted by $h_0^F$.

\begin{theorem}
\label{uniquenessharmonicfunction}
Let $B$ be a closed Riemannian manifold with negative sectional curvature. Let $F:T^1B\to\R$ be a Hölder continuous potential. Then any $F$-harmonic function on $B$ is obtained from $h_0^F$ by multiplication by a positive constant.
\end{theorem}

\begin{proof}
Theorem \ref{bijectionFgibbsFharmonic} is true even when $X$ is just a singleton: that is to say that there is a bijection between $F$-harmonic measures on $B$ and Gibbs states for the geodesic flow $g_t$ associated to $F$. But such a Gibbs state is unique (see Theorem \ref{lpsgibbsstates}).

Hence the density of the unique $F$-harmonic measure on $B$ is determined, up to multiplication by a positive constant. In other words, there is, up to multiplication by a positive constant, only one $F$-harmonic function on $B$: that is $h_0^F$.
\end{proof}

\paragraph{Remark 4.} Consider a covering $(U_i)_{i\in I}$ of $B$ such that $\Pi^{-1}(U_i)=U_i\times X$. Suppose there is a family $(\nu_i)_{i\in I}$ of holonomy-invariant measures on $X$. In $U_i\times X$, form $h_0^F\Leb_{U_i}\times\nu_i$. Since $(\nu_i)_{i\in I}$ is holonomy-invariant, these measures can be glued together so as to obtain a measure $m$ on $M$ which is $F$-harmonic. Thus in that sense $F$-harmonic measures, like Garnett's harmonic measures \cite{Gar}, generalize invariant measures.

\section{Uniqueness results in the case of projective bundles}
\label{suniqueness}

The goal of this section is to prove uniqueness results in the case of foliated bundles with $\C\PP^1$-fibers \emph{when there is no transverse measure invariant by holonomy}. More precisely, for the foliated geodesic flow, we prove that in this case:
\begin{enumerate}
\item there is a unique Gibbs measure in the sense of Definition \ref{Defigibbs} for the foliated geodesic flow associated to a given Hölder continuous potential;
\item there is a unique $F$-harmonic measure associated to a given Hölder continuous potential.
\end{enumerate}

\subsection{A locally constant projective cocycle}
\label{slocallyconstant}

In \cite{BGV}, Bonatti, G\'omez-Mont and  Viana noticed that the foliated geodesic flow $G_t$ is a \emph{locally constant projective cocycle} over $g_t$. Let us explain what we mean by that. By definition, we know that over a trivializing chart $U$, the orbits of $G_t$ are just copies of that of $g_t$, with the same parametrization. Hence the flow $G_t$ sends fibers on fibers. We shall call the resulting cocycle,
$$A_t(v)=(G_t)_{|X_v}:X_v\To X_{g_t(v)},$$
for any $t\in\R$ and $v\in T^1B$. Here the term \emph{cocycle} refers to the following identity:
$$A_{t_1+t_2}(v)=A_{t_1}(g_{t_2}(v))\circ A_{t_2}(v).$$

Note that if we choose any orbit segment $c=g_{[0,t]}(v)$, then $A_t(v)$ is the holonomy map along the path $c$. The cocycle is therefore projective: the lifted flow $G_t$ sends fibers to fibers as an element of $\PSL_2(\C)$. The cocycle is locally constant in the following sense. Fibers over a trivializing disc in $B$ can all be identified with $\C\PP^1$ so when we read the cocycle in this coordinate, the holonomy over a segment of orbit staying in that disc is the identity.

\subsection{A criterion for the existence of non-zero Lyapunov exponents}
\label{sbogovi}

\paragraph{Lyapunov exponents.} Oseledets' theorem ensures the existence of Lyapunov exponents for this cocycle on a Borel subset that is full for any measure invariant by $g_t$. If a point $v\in T^1B$ belongs to this Borel set, its maximal and minimal Lyapunov exponents are defined by the following formulas
$$\chi^+(v)=\lim_{t\to\infty}\frac{1}{t}\log(||A_t(v)||)\geq 0,$$
$$\chi^-(v)=\lim_{t\to\infty}\frac{1}{t}\log(||A_t(v)^{-1}||^{-1})\leq 0.$$

Here $||.||$ denotes any matrix norm. Note that $A_t(v)$ is a matrix up to a sign so $||A_t(v)||$ is well defined. Note also that we always have $\chi^+(v)=-\chi^-(v)$.

Using the symbolic representation of Anosov flows, Bonatti, G\'omez-Mont and Viana were able to prove:

\begin{theorem}[Bonatti, G\'omez-Mont, Viana]
\label{BoGoVi}
Let $(\Pi,M,B,\C\PP^1,\F)$ be a projective foliated bundle parametrized by a closed negatively curved base $B$. Then the following dichotomy holds:
\begin{itemize}
\item either there exists a probability measure on $\C\PP^1$ invariant under each element of the holonomy group;
\item or there is a Borel set $\X\dans T^1B$, full for all Gibbs states of $g_t$, such that $\chi^+(v)>0$ for every $v\in\X$.
\end{itemize}
\end{theorem}

In the what remains of the section, we shall assume the following.
\begin{center}
\emph{There is no} $\rho(\pi_1(B))$\emph{-invariant probability measure on} $\C\PP^1$.
\end{center}

\paragraph{Lyapunov sections.}

We can suppose that $\X$ trivializes the bundle (for example by taking away the boundaries of a finite number of discs which cover $M$ and trivialize the bundle). 

Hence the cocycle over $\X$ can be linearized: and we have a measurable linear bundle $\overline{\Pi}:\X\times\C^2\to\X$ whose projectivization is $\Pi$  and a linear cocycle $(t,v)\in\R\times\X\mapsto \overline{A}_t(v)$ over the geodesic flow whose projectivization is $A_t(v)$. If $v\in\X$ denote by $\overline{X}_v=\{v\}\times\C^2$, the linear fiber of $v$. Let $|.|$ be any norm in $\C^2$, $||.||$ be the corresponding matrix norm and $\dist_{X_v}$ the associated distance in $X_v\simeq\C\PP^1$.

The next proposition is an application of Theorem \ref{BoGoVi} and of Oseledets' theorem.

\begin{proposition}
\label{lines}
For every $v\in \X$ there exists a splitting of the linear fiber $\overline{X}_v=\overline{\sigma}^+(v)\oplus\overline{\sigma}^-(v)$ such that if $\sigma^{\pm}$ denote the projectivizations of $\overline{\sigma}^{\pm}$ the following assertions hold true
\begin{enumerate}
\item $\sigma^{\pm}$ and $\overline{\sigma}^{\pm}$ vary measurably with $v$;
\item these maps comute with the cocycle: for any $v\in \X$ and $t\in\R$, $A_t(v)\sigma^{\pm}(v)=\sigma^{\pm}(g_t(v))$ and $\overline{A}_t(v)\overline{\sigma}^{\pm}(v)=\overline{\sigma}^{\pm}(g_t(v))$;
\item for any $v\in\X$, we have the following contraction property:
$$\lim_{t\to\infty}\frac{1}{t}\log\,\dist_{X_v}(A_t(v)x,\sigma^+(g_t(v)))=-2\chi^+(v)\,\,\,\,\,for\,\,all\,\,x\in X_v\moins\{\sigma^-(v)\},$$
$$\lim_{t\to\infty}\frac{1}{t}\log\,\dist_{X_v}(A_{-t}(v)x,\sigma^-(g_{-t}(v)))=-2\chi^+(v)\,\,\,\,\,for\,\,all\,\,x\in X_v\moins\{\sigma^+(v)\}.$$
\item we have
$$\lim_{t\to\infty} |\overline{A}_{-t}(v)\omega|=0\,\,\,\,if\,\,and\,\,only\,\,if\,\,\omega\in \overline{\sigma}^+(v),$$
$$\lim_{t\to\infty} |\overline{A}_{t}(v)\omega|=0\,\,\,\,if\,\,and\,\,only\,\,if\,\,\omega\in \overline{\sigma}^-(v).$$
\end{enumerate}
\end{proposition}

\begin{defi}
The sections $\sigma^{\pm}$ of $\Pi$ are called the Lyapunov sections.
\end{defi}

We attribute the next proposition to Bonatti and G\'omez-Mont. They proved it in the case where the base is the unit tangent bundle of a compact hyperbolic surface and $g_t$ is the geodesic flow (see \cite[Proposition 3.1.]{BG}).

\begin{proposition}
\label{sections}
\begin{enumerate}
\item The two sections commute with the flows: $G_t\circ\sigma^{\pm}=\sigma^{\pm}\circ g_t$.
\item The section $\sigma^+$ preserves the unstable foliation, that is $\sigma^+(W^u(v))=\bW^u(\sigma^+(v))$ and commutes with the unstable holonomy maps.
\item The section $\sigma^-$ preserves the stable foliation, that is $\sigma^-(W^s(v))=\bW^s(\sigma^-(v))$ and commutes with the stable holonomy maps.
\end{enumerate}
\end{proposition}

\subsection{Lift the Gibbs states}
\label{sliftgibbsstates}

\paragraph{Cocycle relations.}
Remember the cocycle we defined in relation to the strong stable and unstable foliations of $G_t$ by Formulas (\ref{Eq:cocslift}), (\ref{Eq:coculift}). We will be interested in the ergodic properties of these foliations. More precisely, we study the existence and uniqueness of families of measures on the sets $\TT^+$ and $\TT^-$ of local transversals to $\bcW^{u}$ and $\bcW^{s}$, say $(\nu_{F,T}^+)_{T\in \TT^+}$ and $(\nu_{F,T}^-)_{T\in \TT^-}$, such that:

\begin{equation}
\label{Eq:cocycleweakstablelift}
\frac{d\left[\h u{T_1}{T_2}\,_\ast\nu^+_{F,T_1}\right]}{d\nu^{+}_{F,T_2}}(v)=\widehat{k}_F^u(v,\h u{T_2}{T_1}(v)),
\end{equation}
where $T_1,T_2\in\TT^+$ and $v$ lie in the domain of a holonomy map $\h u{T_2}{T_1}$ and:

\begin{equation}
\label{Eq:cocycleweakunstablelift}
\frac{d\left[\h s{T_1}{T_2}\,_\ast\nu^-_{F,T_1}\right]}{d\nu^{-}_{F,T_2}}(v)=\widehat{k}_F^s(v,\h s{T_2}{T_1}(v)),
\end{equation}
where $T_1,T_2\in\TT^-$ and $v$ lie in the domain of a holonomy map $\h s{T_2}{T_1}$.

\paragraph{Remark.} Suppose we know how to define a family of measures satisfying the cocycle relation (\ref{Eq:cocycleweakstablelift}) (resp. (\ref{Eq:cocycleweakunstablelift})) on a complete family of transversals to the unstable foliation (resp. stable foliation). Then, by an obvious adaptation \cite[Lemma 1.4]{BM}, we deduce a family of measures on the transversals to the unstable foliation (resp. stable foliation).

\paragraph{Families of measures commuting with the Lyapunov sections.} Here is an easy, but very important in the sequel, consequence of Proposition \ref{sections}. Let $F:T^1B\to\R$ be a Hölder continuous potential. Recall that in \S \ref{lfwsehyperbolic} we lifted the families of measures $(\lambda_{F,v}^{\star})_{v\in T^1B}$, $\star=s,u$ via $\Pi_{\ast}$. Note that by Proposition \ref{sections} for every $v\in\X$, the restriction of $\Pi_{\ast}$ to $\bW^u_{loc}(\sigma^+(v))$ (resp. $\bW^s_{loc}(\sigma^-(v))$) is a local inverse of $\sigma^+$ (resp. $\sigma^-$). Thus, we obtain the

\begin{lemma}
\label{sigmamu}
Let $F:T^1B\to\R$ be a Hölder continuous potential. Then the Lyapunov sections commute with the measures in the following sense: for any $v\in\X$,
$$\sigma^+\,_\ast\lambda^u_{F,v}=\widehat{\lambda}^u_{F,\sigma^+(v)}\,\,\,\,and\,\,\,\,\sigma^-\,_\ast\lambda^s_{F,v}=\widehat{\lambda}^s_{F,\sigma^-(v)}.$$
\end{lemma}

\paragraph{Lifts of the Gibbs state.}
Every Gibbs state gives full measure to the Borel set $\X$. So they can all be lifted thanks to the Lyapunov sections. If $\mu_F$ the Gibbs state of $F$, we will define
$$\mu_F^+=\sigma^+\,_\ast\mu_F,\,\,\,\,\,\,\,\,\,\,\,\,\,\,\,\mu_F^-=\sigma^-\,_\ast\mu_F.$$

Since the Lyapunov sections commute with the flows, \emph{these measures are invariant by} $G_t$. We will use the following theorem which is well known by the specialists and whose proof can be found in an Appendix of the author's thesis ( see section VII.4 of \cite{Al5}).
\begin{theorem}
\label{twoergodic}
The measures $\mu^+_F$ and $\mu^-_F$ are the only ergodic measures for $G_t$ which project down to $\mu_F$.
\end{theorem}

\paragraph{Disintegrations on stable and unstable plaques.} We now describe the local structures of $\mu^+$ and $\mu^-$.

Choose a finite cover of $T^1B$ by small open sets of the form $V_i=[W^u_{loc}(v_i),W^{cs}_{loc}(v_i)]$, which trivialize $\Pi_{\ast}:T^1\F\to T^1B$. We saw in Theorem \ref{lpsgibbsstates} that $(d\mu_F)_{|V_i}=(k^u(v,w)d\lambda^u_{v}(w))d\lambda^{cs}_{v_i}(v)$.

The preimages $T^{cs}_i=\Pi_{\ast}^{-1}(W^{cs}_{loc}(v_i))$ form a complete family of transversals of the strong unstable foliation in $T^1\F$. Moreover we can cover $T^1\F$ by the preimages $U_i=\Pi_{\ast}^{-1}(V_i)=\bigcup_{x\in T^{cs}_i} \bW^u_{loc}(x)$. Note that the lifts of the densities to $U_i$ are precisely the $\widehat{k}^u_F(.,.)$.

Of course, by a symmetric argument, if one changes the roles of the strong unstable and of the strong stable foliations, one can see that the sets $U_i$ are also filled with strong stable plaques, on which are defined local densities $\widehat{k}^s_F(.,.)$ and that $T^{cu}_i=\Pi_{\ast}^{-1}(W^{cu}_{loc}(v_i))$ is a complete system of transversals of $\bcW^s$.

Using the local product structure of Gibbs measures decribed above as well as the commutation properties obtained in Proposition \ref{sections}, we get the

\begin{proposition}
\label{lpsreleve}
Let $(U_i)_{i\in I}$ be a finite covering of $T^1\F$ obtained by lifting a covering of $T^1B$ by small open sets with the local product structure as explained above. For $\star=cs$ or $cu$, there exists a family of measures on $T^{\star}_i$ denoted by $\nu^{\star}_{F,v_i}$ such that
 
$$(d\mu_F^+)_{|U_i}=\left(\widehat{k}^u_F(v,w)d\widehat{\lambda}^u_{F,v}(w)\right)\,d\nu^{cs}_{F,v_i},\,\,\,\,\,\,\,\,\,\,\,\,\,\,\,\textrm{and}\,\,\,\,\,\,\,\,\,\,\,\,\,\,\,(d\mu_F^-)_{U_i}=\left(\widehat{k}^s_F(v,w)d\widehat{\lambda}^s_{F,v}(w)\right)\,d\nu^{cu}_{F,v_i}.$$

Moreover, the families $\nu^{cs}_{F,v_i}$ and $\nu^{cu}_{F,v_i}$ respectively satisfy the cocycle relations \eqref{Eq:cocycleweakstablelift} and \eqref{Eq:cocycleweakunstablelift}.
\end{proposition}

\begin{coro}
\label{mugibbsmeasure}
The measure $\mu_F^+$ is a Gibbs measure for the potential $\bF$.
\end{coro}

\paragraph{Remark.} We have been led to make a choice in the definition of a Gibbs measure. The important fact is that no matter the choice in the definition of a Gibbs measure, the measures $\mu^+$ and $\mu^-$ can't be both Gibbs measures, or we would have a measure on the fiber invariant by the action of the holonomy group. This is the object of our next paragraph.

\subsection{Singular disintegration on stable and unstable leaves}

\paragraph{Markov partition adapted to the bundle.} We said that the cocycle defined over $g_t$ is locally constant. This property will be important for us, so let us explain clearly what we mean by that. An idea presented in \cite{BGV} is to use Markov partitions associated to trivializing atlases for $\Pi_{\ast}:T^1\F\to T^1B$.

Take a good trivializing atlas $(V_1,...,V_k)$ in $T^1B$: each $\Pi_{\ast}^{-1}(V_i)$ is a chart that trivializes $\hcF$. Now take a Markov partition adapted to this family (see Theorem \ref{Markovpartition}). The rectangles of this partition are denoted by $R_i=[A^u_i,A^s_i]$ (see Definition \ref{rectangles}) and the cubes are denoted by $C_i$ (see Formula \eqref{cube} for the definition), for $i\in I$. Remember that the cubes are filled with compact unstable domains denoted by $A^u_i(v)$, $i\in I$, which are called Markovian unstable plaques.

By Theorem \ref{Markovpartition}, there are two surjective maps $\alpha,\beta:I\to\{1,...,k\}$, such that when $v\in C_i$, we have $v\in V_{\alpha(i)}$ and $g_{r(v)}(v)\in V_{\beta(i)}$. 

Recall that $\CC$ is the set of points whose returns in $\bigcup R_i$ are always in the interior of the rectangles (see Formula \eqref{Eq:residual}). By Theorem \ref{Markovpartition}, this set is $g_t$-invariant and full for any $g_t$-invariant probability measure which gives positive mass to open sets. 

\begin{center}
\emph{As a consequence, we will assume in the sequel that $\X\dans\CC$.}
\end{center}

Since the boundary of a cube $C_i$ is a union of rectangles (of measure $0$ for any flow invariant measure since they are topologically transverse to the flow) and of subsets of $\CC'=T^1B\moins\CC$, we have the following:

\paragraph{Claim.} \emph{For every Gibbs state $\mu_F$ in $T^1B$, we have $\mu_F(\X\cap\bigcup_i \Int\,C_i)=1$.}

\paragraph{Trivialization over the Markov partition.} Since the interior of a cube $C_i$ is included in the trivializing chart $V_{\alpha(i)}$, there is an identification $\Pi_{\ast}^{-1}(\Int\,C_i)\simeq\Int\,C_i\times\C\PP^1$ which identifies $\hcF_{|\Pi_{\ast}^{-1}(\Int\,C_i)}$ with the partition $(\Int\,C_i\times\{x\})_{x\in\C\PP^1}$. Let us describe  $A_t$ in these new coordinates.

Let $v\in\X\cap\Int\,C_i$. If $g_{[0,t]}(v)\dans\Int\,C_i$ then $A_t(v)$ is the identity. By definition of the Markov partition, $A_{r(v)}(v)$ is the holonomy map $\tau_{\alpha(i)\beta(i)}$. Since $v\in\X\dans\CC$ never meets the boundary of a rectangle, and since rectangles are topologically transverse to the flow, we have $A_{r(v)+\eps}(v)=\tau_{\alpha_i\beta_i}$ for $\eps>0$ sufficiently small. Hence we see that $A_t(v)$ only depends on the sequence of cubes visited by the orbit of $v$. This justifies the terminilogy of locally constant cocycle. We deduce the following

\begin{lemma}
\label{control}
Let $v\in \Int\,C_i$ for $i\in I$. Then for every $t^-\leq 0$, $A_{t^-}(w)$ is independent of $w\in\Int\,A^u_i(v)$. Similarly for every $t^+\geq 0$, $A_{t^+}(w)$ is independent of $w\in\Int\,A^s_i(v)$.
\end{lemma}
\begin{proof}
This is due to the Markov property (see Definition \ref{Markovproperty}). If $w\in\Int\,A^u_i(v)$, then for every $t^-\leq 0$, $g_{[t^-,0]}(v)$ and $g_{[t^-,0]}(w)$ visit exactly the same cubes of the partition. Since the value of the cocycle $A_{t^-}$ depends only on the succession of cubes visited by the trajectory in the past between $0$ and $t^-$, it has to be independent of the choice of $w$.

The same argument in the future gives the proof of the last part of the lemma.
\end{proof}

\paragraph{Lyapunov sections as local graphs.} In our trivialization the two sections $\sigma^+$ and $\sigma^-$ can be written, in restriction to $\Int\,C_i$, as $(Id,s^+_i)$ and $(Id,s^-_i)$, where $s^{\pm}_i:\X\cap\Int\,C_i\to\C\PP^1$ are measurable maps. As a consequence of Proposition \ref{sections}, we get the

\begin{lemma}
\label{codesections}
Let $i\in I$. Then, for $\mu_F$-almost every $v\in\Int\,C_i$, 
\begin{enumerate}
\item the map $s_i^+$ is constant in $W^{cu}_{loc}(v)\cap\Int\,C_i$;
\item the map $s_i^-$ is constant in $W^{cs}_{loc}(v)\cap\Int\,C_i$;
\end{enumerate}
\end{lemma}

\paragraph{Singular disintegrations.} The following proposition is the main result of this section. It states that in the absence of transverse invariant measure, one can't specify the class of conditional measures in both stable and unstable manifolds.

\begin{proposition}
\label{singular}
Let $(\Pi,M,B,\C\PP^1,\F)$ be a projective foliated bundle parametrized by a closed and negatively curved base. Assume that there is no probability measure invariant by the action of the holonomy group. Then:
\begin{enumerate}
\item $\mu^+_{F}$ has a singular disintegration with respect to $(\bcW^s,\widehat{\lambda}^s_{F,v})$.
\item $\mu^-_{F}$ has a singular disintegration with respect to $(\bcW^u,\widehat{\lambda}^u_{F,v})$.
\end{enumerate}
\end{proposition}

We only prove the second assertion: the first one follows from the same argument. We first need the following lemma.

\begin{lemma}
\label{CNS}
The measure $\mu^-_{F}$ has a singular disintegration with respect to $(\bcW^u,\widehat{\lambda}^u_{F,x})$ if and only if for every $i\in I$ and almost every $v\in\X\cap\Int \,C_i$, there is no Borel set $D\dans A_i^u(v)$ with positive $\lambda^u_{F,v}$ measure on which $s^-_i$ is constant.
\end{lemma}

\begin{proof}
Denote by $U_i$ the preimage $\Pi_{\ast}^{-1}(\Int\,C_i)$. Each cube is filled with strong stable and strong unstable plaques. Hence, $U_i$ is filled with the preimages of these plaques by the fibration $\Pi_{\ast}$. One can write
$$U_i=\bigsqcup_{w\in W^{cs}_{loc}(v)\cap \Int\,C_i} T^u_i(w),$$
where by definition, $T^u_i(w)=\Pi_{\ast}^{-1}(W^u_{loc}(w)\cap\Int \,C_i)$.

Since $(U_i)_{i\in I}$ is a partition up to a zero measure set, $\mu_F^-$ has a singular disintegration with respect to $(\bcW^u,\bar{\lambda}^u_{F,x})$ if and only if all the restrictions $(\mu^-_F)_{|U_i}$, $i\in I$ have.

Using Lemma \ref{codesections}, we see that if we choose $i\in I$ and $v_i\in\Int\,C_i\cap\X$, the conditional measures of $(\mu^-_{F})_{|U_i}$ on the $T^u(v)$ for $v\in W^{cs}_{loc}(v_i)\cap\Int\,C_i$, are given by:
$$\sigma^-\,_\ast(k^u_F(v,w)d\lambda^u_{F,v}(w)).$$
The partition of $U_i$ by the local unstable manifolds is a ``subpartition'' of that by the $T^u(v)$. Uniqueness of the disintegration assures that for $\mu_F^-$-almost every $\omega$, the conditional measure of $\mu^-_F$ in $\bW^u_{loc}(\omega)$ coincides with that of $\sigma^-\,_\ast(k^u_F(v,w)d\lambda^u_{F,v}(w))$, where $\omega\in T^u(v)$.

Now, we are led to disintegrate a measure supported by a graph. An elementary result of disintegration theory shows that since, $\sigma^-_{|C_i}=(Id,s^-_i)$, the following equivalence holds. The measure $\sigma^-\,_\ast(k^u_F(v,w)d\lambda^u_{F,v}(w))$ has singular disintegration with respect to $(\bcW^u,\widehat{\lambda}^u_{F,v})$ if and only if there is no $D\dans A^u_i(v)$ with positive $\lambda^u_{F,v}$ measure on which $s^-_i$ is constant.
\end{proof}

The proof of the following lemma is a bit technical and will be given in the next paragraph. After the statement of this lemma, we show how it implies the main proposition.

\begin{lemma}
\label{lemmaprincipal}
Assume that there exists $i_0\in I$, $v\in\Int\,C_{i_0}\cap\X$ and a Borel set $D\dans A^u_{i_0}(v)$ such that:
\begin{itemize}
\item $\lambda^u_{F,v}(D)>0$;
\item $s^-_{i_0}$ is constant on $D$.
\end{itemize}
Then, for all $i\in I$, $s^-_i$ is constant $\mu_F^-$-almost everywhere on $\Int\,C_i$.
\end{lemma}

\paragraph{End of the proof of Proposition \ref{singular}.} We proceed by contradiction and prove that if the conclusion of the proposition fails to be true, then the disintegration of $\mu_F^-$ in the fibers of $\Pi_{\ast}:T^1\F\to T^1B$ is invariant by all the holonomy maps of the foliation $\hcF$. In particular, there exists a measure on $\C\PP^1$ invariant by the action of $\rho(\pi_1(B))$, contradicting the hypothesis.

Consider the disintegration $(\theta_v)_{v\in\X}$ of $\mu_F^-$ on the fibers $X_v$, $v\in\X$. The measures $\theta_v$ are defined for $\mu_F$-almost every $v\in\X$ and are equal to the Dirac mass $\delta_{\sigma^-(v)}$. By Proposition \ref{sections}, we know that this family of measures is invariant by any holonomy map over any path inside $\W^{cs}$.

Now assume that the conclusion of Proposition \ref{singular} does not hold: the disintegration of $\mu^-_F$ on unstable plaques is not everywhere singular with respect to $(\widehat{\lambda}^u_{F,v})$. By Lemma \ref{CNS}, this implies the existence of $i_0\in I$ and of a Borel subset $D$ included in some unstable domain of $\Int\,C_{i_0}$ with nonzero $\lambda^u_{F,v}$ measure on which $s_{i_0}^-$ is constant. Finally, by Lemma \ref{lemmaprincipal}, it implies that on each $\Int\,C_i$, the map $s_i^-$ is constant almost everywhere. If one prefers, for every $i$ the conditional measures $(\theta_v)_{v\in\X\cap \Int\,C_i}$ are invariant by the holonomy maps over paths staying in $\Int\,C_i$.

Note that any path of $T^1B$ is homotopic to a concatenation of paths inside $\W^{cs}$ and paths inside $\W^u$. Hence, in order to prove that the family $(\theta_v)_{v\in\X}$ is holonomy-invariant, it is enough to prove that it is invariant by holonomy maps over paths inside $\W^u$.

It is classical that $\mu_F$ is ergodic and charges all open sets. Hence a $\mu_F$-typical orbit is dense. We deduce that $\mu_F$ gives zero measure to the set of periodic orbits. Finally we can assume that $W^{cu}(v)$ is simply connected for every $v\in\X$.

Let $v,w\in\X$ lying in the same unstable manifold and consider any unstable path $\gamma$ between them. Its image by a sufficiently long iteration of the flow in the past, is a path $\gamma'=g_{-t}(\gamma)$ lying in the interior of a single cube $C_i$. Since $W^{cu}(v)$ is simply connected, the holonomy of $\hcF$ over $\gamma$ is the composition of the holonomy maps over $g_{[0,t]}(g_{-t}(w))$, $\gamma'$ and $g_{[-t,0]}(v)$. Since the holonomy maps over orbit segments leave invariant the disintegration of $\mu_F^-$ and that the holonomy maps over paths inside $C_i$ leave invariant the family $(\theta_v)_{v\in\X\cap C_i}$, we proved that the holonomy along $\gamma$ sends $\theta_v$ on $\theta_w$, thus concluding the proof of the proposition. \quad \hfill $\square$.

\subsection{Proof of Lemma \ref{lemmaprincipal}}

\paragraph{Control of distortion.} In the sequel, we will need the following classical distortion lemma, which we will state without proof.

\begin{lemma}
\label{distortionlemma}
Let $F:T^1B\to\R$ be a Hölder continuous potential. Let $V$ be an unstable open set. Then given a positive number $\Delta$, there is a constant $K_0>1$ such that for any $t\geq 0$, any open subset $O\dans g_t(V)$ of diameter smaller than $\Delta$ and any Borel subset $Z\dans O$, we have, if $w\in O$:
$$\frac{1}{K_0}\frac{\lambda^u_{F,g_t(w)}(Z)}{\lambda^u_{F,g_t(w)}(O)}\leq\frac{g_t\,_\ast \lambda^u_{F,w}(Z)}{g_t\,_\ast\lambda^u_{F,w}(O)}\leq K_0\frac{\lambda^u_{F,g_t(w)}(Z)}{\lambda^u_{F,g_t(w)}(O)}.$$
\end{lemma}

\paragraph{Selection of good markovian plaques.}
Until the end of this paragraph, we assume that there is a set $D\dans A^u_{i_0}(v)$, for $v\in\Int\,C_{i_0}$ which is positive for $\lambda^u_{F,v}$ and on which $s^-_{i_0}$ is constant.

Fix a number $\delta>0$ smaller than all the $(\mu_F(C_i)/4)^2$, $i\in I$. By the Borel density theorem, there exists a small unstable ball $V\dans A^u_{i_0}(v)$ such that $\lambda^u_{F,v}(V\moins D)/\lambda^u_{F,v}(V)\leq \frac{9}{10}\delta$. Moreover, we can choose $V$ small enough in such a way that $\lambda^u_{F,v}(\partial V)=0$.

\begin{lemma}
\label{unionplaques}
There is a positive number $T_0$ such that for $t\geq T_0$, there exists $V_t\dans V$ such that
\begin{enumerate}
\item $\lambda^u_{F,v}(V_t)/\lambda^u_{F,v}(V)\leq\frac{9}{10}$;
\item $g_t(V_t)$ is a disjoint union of unstable markovian plaques $A^u_i(w)$.
\end{enumerate}
\end{lemma}

\begin{proof}
Let $i\in I$ and $t>0$. We say that a connected component of $\Int\,C_i\cap g_t(V)$ is \emph{markovian} if it is an unstable markovian plaque $A^u_i(w)$ for some $w$. It means that it crosses the cube $C_i$. Then,
$$\Int\,C_i\cap g_t(V)=\widehat{V}_{M,i}\sqcup\widehat{V}_{NM,i},$$
where $\widehat{V}_{M,i}$ denotes the union of markovian components of $\Int\,C_i\cap g_t(V)$ and $\widehat{V}_{NM,i}$, its complement (the union of non markovian components of the intersection).

Let $\Delta>0$ be greater than all diameters of the markovian plaques. The set $\bigcup_i\widehat{V}_{NM,i}$ is included in the $\Delta$-neighbourhood of $\partial g_t(V)$. Since this set is included in an unstable leaf, $g_{-t}\left(\bigcup_i\widehat{V}_{NM,i}\right)$ is included in the $\Delta (b/a) e^{-at }$-neighbourhood of $\partial V$. Since the decreasing intersection of these neighbourhoods is $\partial V$ and is null for $\lambda_{F,v}^u$, we can take $T_0>0$ such that for any $t\geq T_0$,
$$\lambda^u_{F,v}\left[g_{-t}\left(\bigcup_i\widehat{V}_{NM,i}\right)\right]\leq \frac{\lambda^u_{F,v}(V)}{10}.$$
It is now obvious that if $t\geq T_0$, the set $V_t=g_{-t}\left(\bigcup_i\widehat{V}_{M,i}\right)$ suits.
\end{proof}

\paragraph{Remark 1.} Since the flow preserves the class of $\lambda^u_{F,v}$ and each $\partial A^u_i(v)$ is null, for any $t\geq T_0$, there exists a collection $(D_j)_{j\in J}$ of disjoint open subsets of $V_t$ such that:
\begin{itemize}
\item for all $j\in J$, $g_t(D_j)$ is a markovian plaque;
\item $\lambda^u_{F,v}\left(V_t\moins\bigcup_j D_j\right)=0$.
\end{itemize}

\paragraph{Remark 2.} The first assertion of Lemma \ref{unionplaques} and the choice of $V$, imply that when $t\geq T_0$
$$\frac{\lambda^u_{F,v}(V_t\moins D)}{\lambda^u_{F,v}(V_t)}\leq\delta.$$

\begin{lemma}
\label{mixing}
There exists $T_1\geq T_0$ such that for any $t\geq T_1$ and any $i\in I$,
$$\frac{\lambda^u_{F,v}(V_t\cap g_{-t}(C_i))}{\lambda^u_{F,v}(V_t)}\geq\frac{\mu_F(C_i)}{2}.$$
\end{lemma}

\begin{proof}
We know that the following family of measures
$$\mu_t=\frac{g_t\,_\ast(\lambda^u_{F,v})_{|V}}{\lambda^u_{F,v}(V)},$$
converges to $\mu_F$ (see Theorem \ref{reconstructthegibbsstate}) as $t$ tends to infinity. In particular, since for every $i\in I$, $\mu_F(\partial C_i)=0$, we have $\lim_{t\to\infty}\lambda^u_{F,v}(V\cap g_{-t}(C_i))/\lambda^u_{F,v}(V)=\mu_F(C_i)$. But the following inequality holds for any $t>0$
$$\frac{\lambda^u_{F,v}(V_t\cap g_{-t}(C_i))}{\lambda^u_{F,v}(V_t)}\geq\frac{\lambda^u_{F,v}(V\cap g_{-t}(C_i))}{\lambda^u_{F,v}(V)}-\frac{\lambda^u_{F,v}(V\moins V_t)}{\lambda^u_{F,v}(V)}.$$
The second term of the difference tends to zero (see the proof of Lemma \ref{unionplaques}) as the first one converges to $\mu_F(C_i)$. Hence, the lemma follows.
\end{proof}

Now we choose $t\geq T_1$. Recall the content of Remark 1: there is a partition of $V_t$ modulo $\lambda^u_{F,v}$ denoted by $(D_j)_{j\in J}$ such that all $g_t(D_j)$ are markovian plaques.

\begin{lemma}
\label{Markovinequality}
Let $J_0$ denote the set of $j\in J$ such that $\lambda^u_{F,v}(D_j\moins D)/\lambda^u_{F,v}(D_j)\geq\sqrt{\delta}$. Then
$$\frac{\lambda^u_{F,v}\left(\bigcup_{j\in J_0} D_j\right)}{\lambda^u_{F,v}(V_t)}\leq\sqrt{\delta}.$$
\end{lemma}

\begin{proof}
The proof of this lemma is a simple application of an inequality ``à la Markov''. Since $(D_j)_{j\in J}$ is a partition of $V_t$ modulo $\lambda^u_{F,v}$, we have
$$\frac{\sum_{j\in J} b_j a_j}{\sum_{j\in J} b_j}=\frac{\lambda^u_{F,v}(V_t\moins D)}{\lambda^u_{F,v}(V_t)},$$
where for $j\in J$, $a_j=\lambda^u_{F,v}(D_j\moins D)/\lambda^u_{F,v}(D_j)$ and $b_j=\lambda^u_{F,v}(D_j)$. In particular, this quotient is smaller than $\delta$.

Now, since $J_0$ consists of those $j$ such that $a_j\geq\sqrt{\delta}$, we have the following chain of inequalities
$$\sqrt{\delta}\frac{\sum_{j\in J_0} b_j}{\sum_{j\in J} b_j}\leq\frac{\sum_{j\in J_0} b_j a_j}{\sum_{j\in J}b_j}\leq\delta.$$
The lemma follows because $\sum_{j\in J_0} b_j=\lambda^u_{F,v}\left(\bigcup_{j\in J_0} D_j\right)$.
\end{proof}

\begin{lemma}
\label{findarectangle}
There is a constant $K_0>1$, independent of $t$ and $\delta$, such that for any $i\in I$ and $t\geq T_1$, there exists $v_i\in\Int\,C_i$ such that
$$\frac{\lambda^u_{F,v_i}(A_i^u(v_i)\moins g_t(D))}{\lambda^u_{F,v_i}(A^u_i(v_i))}\leq K_0\sqrt{\delta}.$$
\end{lemma}

\begin{proof}
Let $i\in I$ and $t\geq T_1$. From Lemma \ref{mixing}, we know that inside $V_t$ the proportion of sets $D_j$ whose image by $g_t$ is a markovian plaque of $C_i$ is more than $\mu_F(C_i)/2$. From Lemma \ref{Markovinequality}, we also know that inside $V_t$ the proportion of sets $D_j$ whose intersection with the complement of $D$ weights more than $\sqrt{\delta}$ of its total mass, is less than $\sqrt{\delta}$. Moreover, we have chosen a $\delta$ in such a way that $\sqrt{\delta}\leq\mu_F(C_i)/4$.

From this, we deduce that there exists a set $D_j$ such that
\begin{itemize}
\item there is $v_i\in\Int\,C_i$ such that $g_t(D_j)=A^u_i(v_i)$;
\item $\lambda^u_{F,v}(D_j\moins D)/\lambda^u_{F,v}(D_j)\leq\sqrt{\delta}$.
\end{itemize}
Now, the distortion lemma \ref{distortionlemma} allows us to conclude
$$\frac{\lambda^u_{F,v_i}(A^u_i(v_i)\moins g_t(D))}{\lambda^u_{F,v_i}(A^u_i(v_i))}\leq K_0\frac{g_t\,_\ast\lambda^u_{F,v}(A^u_i(v_i)\moins g_t(D))}{g_t\,_\ast\lambda^u_{F,v}(A^u_i(v_i))}=K_0\frac{\lambda^u_{F,v}(D_j\moins D)}{\lambda^u_{F,v}(D_j)}\leq K_0\sqrt{\delta}.$$
The proof is now over.
\end{proof}

\begin{lemma}
\label{constantsection}
Let $t\geq T_1$, $i\in I$ and $v_i\in\Int\,C_i$ as in Lemma \ref{findarectangle}. The map $s^-_i$ is constant on $A_i^u(v_i)\cap g_t(D)$.
\end{lemma}

\begin{proof}
In order to see that, use Lemma \ref{control}. The value of the cocycle $A_{-t}$ is constant on $A_i^u(v_i)$, for any $w\in \Int\,C_i$. Hence if $w_1,w_2\in A_i^u(v_i)\cap g_t(D)$, we have
$$s^-_i(w_1)=A_{-t}(w_1)^{-1}s_{i_0}^-(g_{-t}(w_1))=A_{-t}(w_2)^{-1}s^-_{i_0}(g_{-t}(w_2))=s^-_i(w_2).$$
\end{proof}

\paragraph{An argument à la Hopf.} We know that in each cube $C_i$ the map $s_i^-$ is constant on a large proportion of some markovian unstable plaque. Using the absolute continuity of the center stable foliation, as well as the invariance of $\sigma^-$ by the holonomy maps along small center stable paths, we will show that the map $s_i^-$ is constant on a large proportion of $C_i$. This is the content of the following

\begin{lemma}
\label{saturatecs}
There exists a constant $K_1>1$, independent of the number $\delta$, such that for $i\in I$, there is a Borel set $O_i\dans \Int\,C_i$ verifying
\begin{enumerate}
\item $O_i$ is saturated in $C_i$ in the center stable direction;
\item $s^-_i$ is constant on $O_i$;
\item for any $w\in O_i$, we have
$$\frac{\lambda^u_{F,w}(A_i^u(w)\moins O_i)}{\lambda^u_{F,w}(A^u_i(w))}\leq K_1\sqrt{\delta}.$$
\end{enumerate}
\end{lemma}

\begin{proof}
Let $i\in I$ and $t\geq T_1$. By Lemma \ref{findarectangle}, there is a point $v_i\in\Int\,C_i$ such that the $\lambda^u_{F,v_i}$-proportion in the markovian plaque $A^u_i(v_i)$ of the complement of $g_t(D)$ is smaller than $K_0\sqrt{\delta}$ for some $K_0>1$. We can define $O_i\dans \Int\,C_i$ to be the saturated set of $A^u_i(v_i)\cap g_t(D)$ in the center stable direction.

By Lemma \ref{sections}, the map $s^-_i$ is constant in the center stable plaques of $C_i$. Hence since it is also constant on $A^u_i(v_i)\cap g_t(D)$, it is constant on $O_i$.

Finally, if $w\in W^{cs}_{loc}(v_i)$, we have
$$A_i^u(w)\moins O_i=\h {cs}{v_i}q(A^u_i(v_i)\moins g_t(D)).$$

Moreover, the holonomy maps of $\W^{cs}$ satisfy the following property of absolute continuity. Assume that $w'\in A^u(w)$ and set $w''=\h {cs}w{v_i}(w')$. Denote by $w'_s$ the projection on $W^s_{loc}(w'')$ along the flow. Take $T\in\R$ such that $g_T(w')=w'_s$. We have

$$\frac{d\left[\h {cs}{v_i}w\,_\ast\lambda^u_{F,v_i}\right]}{d\lambda^u_{F,w}}(w')=\exp\left[\int_0^{\infty}(F\circ g_{-t}(w'')-F\circ g_{-t}(w'_s))dt\right]\exp\left[\int_0^T(F\circ g_t(w')-P(F))dt\right].$$

Since the diameters of center stable plaques inside $C_i$ are uniformly bounded and since $F$ is bounded and uniformly Hölder continuous, the Radon-Nikodym derivatives above are uniformly bounded. Consequently, the third assertion follows.
\end{proof}

\begin{lemma}
\label{bigOi}
There exists a constant $K_2>1$, independent of $t$ and $\delta$, such that for $i\in I$ and $t\geq T_1$, if $O_i$ is the Borel set constructed in \ref{saturatecs}, we have
$$\frac{\mu_F(C_i\moins O_i)}{\mu_F(C_i)}\leq K_2\sqrt{\delta}.$$
\end{lemma}

\begin{proof}
Remember that in $C_i$, $\mu_F$ has a local product structure (see Theorem \ref{lpsgibbsstates}): it is obtained by integration of the $k^u_F(w,.)\lambda^u_{F,w}$ against $\lambda^{cs}_{F,v_i}$. Since all $k^u_F(w,.)$ are uniformly log-bounded on local unstable manifolds (see Theorem \ref{lpsgibbsstates}), the latter measures are equivalent to the normalized restriction $(\lambda^u_{F,w})_{|A^u_i(w)}/\lambda^u_{F,w}(A^u_i(w))$ with uniformly log-bounded densities.

Now a Fubini argument allows us to conclude the proof of the lemma.
\end{proof}

\paragraph{End of the proof of Lemma \ref{lemmaprincipal}}
Lemmas \ref{saturatecs} and \ref{bigOi} give for any $\delta>0$, a Borel set $O_i$ in each cube $C_i$ such that
\begin{itemize}
\item $s^-_i$ is constant on $O_i$;
\item $\mu_F(C_i\moins O_i)/\mu_F(C_i)\leq K_2\sqrt{\delta}$ for a number $K_2>1$ independent of $\delta$.
\end{itemize}

We see that the value of the constant does not depend on $\delta$ sufficiently small, since two different $O_i$ have to intersect, their complementaries weighting less than $1/2$ of the mass of $C_i$ when $\delta$ is small enough. Call $s^*_i$ this constant. The complement of $(s_i^-)^{-1}(s_i^*)$ is of arbitrarily small measure: it has zero measure for $\mu_F$ and the lemma is proven. \quad \hfill $\square$.

\subsection{Uniqueness of Gibbs and $F$-harmonic measures}

\paragraph{Uniqueness of the Gibbs measure.} We can end the proof of Theorem \ref{tripleuniquegibbs}. Assume that there is no transverse holonomy invariant measure. Let $F:T^1B\to\R$ be a Hölder continuous potential. We know from Corollary \ref{mugibbsmeasure} that $\mu_F^+$ is a Gibbs measure associated to $\bF$ as defined in Definition \ref{Defigibbs}.

Now, any Gibbs measure $\mu$ for $G_t$ associated to $\bF$ projects down to $\mu_F$. Since the ergodic components of any Gibbs measure associated to $\bF$ are also Gibbs measures (see Theorem \ref{characgibbsmeasures}), Theorem \ref{twoergodic} implies the following alternative. Either $\mu^-_F$ is a Gibbs measure associated to $\bF$, or $\mu^+_F$ is the only Gibbs measure.

By Proposition \ref{singular}, the first possibility does not occur since $\mu^-_F$ has a singular disintegration with respect to $(\bcW^u,\bar{\lambda}^u_{F,v})$. This ends the proof of the theorem.

\paragraph{Uniqueness of the $F$-harmonic measure.} The uniqueness of the $F$-harmonic measure, as defined in Definition \ref{Fharmonicmeasure}, also follows. Indeed Theorem \ref{bijectionFgibbsFharmonic} provides a bijective correspondence between $F$-harmonic measures and Gibbs measures associated to $\bF$.\quad \hfill $\square$

\subsection{Unique ergodicities of invariant foliations.}  Recall that we are looking for the measures $(\nu^+_T)_{T\in\TT^+}$ and $(\nu^-_T)_{T\in\TT^-}$, where $\TT^+$, $\TT^-$ are respectively the sets of local transversals to the strong unstable and the strong stable foliations, for which Relations \eqref{Eq:cocycleweakstablelift} and \eqref{Eq:cocycleweakunstablelift} hold.

\begin{theorem}
\label{uniqueergodicitiesstableunstable}
Let $(\Pi,M,B,\C\PP^1,\F)$ be a projective foliated bundle parametrized by a closed and negatively curved base $B$. Assume moreover that no probability measure on $\C\PP^1$ is invariant by the action of the holonomy group. Let $F:T^1B\to\R$ be a Hölder continuous potential.

Then, up to a multiplicative constant, there are unique families of measures $(\nu^+_T)_{T\in\TT^+}$ and $(\nu^-_T)_{T\in\TT^-}$ defined respectively on the local transversals to $\W^u$ and $\W^s$ which satisfy Relations \eqref{Eq:cocycleweakstablelift} and \eqref{Eq:cocycleweakunstablelift}.

Moreover, on the complete systems of transversals
$$T^{cs}(v)=\Pi_{\ast}^{-1}(W^{cs}_{loc}(v))\,\,\,\,\,and\,\,\,\,\,T^{cu}(v)=\Pi_{\ast}^{-1}(W^{cu}_{loc}(v)),$$
they are respectively given by $\nu^{cs}_{F,v}=\sigma^+\,_\ast\lambda^{cs}_{F,v}$ and $\nu^{cu}_{F,v}=\sigma^-\,_\ast\lambda^{cu}_{F,v}$ (see Proposition \ref{lpsreleve}).
\end{theorem}

We will proceed by contradiction. Assume for example that another family of measures (i.e. singular to that defined by Proposition \ref{lpsreleve}) denoted by $(\nu^+_T)_{T\in\TT^+}$ exists. In the base $T^1B$, consider a foliated atlas for $\W^u$, denoted by $(V_i,\phi_i)_{i\in I}$, with a complete system of transversals that consists of local center stable manifolds $(W^{cs}_{loc}(v_i))_{i\in I}$. Assume moreover that for $0\leq t\leq 1$, the preimages $g_{-t} (V_i)$ also form a foliated atlas of $\W^u$ and trivialize the bundle. Then, if $U_i=\Pi_{\ast}^{-1}(V_i)$ and $T_i=\Pi_{\ast}^{-1}(W^{cs}_{loc}(v_i))$, every $G_{-t}(U_i)$, $0\leq t\leq 1$ forms a foliated atlas of $\bcW^u$ with systems of transversals $G_{-t}(T_i)$.

\begin{lemma}
\label{propertymu}
There exists a probability measure $\mu$ on $T^1\F$ which in restriction to $U_i$ is obtained by integrating against $d\nu^+_{T_i}$ the measures $\widehat{k}^u_F(w_i,w)d\widehat{\lambda}^u_{F,w_i}(w)$ ($w_i\in T_i$). Moreover, such a measure $\mu$ projects down to $\mu_F$, the Gibbs measure for $g_t$ in $T^1B$ associated to $F$ (compare with Proposition \ref{lpsreleve}). Finally it is singular with respect to $\mu^+_F$.
\end{lemma}

\begin{proof}
We can define a measure on each $U_i$ as suggested in the lemma. The issue is to see that all these measures can be glued together. The measures $\nu^+_{T_i}$ satisfy the cocycle condition \eqref{Eq:cocycleweakstablelift} so we have for every $i,j$ such that $U_i\cap U_j\neq\vide$, and $w_j$ lying in the domain of the holonomy map $\h u{T_j}{T_i}$
$$\frac{d\left[\h u{T_i}{T_j}\,_\ast\nu^+_{T_i}\right]}{d\nu^+_{T_j}}(w_j)=\widehat{k}^u_F(w_j,\h u{T_j}{T_i}(w_j))=\frac{\widehat{k}^u_F(w_j,.)}{\widehat{k}^u_F(\h u{T_j}{T_i}(w_j),.)}.$$
Consequently, these measures can be glued and $\mu$ is well defined.

Now we have to see that $(\Pi_{\ast})_\ast\mu=\mu_F$. This is due to the first uniqueness result of Theorem \ref{familiesofmeasures}. Indeed, the fibration $\Pi_\ast$ commutes with unstable holonomies, so the projection of $(\nu^+_{T_i})_{i\in I}$ satisfies the cocycle relation \eqref{Eq:cocycleweakstable}. We know that in this case, this family is proportional to $(\lambda^{cs}_{|W^{cs}_{loc}(v_i)})_{i\in I}$ (after renormalization, we can assume it is equal). By definition of $\mu$, its projection on $T^1B$ is now locally defined by integration against $d\lambda^{cs}_{|W^{cs}_{loc}(v_i)}$ of the measures $k^u_{F}(v_i,v)\,d\lambda^u_{F,v_i}(v)$. This is the local product structure of $\mu_F$. We have what we wanted.

Finally, since by hypothesis the two measures $\mu$ and $\mu_F^+$ induce singular measures in a complete system of transversals to $\bcW^{u}$, they have to be singular. So we can conclude the proof.
\end{proof}

Even if the measure $\mu$ is not a priori invariant by the flow $G_t$, we will see that for any $t\geq 0$, $G_{-t}\ast\mu$ also satisfies the properties stated in Lemma \ref{propertymu}.

\begin{lemma}
\label{lemma1}
For any $t\geq 0$, there is a family of measures $(\nu_{-t,T}^+)_{T\in\TT^+}$ defined on the set $\TT^+$ of local transversals to $\bW^u$ such that:
\begin{enumerate}
\item $(\nu_{-t,T}^+)_{T\in\TT^+}$ satisfies Relation \eqref{Eq:cocycleweakstablelift} for any couple $T_1,T_2\in\TT^+$;
\item Locally, $G_{-t}\,_\ast\mu$ is obtained by integration against $d\nu_{-t,T}^+$ of the measures $\widehat{k}^u_{F}(w_T,w)d\widehat{\lambda}^u_{F,w_T}(w)$, $w_T\in T$.
In particular, $G_{-t}\,_\ast\mu$ has an absolutely continuous disintegration with respect to $(\bcW^u,\widehat{\lambda}^u_{F,w})$ and the local densities are uniformly log-bounded.
\end{enumerate} 
\end{lemma}

\begin{proof} First, choose $t\in[0,1]$. We can disintegrate $G_{-t}\,_\ast\mu$ in the plaques of $G_{-t}(U_i)$ with respect to the measure:
$$\exp\left[\int_0^t (F\circ G_s(w_i)-P(F))ds\right]\,d\left(G_{-t}\,_\ast\nu^+_{T_i}\right)(w_i)=d\nu^+_{-t,G_{-t}(T_i)}(w_i).$$

Note that the family of measures $(\nu_{-t,G_{-t}(T_i)}^+)_{i\in I}$ satisfies Relation \eqref{Eq:cocycleweakstablelift}.

Using the facts that the fibration $\Pi_{\ast}$ commutes with the flows, as well as with the unstable holonomy maps and that $\mu_F$ is invariant by $g_t$, we infer that $G_{-t}\,_\ast\mu$ projects down onto $\mu_F$ and that its disintegration with respect to $\nu_{-t,G_{-t}(T_i)}^+$ (which satisfies Relation \eqref{Eq:cocycleweakstablelift}) is given by lifting the local product structure of $\mu_F$.

It comes that for any $t\in[0,1]$, $G_{-t}\,_\ast\mu$ is obtained locally by integration against $d\nu_{-t,G_{-t}(T_i)}^+$ of the measures $\widehat{k}^u_{F}(w_i,w)\,d\widehat{\lambda}^u_{F,w_i}$. Since the transversals $(G_{-t}(T_i))_{i\in I}$ form a complete system, we can construct the desired family $(\nu_{-t,T})_{T\in\TT^+}$ by an immediate adaptation of Lemma 1.4 of \cite{BM}. Remark that the densities on the unstable plaques are log uniformly bounded independently of $t\in[0,1]$.

Now, proceed by induction on $n$ to prove that for any $t\in[n,n+1]$, such a family of measures exists. The heredity is straightforward: once we know the existence for $n$, we construct a family $(\nu^+_{-(n+t),G_{-t}(T_i)})_{i\in I}$ which makes the deal on the complete system of transversals $G_{-t}(T_i)$ for any $t\in[0,1]$. Then another adaptation of Lemma 1.4 of \cite{BM} provides the desired family $(\nu_{-(n+t),T})_{T\in T^+}$, thus concluding the proof.
\end{proof}

\begin{lemma}
\label{lemma2}
The measures $\mu_T$ converge to $\mu_F^-$ as $T$ tends to infinity, where:
$$\mu_T=\frac{1}{T}\int_0^T\left(G_{-t}\,_\ast\mu\right)\,dt.$$
\end{lemma}

\begin{proof}
The measure $\mu$ is singular with respect to $\mu_F^+$ and both of them project down to $\mu_F$ on $T^1B$. This means that their conditional measures on the fibers of $\Pi_{\ast}:T^1\F\to T^1B$ are singular. If $(\mu_v)_{v\in T^1B}$ denotes the disintegration of $\mu$ on the fibers, this implies that for $\mu_F$-almost every $v\in T^1B$, $\mu_v(\sigma^+(v))=0$. By Theorem \ref{twoergodic}, for  $\mu_F$-almost every $v\in T^1B$ and every point $w\in X_v\moins\{\sigma^+(v)\}$, $\dist_{X_v}(G_{-t}(w),G_{-t}(\sigma^-(v)))$ tends to $0$ as $t$ tends to $\infty$. Hence as $T$ goes to infinity, the past Birkhoff average of $w$ (i.e. $1/T\int_0^T\delta_{G_{-t}(w)}\,dt$) approaches that of $\sigma^-(v)$ and therefore must converge to $\mu_F^-$. Since $\mu_v(\sigma^+(v)))=0$ this property holds for $\mu_v$-almost every $w\in X_v$.

This proves the convergence of past Birkhoff averages of $\mu$-almost every $w\in T^1\F$ to $\mu_F^-$. It is now easy to conclude the proof by use of dominated convergence.
\end{proof}
\paragraph{End of the proof of Theorem \ref{uniqueergodicitiesstableunstable}.} The two lemmas above directly lead to a contradiction. Indeed, by Lemma \ref{lemma1}, all the measures $\mu_T$ defined in Lemma \ref{lemma2} have an absolutely continuous disintegration with respect to $(\bcW^u,\widehat{\lambda}^u_{F,v})$ and the local densities are uniformly log-bounded.

This implies that $\mu^-_F$, which is the limit measure, also has such an absolutely continuous disintegration, contradicting Proposition \ref{singular}. Hence, Theorem \ref{uniqueergodicitiesstableunstable} is proven.\quad \hfill $\square$.

\subsection{Disintegration of the $F$-harmonic measure}

We end this section by identifying the conditional measures in the fibers of $\Pi:M\to B$ of the unique $F$-harmonic measure for $\F$, denoted hereafter by $m_F$.

\paragraph{Structure of the $F$-harmonic measure.} For $z\in N$ define the homeomorphism $\check{\pi}_z:T^1_z N\mapsto N(\infty)$ sending every unit vector $v$ based at $z$ on $c_v(-\infty)$. It satisfies the invariance property: $\gamma\circ\check{\pi}_z=\check{\pi}_{\gamma z}\circ D\gamma$.

Consider the family $(\nu_z^F)_{z\in N}$ of Ledrappier measures as well as the family $(\omega_z^F)_{z\in N}$ of measures on $T^1_zN$ defined as $\omega_z^F=\check{\pi}^{-1}_z\,_{\ast}\nu_z^F$. These measures clearly satisfy $\omega_{\gamma z}^F=D_z\gamma\,_{\ast}\omega_z^F$ so they descend to a family $(\omega_p^F)_{p\in B}$ of measures on the $T^1_pB$.

The Lyapunov section lifts as a measurable section $\widetilde{\sigma}^+:T^1N\to T^1N\times\C\PP^1$  
satisfying $\widetilde{\sigma}^+\circ D\gamma=\alpha(\gamma)\circ\widetilde{\sigma}^+$ where $\alpha(\gamma)(v,x)=(D\gamma v,\rho(\gamma)x)$ 
denotes the diagonal action. It yields a measurable map $\widetilde{s}_z^+:T^1_zN\mapsto\{z\}\times\C\PP^1$ 
by $\widetilde{s}_z^+(v)=(z,pr_2(\widetilde{\sigma}^+(v))$ where $pr_2:T^1 N\times\C\PP^1\to\C\PP^1$ 
is the projection onto the second factor. It clearly satisfies 
$\widetilde{s}^+_{\gamma z}\circ D_z\gamma(v)=(\gamma z,\rho(\gamma)\circ\widetilde{s}^+_z(v))$. Hence the family $(\widetilde{s}^+_z)_{z\in N}$ descend as a family of measurable maps $s_p^+:T^1_pB\to X_p$, $p\in B$.

\begin{maintheorem}
\label{disintegrationfharmonic}
Let $(\Pi,M,B,\C\PP^1,\F)$ be a projective foliated bundle parametrized by a closed Riemannian manifold $B$ with negative sectional curvature. Assume that no probability measure on $\C\PP^1$ is invariant by the holonomy group. Let $F:T^1B\to\R$ be a Hölder continuous potential and denote by $m_F$ the unique $F$-harmonic measure. Consider the families $(\omega_p^F)_{p\in B}$ and $(s^+_p)_{p\in B}$ defined above and $h_0$ the unique $F$-harmonic function on $B$. Then
\begin{enumerate}
\item $m_F$ projects down to $h_0\,\Leb$;
\item the system of conditional measures of $m_F$ in the fibers $X_p$ is given for any $p$ by
$$m_{p,F}=s^+_p\,_\ast\left[\frac{\omega_p^F}{\mass(\omega_p^F)}\right].$$
\end{enumerate}
\end{maintheorem}

\paragraph{Trivialization of the center unstable foliation.} To prove this theorem, it is useful to work in $T^1N\times\C\PP^1$ with coordinates that trivialize the center unstable foliation. Consider the identification $\Phi:T^1N\to N\times N(\infty)$ such that $\Phi(z,v)=(z,\check{\pi}_z(v)))$. This identification sends $\tcW^{cu}$ on the partition $(N\times \{\xi\})_{\xi\in N(\infty)}.$

The lifted Lyapunov section commute with the center unstable foliation. As a consequence it reads in the coordinates given by $\Phi$ as $\widetilde{\sigma}^+(z,\xi)=(z,\xi,\widetilde{s}^+(\xi))$ where $\widetilde{s}^+:N(\infty)\to\C\PP^1$ is a measurable map satisfying $\widetilde{s}^+\circ \gamma=\rho(\gamma)\circ\widetilde{s}^+$. Given $z\in N$, define $\widetilde{\sigma}_z^+=\widetilde{\sigma}^+(z,.):N(\infty)\mapsto\{z\}\times N(\infty)\times\C\PP^1$.

\begin{lemma}
\label{conditionalmeasuinfharmonic}
Let $\widetilde{m}_F^+$ be the measure on $N\times N(\infty)\times\C\PP^1$ defined by integration of the measures $\widetilde{\sigma}^+_z\,_\ast\nu_z^F$ against $\Leb(z)$. Let $\widetilde{m}_F$ be the projection of $\widetilde{m}_F^+$ by the canonical projection $N\times N(\infty)\times\C\PP^1\to N\times\C\PP^1$.
\begin{enumerate}
\item $\widetilde{m}_F^+$ and $\widetilde{m}_F$ are respectively invariant by the diagonal actions of $\pi_1(B)$ on $N\times N(\infty)\times\C\PP^1$ and on $N\times\C\PP^1$;
\item the quotient measure of $\widetilde{m}_F$ on $M$ is the unique $F$-harmonic measure for $\F$;
\item $\widetilde{m}_F$ is obtained by integration of $\widetilde{s}^+_z\,_\ast\omega^F_z$ against $\Leb(z)$.
\end{enumerate}
\end{lemma}

\begin{proof}

Denote by $\beta$ the diagonal action on $N\times N(\infty)\times\C\PP^1$. By the equivariance properties of $\widetilde{\sigma}_z^+$ and $\nu_z^F$, one easily checks that the equivariance relation $\beta(\gamma)_{\ast}[\widetilde{\sigma}^+_z\,_{\ast}\nu^F_z]=\widetilde{\sigma}^+_{\gamma z}\,_{\ast}\nu^F_{\gamma z}$ holds for every $\gamma\in\pi_1(B)$.

Since $\pi_1(B)$ acts on $N$ by isometries we have $\gamma_{\ast}\Leb=\Leb$ for every $\gamma$. This proves that $\widetilde{m}_F^+$ is invariant by every $\beta(\gamma)$.

The invariance of $\widetilde{m}_F$ follows from the fact that if $P_1:N\times N(\infty)\to N$ denotes the projection on the first factor, we have $P_1\circ\gamma=\gamma\circ P_1$ (here $\pi_1(B)$ acts on $N\times N(\infty)$ by the diagonal action).

Now let us prove the second assertion. Let us fix a point $o\in N$. The measure $\widetilde{\sigma}_o^+\,_\ast\nu_o^F$ lives on $\{o\}\times N(\infty)\times\C\PP^1$, identified with $N(\infty)\times\C\PP^1$, which is supported by the graph of $\widetilde{s}^+$. We can disintegrate this measure in the $(N(\infty)\times\{x\})_{x\in\C\PP^1}$: $d[\widetilde{\sigma}_o^+\,_\ast\nu_o^F]=(d\eta_x)\,\,d[\widetilde{s}^+\,_{\ast}\nu_o^F](x)$, where $(\eta_x)_{x\in\C\PP^1}$ is a family of measures on the $N(\infty)\times\{x\}$.

The section has the special form $\widetilde{\sigma}^+(z,\xi)=(z,\xi,\widetilde{s}^+(\xi))$  and we have the cocycle relation $d\nu_z^F(\xi)=k^F(o,z;\xi)\,d\nu^F_o(\xi)$ so we have
\begin{eqnarray*}
d\widetilde{m}_F^+&=&d\left[\widetilde{\sigma}^+_z\,_{\ast}\nu_z\right](\xi,x)\,d\Leb(z)\\
                  &=&k^F(o,z;\xi)\,d\left[\widetilde{\sigma}_o^+\,_{\ast}\nu^F_o\right](\xi,x)\,d\Leb(z)\\
                  &=&\left(k^F(o,z;\xi)d\Leb(z)\,d\eta_x(\xi)\right)d\left[\widetilde{s}^+\,_\ast\nu_o^F\right](x).
\end{eqnarray*}

This proves that the densities with respect to Lebesgue of the conditional measures of $\widetilde{m}_F^+$ in the $N\times\{\xi\}\times\{x\}$ are given by $k^F(o,z;\xi)$. By projecting along the $N(\infty)$-factor we  see that the conditional measures of $\widetilde{m}_F$ in the $N\times\{x\}$ have a density with respect to Lebesgue given by integration of $k^F(o,z;\xi)$ against $\eta_x$. These are $F$-harmonic functions.

This proves that the quotient measure $m_F$ on $M$ has to be a $F$-harmonic measure. Since there is a unique one, this is $m_F$.

Now if one reads the measure $\widetilde{m}_F^+$ in $T^1N\times\C\PP^1$ (using $\Phi^{-1}$ in slices $N\times N(\infty)\times\{x\}$, $x\in\C\PP^1$) one sees that it is obtained as the integration against $\Leb(z)$ of measures $\widetilde{\sigma}^+_z\,_\ast\omega^F_z$. Denoting by $pr_z:T^1_zN\times\C\PP^1\to \{z\}\times\C\PP^1$ the basepoint projection, one has $\widetilde{s}^+_z=pr_z\circ\widetilde{\sigma}_z^+$, in such a way that $pr_z\,_\ast[\widetilde{\sigma}_z\,_\ast\,\omega^F_z]=\widetilde{s}^+_z\,_\ast\omega^F_z$. Hence the latter measures are the conditional measures of $\widetilde{m}_F$ on the fibers $\{z\}\times\C\PP^1$. The proof of the lemma is now over.
\end{proof}

\paragraph{End of the proof of Theorem \ref{disintegrationfharmonic}.} The unique $F$-harmonic measure $m_F$ is the quotient measure of $\widetilde{m}_F$. We have shown how the families $(\widetilde{s}^+_z)_{z\in N}$ and $(\omega_z^F)_{z\in N}$ induce maps $s_p^+:T^1_pB\to X_p\simeq\C\PP^1$, as well as measures $\omega^F_p$ on $T^1_pB$. Using Lemma \ref{conditionalmeasuinfharmonic} and descending to the quotient shows that the following disintegration holds $dm_F=d[s^+_p\,_\ast\omega^F_p]\,d\Leb(p).$

Now, remember that we defined the $F$-harmonic function $h_0(p)=\mass(\omega^F_p)$ on $B$. The measures $s^+_p\,_\ast\omega^F_p/h_0(p)$ are probability measures, so $\Pi_\ast m_F=h_0\Leb$ and
$$m_F=s^+_p\,_\ast\left[\frac{\omega_p^F}{\mass(\omega^F_p)}\right]\,h_0(p)\Leb(p).$$ \quad \hfill $\square$.

\section{Limits of large balls}
\label{slargeballs}

In this section we wish to prove that $F$-harmonic measures are limits of distributions of large balls \emph{weighted by a potential}. 

Before stating the precise result, we will need to set some notations.

\subsection{Weighted averages on large balls}
The leaf $L_x$ of a point $x\in M$ is a Riemannian cover of the base $B$ hence its Riemannian universal cover is given by $N$. We will denote the canonical projection $\proj_x:N\to L_x$.

\paragraph{The weight.} Let $F:T^1B\to\R$ be a Hölder continuous function and $\wF:T^1N\to\R$ be its lift. Denote by $S(o,R)\dans N$ the sphere centered at $o$ with radius $R$.

Consider the function defined on $N\times N$ by the following formula:

\begin{equation}
\label{Eq:Ffunction}
\kappa^F(o,z)=\exp\left[\int_{o}^z\wF\right],
\end{equation}

Let us give a theorem which is a weighted version of a theorem due to Margulis. It is due to Ledrappier \cite{L4}, in a slightly more general form.

\begin{theorem}
\label{kappaspheres}
For every $o\in N$, there is a number $c(o)$, such that:
$$\int_{S(o,R)}\kappa^F(o,y)d\Leb(y)\sim c(o)e^{RP(F)},$$
as $R$ goes to infinity.
\end{theorem}

\paragraph{Weighted spherical averages.} First we define the \emph{spherical averages}. If $x\in L_x$ and $R>0$ we choose $o\in\proj_x^{-1}(x)$.

\begin{equation}
\label{Eq:sphericalaverages}
m^F_{x,R}=\proj_x\,_\ast\left(\frac{\kappa^F(o,y)\,d\Leb_{|S(o,R)}(y)}{\int_{S(o,R)}\kappa^F(o,y) d\Leb(y)}\right).
\end{equation}
Note that this measure is independent of the choice of $o\in\proj_x^{-1}(x)$. The main result of this section is the

\begin{theorem}
\label{limitoflargespheres}
Let $(\Pi,M,B,\C\PP^1,\F)$ be a projective foliated bundle parametrized by a closed and negatively curved base $B$. Assume that there is no probability measure on $\C\PP^1$ which is invariant by the holonomy group. Then, for every Hölder continuous potential $F:T^1B\to\R$ and every sequences $(x_n)_{n\in\N}\in M^{\N}$ and $(R_n)_{n\in\N}$ tending to infinity, the sequence of measures $m^F_{x_n,R_n}$ defined above converges to the unique $F$-harmonic measure.
\end{theorem}

\paragraph{Limit of large balls.} We will now consider the weighted averages on large balls $B(o,R)$. For $o\in\proj_x^{-1}(x)$ and $R>0$ set

\begin{equation}
\label{bigballs}
\mu^F_{x,R}=\proj_x\,_\ast\left(\frac{\kappa^F(o,y)\,d\Leb_{|B(o,R)}(y)}{\int_{B(o,R)}\kappa^F(o,y) d\Leb(y)}\right).
\end{equation}

We see these measures as multidimensional and weighted analogues of Birkhoff averages. Here again they are independent of the choice of $o\in\proj_x^{-1}(x)$. One of the main result of the papers is
\begin{maintheorem}
\label{limitlargeballs}
Let $(\Pi,M,B,\C\PP^1,\F)$ be a projective foliated bundle parametrized by a closed Riemannian manifold $B$ with negative sectional curvature. Assume that no probability measure on $\C\PP^1$ is invariant by the holonomy group. Let $F:T^1B\to\R$ be a Hölder continuous potential with $P(F)>0$. Then, for every sequences $(x_n)_{n\in\N}\in M^{\N}$ and $(R_n)_{n\in\N}$ tending to infinity, the sequence of measures $\mu^F_{x_n,R_n}$ defined above converges to the unique $F$-harmonic measure.
\end{maintheorem}

\paragraph{Proof of Theorem \ref{limitlargeballs}.} Let us emphasize on the hypothesis on the pressure of $F$ in order to deal with weighted averages on large balls. Because of Theorem \ref{kappaspheres}, the integral $\int_{S(o,R)}\kappa^F(o,y) d\Leb(y)$ behaves like $e^{RP(F)}$: in particular, it tends to infinity only when $P(F)>0$. This property is satisfied in the case where the potential was the null function, whose pressure is equal to the topological entropy of $g_t$.

The argument follows the lines of that of Bonatti and G\'omez-Mont (see \cite[Proposition 0.1]{BG}): let us sketch the proof. Assume that no measure on $\C\PP^1$ is invariant by the holonomy group. By Theorem \ref{limitoflargespheres}, the weighted spherical averages converge to the unique $F$-harmonic measure for $\F$. Take the crown $C_R=B(o,R)\moins B(o,R/2)$. Since the pressure of $F$ is positive, we know that the integral of the weight on $C_R$ becomes closer and closer to that on $B(o,R)$. Hence, it is sufficient to show that the normalized restriction of $\mu^F_{x,R}$ to the projection of $C_R$ approaches the unique $F$-harmonic measure.

Now, the latter measure is an average of weighted spherical averages, each of them approaching the $F$-harmonic measure (by Theorem \ref{limitoflargespheres}). We conclude that weighted averages on large balls also converge to the same limit. \quad \hfill $\square$.

\paragraph{Remark.} Let us make a comment on our pressure assumption. First note that $P(F+c)=P(F)+c$ and that $k^F=k^{F+c}$. This proves that every $F$-harmonic measure and $F+c$-harmonic measures coincide. As a result, every $F$-harmonic measure is the limit of weighted averages on large balls .

Without the assumption $P(F)>0$, we can priori give an answer only for the convergence of weighted averages on large crowns $B(o,R)\moins B(o,\phi(R))$ where $\phi$ is a function tends to infinity.

For $P(F)<0$, the integral $\int_{N}\kappa^F(o,y)d\Leb(y)$ converges and the question about weighted averages on balls is not very interesting. However when $P(F)=0$ we can use Patterson's trick (see \cite{Pat}) to ensure the convergence. It consists in considering a continuous function $h:\R^+\to\R^+$ which is non-decreasing, has slow growth (i.e. for every $d>0$, $|h(R+d)/h(R)-1|<\eps$ when $R$ large) and satisfies 
$$\lim_{R\to\infty}h(R)\int_{S(o,R)}\kappa^F(0,y)d\Leb(y)=\infty.$$
Since $\int_{S(0,R)}\kappa^F(0,y)d\Leb(y)$ has a nonzero limit when $R$ goes to infinity such an $h$ is very easy to construct ($h(R)=\log(1+R)$). The proof of theorem \ref{limitlargeballs} adapts to give convergence of weighted averages of large balls for the new weight $h(\dist(o,y))\kappa^F(o,y)$ when $P(F)=0$.

\subsection{The horospherical foliation as a limit}
\label{shorosphericlimit}

We will now show how to prove Theorem \ref{limitoflargespheres}. The idea of the proof is a weighted and higher dimensional version of \cite[Proposition 0.1]{BG}. We will use that large spheres look like horospheres and that any accumulation point of the family of measures $m^F_{x,R}$ induces a family of quasi-invariant measures transverse to the unstable foliations satisfying Relation \eqref{Eq:cocycleweakunstablelift}. We will then use the uniqueness of this family.

Such a line of reasoning is not new. It has been used by Knieper in \cite{Kn} in order to prove the convergence of spherical means to the horospherical measure. We also cite \cite{AR,K2,KP,Sc1,Sc2} for works using this idea.

\paragraph{Horospheres and large spheres.} For $z\in N$ and $R>0$, the sphere $S(z,R)$ can be embedded in $T^1N$ by attaching to each point $y$ of the sphere, the outward normal vector based at $y$. We denote this embedded sphere by $S^+(z,R)$. Such spheres form a foliation of $T^1N$ denoted by $\widetilde{\W}^+_R$ which satisfies $G_t(\widetilde{\W}^+_R)=\widetilde{\W}^+_{R+t}$. The classical inclination lemma for compact hyperbolic sets (see chapter 9 of \cite{Sh}) implies that, if $\W^+_R$ denotes the projection of $\widetilde{\W}^+_R$ to $T^1B$, $\W^+_R$ converges to $\W^u$ in the $C^0$-topology of plane fields (the convergence is uniform thanks to the compactness of $T^1B$). This implies that the induced leafwise metric on $\W^+_R$ converges uniformly to that of $\W^u$. It also implies the following

\begin{proposition}
\label{cvfoliations}
For every small $\eps>0$, there is a positive number $R_0$ and a finite family of disjoint small embedded discs $T_i\dans T^1B$, $i\in I$ such that for every $R>R_0$
\begin{enumerate}
\item $(T_i)_{i\in I}$ is both a complete system of transversal for $\W^u$ and for $\W^+_R$;
\item for each holonomy map $\h u{S_i}{S_j}:S_i\to S_j$ along a path included in $W^u(v)$ of length $\leq 1$, where $S_i\dans T_i$ and $S_j\dans T_j$ are relatively compact open sets with $v\in S_i$, there are open sets $S_i'\dans S_i$, $S_j'\dans S_j$ with $v\in S_i'$, a path $c$ in $W^+_R(v)$ of length $<2$ and a holonomy map for $\W^+_R$ along $c$, $\tau_{R,c}:S_i'\to S_j'$ which is $\eps$-close in the $C^0$-topology to the restriction to $S_i'$ of holonomy map $\h u{S_i}{S_j}$.
\end{enumerate}
\end{proposition}

We can lift the foliations $\W^+_R$ to the leaves of $\widehat{\F}$ as subfoliations that we denote by $\bcW^+_R$. They converge to $\bcW^u$ in the topology of plane fields: and the proposition above holds in this context as well.

\paragraph{A geometric estimate.}
We will need the following distortion lemma, which is proven by using the CAT inequalities.
\begin{lemma}
\label{distortioncomparison}
There are positive constants $C>0$ and $\alpha>0$, such that for all $o\in N$, $R>0$ and $y,z\in S(o,R)$, we have
$$\frac{\kappa^F(o,z)}{\kappa^F(o,y)}\leq\exp\left(C\dist_{S(o,R)}(y,z)^{\alpha}\right),$$
where $\dist_{S(o,R)}$ is the distance function coming from the Riemannian structure induced on $S(o,R)$.
\end{lemma}

\begin{proof}
First note that the following formula holds
\begin{equation}
\label{Eq:relationFfunction}
\frac{\kappa^F(o,z)}{\kappa^F(o,y)}=\exp\left[\int_o^z\widetilde{F}-\int_o^y\widetilde{F}\right].
\end{equation}

Let $o\in N$, $R>0$ and $y,z\in S(o,R)$. We parametrize the geodesic segments $[o,y]$ and $[o,z]$ by arc length: these parametrizations are respectively denoted by $y(r)$ and $z(r)$. Let $c(r)$ be a minimizing geodesic in $S(o,r)$ between $y(r)$ and $z(r)$ and $l(r)$ be its length. We state without proof the following claim which follows classically from a use of CAT($-a^2$) inequalities.

\paragraph{Claim.} \emph{For every} $r\leq R$, \emph{the following holds true}
$$\frac{l(r)}{l(R)}\leq\frac{\sinh(ar)}{\sinh(aR)}.$$

Let us prove that the claim implies the lemma. We consider the unit vectors $v_{oy}$ and $v_{oz}$ based respectively at $y$ and $z$, which are normal to the sphere and point outside. Since the transport on $S(o,R)$ for the induced connexion of the orthogonal vector field is parallel along $c(R)$, the Sasaki distance between these two vectors is equal to $l(R)$ (this argument is obviously valid for any choice of $r\leq R$). The lemma then follows from the Hölder continuity of $F$.
\end{proof}

\paragraph{Uniform convergence to the cocycle.} Now let us show how to approach uniformly on compact sets the Gibbs kernel $k^F$, defined by Formula \eqref{gibbskernel}, thanks to the weight $\kappa^F$.

\begin{lemma}
\label{convergencetococycle}
Let $K\dans N$ be a compact set and $\xi\in N(\infty)$. Then
$$\lim_{o\to\xi}\frac{\kappa^F(o,z)}{\kappa^F(o,y)}=k^{F}(y,z;\xi),$$
uniformly in $y,z\in K$ lying in the same horosphere centered at $\xi$, as $o$ converges nontangentially to $\xi$ in the cone topology.
\end{lemma}

\begin{proof}
First recall that the weight $\kappa^F$ satisfies Relation \eqref{Eq:relationFfunction} and that for $y,z$ on the same horosphere at $\xi$, $\beta_{\xi}(y,z)=0$.

Now, let $\xi\in N(\infty)$, $K\dans N$ compact and $y,z\in K$ lying on the same horosphere centered at $\xi$. It is enough, by definition of the nontangential limit, to assume that $o$ tends to $\xi$ by staying on some geodesic ray $c$ with $\xi$ as an extremal point. Note that the difference

$$\left(\int_{c(T)}^z\widetilde{F}-\int_{c(T)}^y\widetilde{F}\right)-\left(\int_{\xi}^z\widetilde{F}-\int_{\xi}^y\widetilde{F}\right)$$
can be broken into four pieces
$$\left(\int_{c(T)}^z\widetilde{F}-\int_{z(T)}^z\widetilde{F}\right)+\left(\int_{\xi}^{c(T)}\widetilde{F}-\int_{\xi}^{z(T)}\widetilde{F}\right)+\left(\int_{\xi}^{y(T)}\widetilde{F}-\int_{\xi}^{c(T)}\widetilde{F}\right)+\left(\int_{y(T)}^y\widetilde{F}-\int_{c(T)}^y\widetilde{F}\right).$$

Using the classical distortion control, as well as Lemma \ref{distortioncomparison}, we see that each of these terms can be controlled by quantities tending exponentially fast to $0$ independently of the choice of $y,z$ \emph{in the compact} $K$.

We finish the proof by taking the exponential of these quantities.
\end{proof}

\subsection{Proof of Theorem \ref{limitoflargespheres}}
\label{slimitsmeasureslargeballs}

We assume that there is no probability measure on $\C\PP^1$ invariant by the holonomy group of the projective foliated bundle $(\Pi,M,B,\C\PP^1,\F)$.

\paragraph{Lift to the unit tangent bundle.} The measures $m^F_{x,R}$ can obviously be lifted to $T^1\F$ through the embedding $S(x,R)\to S^+(x,R)$ (remember that $S^+(x,R)$ denotes the image of $T^1_x\F$ by $G_R$). We will denote by  $m^{F,+}_{x,R}$ the resulting measure.

Since the foliation $\bcW^+_R$ converges to $\bcW^u$ uniformly when $R$ goes to infinity (see Proposition \ref{cvfoliations}), we can find $R_0>0$ and a foliated atlas $\A=(U_i,\phi_i)_{i\in I}$ for the strong unstable foliation $\bcW^u$, such that $\A_R=(U_i,\phi_{R,i})_{i\in I}$ is a foliated atlas for $\bcW^+_R$, $R>R_0$. Denote by $(T_i)_{i\in I}$ a corresponding complete system of transversals and by $(P_i(w))_{w\in T_i}$ the unstable plaques. We intend to prove the following proposition.

\begin{proposition}
\label{accumulationharmonic}
Let $m^+$ be an accumulation point of the family $m^{F,+}_{x,R}$. Its restriction to a chart $U_i$ has a disintegration in the unstable plaques of the form $[H_i(v)\Leb^u_{|P_i(w)}(v)]\,\nu_i^+(w)$, where:
\begin{enumerate}
\item $H_i$ is a positive and measurable function which is uniformly log-bounded and which satisfies the following formula for $v,w$ on the same plaque $P_i$
\begin{equation}
\label{Eq:densitieslargespheres}
\frac{H_i(v)}{H_i(w)}=k^F(w,v;\xi),
\end{equation}
where $\xi\in N(\infty)$ is the common limit of the backward iterations on the unstable manifold of $v$ and $w$.
\item $(\nu_i^+)_{i\in I}$ is a family of finite measures defined on the family of transversals $(T_i)_{i\in I}$ satisfying
\begin{equation}
\label{Eq:transvmeaslargespheres}
\frac{d[\h u{T_i}{T_j}\ast\nu_i^+]}{d\nu_j^+}(v)=k^F(v,\h u{T_j}{T_i}(v);\xi),
\end{equation}
where $v$ belongs to the domain of $\h u{T_j}{T_i}$ and $\xi$ is the common limit of the backward iterations on the unstable manifold of $v$.
\end{enumerate}
\end{proposition}

\paragraph{Remark 1.}
The notation used here is a little bit abusive because $v,w$ lie in $T^1 P_i$ and $k^F$ has been defined on $N\times N\times N(\infty)$: we have evaluated it on lifts to the universal cover of the base points of the vectors $v$ and $w$ and on the common limit point in the past.
\paragraph{Remark 2.} The two properties are equivalent in the sense that if one can glue together the measures $[H_i(v)\Leb^u_{|W^u_{loc}(w)(v)}]\,\nu_i^+(w)$ and form a measure of $T^1\F$, then if one of the two properties holds, the other holds too.

\paragraph{Proof of Theorem \ref{limitoflargespheres} from Proposition \ref{accumulationharmonic}.} Assume that the proposition above is true and take an accumulation point $m^+$. In that case, the family $(\nu^+_i)_{i\in I}$ satisfies the cocycle relation \eqref{Eq:cocycleweakstablelift} and by the unique ergodicity properties of the unstable foliation it is the family described in Theorem \ref{uniqueergodicitiesstableunstable}. By the results of the previous section, $m^+$ has to be the canonical lift of the unique $F$-harmonic measure: if one prefers, its projection on $M$ is the unique $F$-harmonic measure. This ends the proof of Theorem \ref{limitoflargespheres}. \quad \hfill $\square$

Consequently, in order to conclude the proof of this theorem, Proposition \ref{accumulationharmonic} only remains to be proven.

\subsection{Proof of Proposition \ref{accumulationharmonic}.}
Let $m^+$ be an accumulation point of $m^{F,+}_{x,R}$. We choose sequences $(x_n)_{n\in\N}\in M^{\N}$, $(R_n)_{n\in\N}$ which tends to infinity, such that $m^{F,+}_{R_n,x_n}$ converges to $m^+$. Remember that we have a foliated atlas $\A=(U_i,\phi_i)_{i\in I}$ for the strong unstable foliation $\bcW^u$ and a number $R_0>0$ such that when $R>R_0$, $\A_R=(U_i,\phi_{R,i})_{i\in I}$ is a foliated atlas for $\bcW^+_R$. We can refine the atlas $\A$ in such a way that $m^+(\partial U_i)=0$ for any $i\in I$ (so, we always have $\lim_{n\to\infty}m^{F,+}_{x_n,R_n}(U_i)=m^+(U_i)$). The corresponding complete system of transversals is denoted by $(T_i)_{i\in I}$. The corresponding unstable plaques are denoted by $(P_i(v))_{v\in T_i}$ and the plaques of the foliation $\bcW^+_{R_n}$ are denoted by $(P_{i,n}(v))_{v\in T_i}$. Let $i\in I$ be such that $m^+(U_i)>0$.

\begin{lemma}
\label{transversallargesphere}
The projection on the transversal $T_i$ of the restriction of $m^{F,+}_{x_n,R_n}$ to $U_i$ is equivalent to the following measure
$$\nu^+_{i,n}=\frac{1}{\int_{S(o,R_n)}\kappa^F(o,y) d\Leb(y)}\sum_{v\in T_i\cap S^+(x_n,R_n)}\kappa^F(x_n,v)\delta_v$$
with a Radon-Nikodym derivative which is log-bounded independently of $n$.
\end{lemma}

\paragraph{Remark 3.} Here again, the notation is a bit abusive. We have to lift $v\in T_i$ and $x_n$ to the universal cover and to evaluate $\kappa^F$ on the lift of $x_n$ and the base point of the lift of $v$. If we impose that this base point stays in a given fundamental domain, everything happens as if $x_n$ tends to infinity in the past while staying on the geodesic directed by $v$.

\begin{proof}
The projection of $(m^{F,+}_{x_n,R_n})_{|U_i}$ on $T_i$ is given by the following counting measure
$$ \frac{1}{\int_{S(o,R_n)}\kappa^F(o,y) d\Leb(y)}\sum_{v\in T_i\cap S^+(x_n,R_n)}m^{F,+}_{x_n,R_n}(P_{i,n}(v))\delta_v.$$

Let $\Delta$ be an upper bound of the diameters of all plaques of $\bcW^+_{R_n}$ and $K^{\pm 1}$ be numbers which bound their volumes respectively from above and below. Using distortion lemma \ref{distortioncomparison}, we see that for any $v\in T_i\cap S^+(x_n,R_n)$, $m^{F,+}_{x_n,R_n}(P_{i,n}(v))\in [C_0^{-1},C_0]\kappa^F(x_n,v)$, where $C_0=Ke^{C\Delta^{\alpha}}$.

Of course, $\Delta$ and $V$ can be chosen independently of $n$ since the foliation $\bcW^+_R$ converges uniformly to $\bcW^u$: in particular bounds on diameter and volumes of plaques of $\bcW^u$ give bounds on those of plaques of $\bcW^+_{R_n}$, at least for $n$ large enough. Hence the lemma is proven.
\end{proof}

\begin{lemma}
\label{disintegrationlargespheres}
The restriction of $m^{F,+}_{x_n,R_n}$ to a chart $U_i$ has a disintegration in the unstable plaques of the form $[H_{i,n}(v)\Leb_{|P_{i,n}(w)(v)}]\,\nu_{i,n}^+(w)$, where
\begin{enumerate}
\label{approachdensitieslargespheres}
\item $H_{i,n}$ is a positive and measurable function which is uniformly log-bounded independently of $n$ and which satisfy the following formula for $v,w$ on the same plaque $P_{i,n}$
\begin{equation}
\frac{H_{i,n}(v)}{H_{i,n}(w)}=\frac{\kappa^F(x_n,v)}{\kappa^F(x_n,w)},
\end{equation}
where  $v$ and $w$ are normal to the sphere of center $x_n$ and of radius $R$.
\item $(\nu_{i,n}^+)_{i\in I}$ is the family of finite measures defined on the family of transversals $(T_i)_{i\in I}$ by Lemma \ref{transversallargesphere}. They satisfy
\begin{equation}
\label{Eq:approachtransvmeaslargespheres}
\frac{d[\tau_{ij,n}\ast\nu^+_{i,n}]}{d\nu^+_{j,n}}(v)=\frac{\kappa^F(x_n,\tau_{ji,n}(v))}{\kappa^F(x_n,v)},
\end{equation}
where $\tau_{ij,n}$ is a holonomy map of $\W^+_{R_n}$ whose domain is an open subset of $T_i$ and its range, an open subset of $T_j$ and where $v$ belongs to the domain of $\tau_{ji,n}$.
\end{enumerate}
\end{lemma}

\begin{proof}
Since by definition, $m^{F,+}_{x_n,R_n}$ has Lebesgue disintegration in the leaves of $\bcW^+_{R_n}$, the disintegration of $(m^{F,+}_{x_n,R_n})_{|U_i}$ with respect to $\nu_{i,n}^+$ has the same form as stated in the lemma. We have to see that the densities $H_{i,n}$ and the measures $\nu_{i,n}$ satisfy the stated properties.

The second property is immediate by the definition of $(\nu_{i,n}^+)_{i\in I}$: see Lemma \ref{transversallargesphere}.

The fact that densities $H_{i,n}$ are log-bounded independently of $n$ follows from the fact that $\nu_{i,n}^+$ is equivalent to the projection on $T_i$ of $(m^{F,+}_{x_n,R_n})_{|U_i}$ with a Radon-Nikodym derivative which is log-bounded independently of $n$ (see Lemma \ref{transversallargesphere}). In order to see that these densities satisfy the desired relation, one just has to use the fact that this is the case for the family $(\nu_{i,n}^+)_{i\in I}$ and that the integrated measures on $U_i$ can be glued together nicely (see Remark 2 above). The proof is therefore over.
\end{proof}

\paragraph{End of the proof of Proposition \ref{accumulationharmonic}}
By the choice of charts $U_i$, we know that $(m^{F,+}_{x_n,R_n})_{|U_i}$ converges to $m^+_{|U_i}$ as $n$ tends to infinity. We want to deduce that $m^+$ disintegrates in the plaques $P_i$ as desired. Take two charts $U_i$ and $U_j$ with non empty intersection. We can always assume that when $v$ lies in the domain of a holonomy map $\h u{T_i}{T_j}$, we have $\dist_u(v,\h u{T_i}{T_j}(v))\leq 1$.

First, one can assume, by taking subsequences if needed, that $\nu_{n,i}^+$ converges to some finite measure $\nu^+_i$ on $T_i$. We know from Proposition \ref{cvfoliations} that for any $v$ lying in the domain of $\h u{T_j}{T_i}$, there exists an open set $S_j\dans T_j$ containing $v$, such that the restriction of $\h u{T_j}{T_i}$ to $S_i$ is a uniform limit of holonomy maps $\tau_{ij,n}$. Now, because of Lemma \ref{convergencetococycle}, we have for any $v\in S_j$
$$\lim_{n\to\infty}\frac{\kappa^F(x_n,\tau_{ji,n}(v))}{\kappa^F(x_n,v)}=k^F(v,\h u{T_j}{T_i}(v);\xi),$$
where $\xi\in N(\infty)$ is the limit of backward iterations of $W^u(v)$ by the flow. Therefore, since, by Lemma \ref{disintegrationlargespheres}, the measures $\nu_{n,i}^+$ satisfy Relation \eqref{Eq:approachtransvmeaslargespheres}, the limit measures $\nu^+_i$ satisfy Relation \eqref{Eq:transvmeaslargespheres}.

In order to finish the proof, we have to show that $m^+$ has a Lebesgue disintegration. This is true indeed because on the one hand, each $m^{F,+}_{x_n,R_n}$ has Lebesgue disintegration with local densities log bounded independently of $n$ and on the other hand, the Riemannian structure and hence the Lebesgue measure, induced on the leaves of $\bcW^+_{R_n}$ converges to that of $\bcW^u$. Finally, in order to see that the densities $H_i$ satisfy the desired relation, we can proceed as in the end of the proof of Lemma \ref{approachdensitieslargespheres}. The proof is now over.\quad \hfill $\square$.

\section{An equidistribution result}
\label{sequidistribution}

In this section, we will give an interpretation of our results foliations transverse to $\C\PP^1$-bundles from the point of view of group actions parametrized with a Riemannian metric of negative curvature. In what follows, we consider a projective representation $\rho:\pi_1(B)\to \PSL_2(\C)$ which leaves no measure invariant on $\C\PP^1$.

Given $o\in N$ we get a distance function on $\pi_1(B)$ by setting $d(\gamma_1,\gamma_2)=\dist(\gamma_1 o,\gamma_2 o)$. Remark that since $B$ is compact, this distance is quasi-isometric to the word distance associated to any system of generators. Let $B_R$ denote the corresponding ball of radius $R>0$. Other links between continuous and discrete equidistribution results may be found in \cite{Rob,Sc1}.

\subsection{Weighted counting measures}
\label{sweightedcm}
Under our hypothesis, Theorem \ref{tripleuniquegibbs} states for every $F:T^1B\to\R$ Hölder continuous, there is a unique $F$-harmonic measure for the suspended foliation $\F$, denoted by $m_F$. By Theorem \ref{disintegrationfharmonic}, this measure can be disintegrated with respect to $h_0\Leb$ in the base and we denote by $(m_{F,p})_{p\in B}$ the system of conditional measures on the fibers of $\Pi:M\to B$. Remember that for any $p\in B$, $h_0(p)=\mass(\omega_p^F)$, where $\omega_p^F$ stands for the Ledrappier measure on $T^1_pB$ associated to the potential $F$. Remember also that $m_{F,p}=s^+_p\,_\ast\omega_p^F/\mass(\omega_p^F)$.

Here we want to show that these measures, which are defined on $\C\PP^1$, can be obtained as limits of weighted counting measures. Let $o\in N$. Give the element $\gamma\in\pi_1(B)$ the following weight
$$\kappa^F_o(\gamma)=\kappa^F(o,\gamma o).$$

Note that $\kappa^F_o=\kappa^F_{\eta o}$ for every $\eta\in\pi_1(B)$.

\paragraph{Limits of weighted counting measures.} Fix $p\in B$, $x\in X_p$ and $o\in\proj_x^{-1}(x)$. We are interested in the \emph{weighted counting measures} defined by
$$\theta_{F,R}=\frac{1}{\sum_{\gamma\in B_R}\kappa^F_o(\gamma)}\sum_{\gamma\in B_R}\kappa_o^F(\gamma) \delta_{\rho(\gamma)^{-1} x}.$$

The main goal of this section is to prove the following

\begin{maintheorem}
\label{countingmeasures}
Let $B$ be a closed Riemannian manifold with negative sectional curvature. Consider a projective representation $\rho:\pi_1(B)\to \PSL_2(\C)$ which leaves no probability measure invariant on $\C\PP^1$ and consider a Hölder continuous potential $F:T^1B\to\R$. Assume moreover that the potential $F$ has positive pressure. Then the measure $\theta_{F,R}$ converges to $m_{F,p}$ as $R$ tends to infinity.
\end{maintheorem}

\paragraph{Remark.} One can see $\theta_{F,R}$ in a more geometric way. Consider the counting measure on $N$ defined as
$$\tilde{\theta}_{F,R}=\frac{1}{\sum_{\gamma\in B_R}\kappa^F_o(\gamma)}\sum_{\gamma\in B_R}\kappa^F_o(\gamma) \delta_{\gamma o}.$$
The measure $\theta_{F,R}$ is then defined as $\theta_{F,R}=\proj_x\ast\tilde{\theta}_{F,R}$.

We have the following interpretation of $\theta_{F,R}$. Consider a large ball centered at $x$ tangent to its leaf. If it intersects the fiber $X_p$ at a point $y$, give to this point the weight $\kappa^F(x,y)$ (if the ball intersects the fiber at some point multiple times, its weight will be the sum of the corresponding $\kappa^F(x,y)$). The measure $\theta_{F,R}$ is the average of the Dirac masses at the intersection between the ball and the fiber, weighted by these $\kappa^F$.

\subsection{Proof of Theorem \ref{countingmeasures}}

\paragraph{Restriction to a small cylinder.} Under the hypothesis of Theorem \ref{countingmeasures}, Theorem \ref{limitlargeballs} states that the unique $F$-harmonic measure $m_F$ is obtained as the limit of weighted averages of larges balls inside the leaves. The idea is quite simple: inside a small neighbourhood of the fiber, the restriction of weighted averages in the balls has to converge to the restriction of $m_F$. If we manage to compare this normalized restriction with $\theta_{F,R}$, we should be able, by letting the size of the neighbourhood tend to zero, to see the convergence to the conditional measure of $m_F$ on the fiber.

Choose $\eps>0$ small enough so that the ball $B(p,\eps)$ trivializes the fiber bundle, as well as the cover $N\to B$. Consider the cylinder $K_{p,\eps}=\Pi^{-1}(B(p,\eps))\simeq B(p,\eps)\times\C\PP^1$. It comes with a natural projection map along the leaves $pr_{\eps}:K_{p,\eps}\to X_p$. The boundary of $K_{p,\eps}$ has $m_F$-measure $0$ since projection of $m_F$ on the base is equivalent to the volume and the projection of $\partial K_{p,\eps}$ is the sphere $S(p,\eps)$, which has zero volume. As a consequence, we obtain the following
\begin{lemma}
\label{convergenceinrestriction}
In restriction to $K_{p,\eps}$, the measure $\mu^F_{x,R}$ defined by \eqref{bigballs} converges to $(m_F)_{|K_{p,\eps}}$ as $R$ tends to infinity.
\end{lemma}

The intersection of $B(x,R)$ with the small cylinder is not markovian. Some connected components of the intersection are strictly included in the corresponding connected component of $L_x\cap K_{p,\eps}$. As we have chosen a potential with positive pressure, these non-markovian components can't be neglected. Thus, we will be led to introduce two ``markovian'' measures.

\paragraph{Notation.} We will use the two following convenient notations for two quantities that we shall compare later
\begin{equation}
\label{Eq:integralweight}
I_{F,R}=\int_{B(o,R)}\kappa^F(o,y) d\Leb(y),
\end{equation}

\begin{equation}
\label{Eq:sumweight}
J_{F,R}=\sum_{\gamma\in B_R}\kappa^F_o(\gamma).
\end{equation}

\paragraph{Markovian measures.} Denote by $\mu^F_{x,R,\eps}$ the restriction of $\mu^F_{x,R}$ to $K_{p,\eps}$. Note that it can be obtained as the projection on $M$ by the canonical projection $\proj_x:(N,o)\to (L_x,x)$ of the following measure
\begin{equation}
\label{Eq:mufr}
\tilde{\mu}^F_R=\frac{1}{I_{F,R}}\sum_{\gamma\in B_R}\kappa^F(o,y)\Leb_{|B(o,R)\cap B(\gamma o,\eps)}(y).
\end{equation}
Now consider the two measures $\mu^{F\,\pm}_{x,R,\eps}$ defined as $\proj_x\,_\ast(\tilde{\mu}^{F\,\pm}_R)$, where

\begin{equation}
\label{Eq:mufrpm}
\tilde{\mu}^{F\,\pm}_R=\frac{1}{I_{F,R}}\sum_{\gamma\in B_{R,\eps}^{\pm}}\kappa^F(o,y)\Leb_{|B(\gamma o,\eps)}(y),
\end{equation}
where
$$B^-_{R,\eps}=\{\gamma\in\pi_1(B)\,|\,B(\gamma o,\eps)\dans B(o,R)\},\,\,\,\,\,\textrm{and},\,\,\,\,\,B^+_{R,\eps}=\{\gamma\in\pi_1(B)\,|\,B(\gamma o,\eps)\cap B(o,R)\neq\vide\}.$$
In other words, $\mu^{F\,-}_{x,R,\eps}$ only charges all connected components of $B(x,R)\cap K_{p,\eps}$ which are entirely included in $K_{p,\eps}$, as $\mu^{F\,+}_{x,R,\eps}$ charges all the small discs that meet a connected component of $B(x,R)\cap K_{p,\eps}$. The inequalities $\mu^{F\,-}_{x,R,\eps}\leq\mu^F_{x,R,\eps}\leq\mu^{F\,+}_{x,R,\eps}$ are then obvious. We have better.
\begin{lemma}
\label{chainofinequalities}
We have the following chain of inequalities
$$\frac{I_{F,R-2\eps}}{I_{F,R}}\mu^F_{x,R-2\eps,\eps}\leq\mu^{F\,-}_{x,R,\eps}\leq\mu^F_{x,R,\eps}\leq\mu^{F\,+}_{x,R,\eps}\leq\frac{I_{F,R+2\eps}}{I_{F,R}}\mu^F_{x,R+2\eps,\eps}.$$
\end{lemma}

\begin{proof}
As mentioned before, the second and third inequalities are obvious. Let us prove the first one, the last one can be proven by the same kind of arguments. It is enough to prove that $B(x,R-2\eps)\cap K_{p,\eps}$ is included in the union of those connected components of $B(x,R)\cap K_{p,\eps}$ which cross the cylinder.

 But this is a rather immediate consequence of the triangular inequality. If $y\in B(x,R-2\eps)\cap K_{p,\eps}$, we project $y$ on $V_p$ by considering $\tilde{p}=pr_{\eps}(y)$: in particular, $\dist_{L_x}(y,\tilde{p})\leq\eps$. If one picks $y'\in B(\tilde{p},\eps)$, one has
$$\dist_{L_x}(y',x)\leq\dist_{L_x}(y',\tilde{p})+\dist_{L_x}(\tilde{p},y)+\dist_{L_x}(y,x)\leq\eps+\eps+R-2\eps=R.$$
Finally, the ball $B(\tilde{p},\eps)$ is entirely included in $B(x,R)$: this allows us to conclude the proof of the lemma.
\end{proof}

\paragraph{Projections on the fiber.} We introduced the measures $\mu^{F\,\pm}_{x,R,\eps}$ with the hope to compare them (their projections on $X_p$ to be more precise) with $\theta_{F,R}$. Hence, let us consider $\theta^{\pm}_{F,R,\eps}=pr_{\eps}\,_{\ast}\theta^{\pm}_{F,R}$ the projections on the fiber $X_p$.

\begin{lemma}
\label{comparisontheta}
There is a number $C(\eps)\geq 1$ which tends to $1$ as $\eps$ tends to $0$ such that the following chain of inequalities (with notations defined in \eqref{Eq:integralweight} and \eqref{Eq:sumweight}) holds true
$$C(\eps)^{-1}\frac{\theta^-_{F,R,\eps}}{\Leb(B(p,\eps))}\leq \frac{J_{F,R}}{I_{F,R}}\,\theta_{F,R}\leq C(\eps)\frac{\theta^+_{F,R,\eps}}{\Leb(B(p,\eps))}.$$
\end{lemma}

\begin{proof}
By definition, $\theta_{F,R}$ is obtained as the projection on $M$ via $\proj_x$ of the averages of the Dirac masses located at $B_R.o$ weighted by $\kappa^F(o,\gamma o)$.

By definition again, $\theta^{\pm}_{F,R,\eps}$ can be described in a similar way. Consider the measure on $N$ which gives to any point $\gamma o$, with $\gamma\in B^{\pm}_{R,\eps}$, the weight $\int_{B(\gamma o,\eps)}\kappa^F(o,y)d\Leb(y)$. Now normalize this measure by $I_{R,F}$ and project it down to $M$ via $\proj_x$.

Firstly, in order to prove the lemma, note that since by definition we have the chain of inclusions $B^-_{R,\eps}\dans B_R\dans B^+_{R,\eps}$, we have $\Supp\theta^-_{F,R,\eps}\dans\Supp\theta_{F,R}\dans\Supp\theta^+_{F,R,\eps}$. .

Now, we have to compare the weights that these measures give to points of their supports. By distortion lemma \ref{distortioncomparison}, there is a number $C(\eps)$ which tends to $1$ as $\eps$ tends to $0$ such that if $y,z\in N$ are distant of at most $\eps$, the quotient $\kappa^F(o,y)/\kappa^F(o,z)$ belongs to $[C(\eps)^{-1},C(\eps)]$. Finally, for all $\gamma\in\pi_1(B)$, the integral with respect to Lebesgue of $\kappa^F(o,.)$ on the ball centered at $\gamma o$ and of radius $\eps$ is, up to $C(\eps)$, close to the value at $\gamma o$ times the volume of the ball. 

Since the volumes of all balls of radius $\eps$ are equal, we can conclude the proof of the lemma by making the suitable normalizations.
\end{proof}

\begin{lemma}
\label{somebounds}
The following assertions hold true.
\begin{enumerate}
\item $I_{R+\eps,\eps}/I_{R,\eps}\to e^{P(F)\eps}$ as $R$ tends to infinity;
\item For $\eps>0$ small enough, we have
$$C(\eps)^{-1}I_{F,R-2\eps}\frac{\mass(\mu^F_{x,R-2\eps,\eps})}{\Leb(B(p,\eps))}\leq J_{F,R}\leq C(\eps)I_{F,R+2\eps}\frac{\mass(\mu^F_{x,R+2\eps,\eps})}{\Leb(B(p,\eps))}.$$
\end{enumerate}
\end{lemma}

\begin{proof}
The first assertion is an immediate consequence of Ledrappier's theorem \ref{kappaspheres}. The second one is obtained by evaluating the masses in Lemmas \ref{comparisontheta} and \ref{chainofinequalities}.
\end{proof}
\paragraph{Remark.} Since, $m_F(\partial K_{p,\eps})=0$, we see, by Lemma \ref{convergenceinrestriction}, that as $R$ grows to infinity, $\mass(\mu^F_{x,R,\eps})$ converges to $m_F(K_{p,\eps})$ which is equal to $\int_{B(p,\eps)}h_0\,d\Leb$.

\paragraph{A key lemma} The next lemma is the key step for proving the convergence of the family $(\theta_{F,R})_{R>0}$. For a continuous function $f:X_p\to\R$, we will pinch its integral against $\theta_{F,R}$ between two quantities which become eventually very close. In what follows, we will use the following notation
$$\Lambda(p)=\int_{B(p,\eps)} h_0\,d\Leb.$$
\begin{lemma}
\label{pinchinglemma}
There is a $C'(\eps)>1$ which tends to $1$ as $\eps$ tends to $0$, as well as a constant $R_0>0$, such that for any continuous function $f:X_p\to\R$, any $R>R_0$ and any small $\eps>0$,
\begin{equation}
\label{encadrement}
C'(\eps)^{-1}\int_{K_{p,\eps}}f\circ pr_{\eps}\,\frac{d\mu^F_{x,R-2\eps,\eps}}{\mass(\mu^F_{x,R-2\eps,\eps})}\leq\int_{X_p}f d\theta_{F,R}\leq C'(\eps)\int_{K_{p,\eps}}f\circ pr_{\eps}\,\frac{d\mu^F_{x,R+2\eps,\eps}}{\mass(\mu^F_{x,R+2\eps,\eps})}
\end{equation}
\end{lemma}

\begin{proof}
We only prove the upper bound: the lower bound follows by the exact same argument. Recall that by definition, $\theta^+_{F,R,\eps}$ is the projection via $pr_{\eps}$ of $\mu^{F\,+}_{x,R,\eps}$. So for any continuous function $f:X_p\to\R$, we have $\int_{X_p} f d\theta^+_{F,R,\eps}=\int_{K_{p,\eps}} f\circ pr_{\eps}d\mu^{F\,+}_{x,R,\eps}$. But now, we can combine Lemmas \ref{chainofinequalities} and \ref{comparisontheta} in order to prove the following inequality, which is valid for any $f$, $R$ and $\eps>0$ sufficiently small
$$\int_{X_p} f d\theta_{F,R}\leq C(\eps)\frac{I_{F,R+2\eps}}{J_{F,R}}\int_{K_{p,\eps}} f\circ pr_{\eps}\,\frac{d\mu^F_{x,R+2\eps,\eps}}{\Leb(B(p,\eps))}.$$
If one uses Lemma \ref{somebounds}, there are two consequences. Firstly, we can use the lower bound of $J_{F,R}$, in order to have the following upper bound
$$\int_{X_p} f d\theta_{F,R}\leq C(\eps)^2\frac{I_{F,R+2\eps}}{I_{F,R-2\eps}}\int_{K_{p,\eps}} f\circ pr_{\eps}\,\frac{d\mu^F_{x,R+2\eps,\eps}}{\mass(\mu^F_{x,R-2\eps,\eps})}.$$

Secondly, we know that $I_{F,R+2\eps}/I_{F,R-2\eps}\to e^{4P(F)\eps}$ as $R$ tends to infinity.

Finally, as explained in the previous remark, when $R\to\infty$, $\mass(\mu^F_{x,R,\eps})\to \Lambda(p)$, so, when $R$ is large enough, we can indistinctly divide the upper bound by $\mass(\mu^F_{x,R-2\eps,\eps})$, or by $\mass(\mu^F_{x,R+2\eps,\eps})$. The existence of a large $R_0>0$ and of a $C'(\eps)$ such that the upper bound of the lemma holds, then follows concluding the proof of the lemma.
\end{proof}

\paragraph{Convergence of the measures.}We are now ready to prove Theorem \ref{countingmeasures}.
\begin{proposition}
The measures $\theta_{F,R}$ converge to the probability measure $m_{F,p}$.
\end{proposition}

\begin{proof}
By compactness of the fibers, one can consider a sequence $(R_n)_{n\in\N}$ such that the sequence $(\theta_{F,R_n})_{n\in\N}$ converges to some $\theta_{F,\infty}$. We must prove that $\theta_{F,\infty}=m_{F,p}$, which will prove the  uniqueness of accumulation point, and thus the convergence, of the family $(\theta_{F,R})_{R>0}$ towards $m_{F,p}.$

For any continuous function $f:X_p\to\R$, we have
$$\lim_{n\to\infty}\int_{X_p}fd\theta_{F,R_n}=\int_{X_p}fd\theta_{F,\infty}.$$

Recall also that by definition, and since the family of conditional measures of $m_F$ is continuous ($F$-harmonic functions are continuous), the conditional measure of $m_F$ on $X_p$ satisfies for any continuous $f:X_p\to\R$
$$\int_{X_p}f\,dm_{F,p}=\lim_{\eps\to 0}\int_{K_{p,\eps}}f\circ pr_{\eps}\,\frac{dm_F}{m_F(K_{p,\eps})}.$$
We will use of Lemma \ref{pinchinglemma}. Write Formula \eqref{encadrement} with $R_n$, and some small $\eps>0$. Now let $n$ go to infinity, and then $\eps$ go to zero. We end up with $\int f\,d\theta_{F,\infty}=\int f\,dm_{F,p}$ for any continuous function $f:X_p\to\R$. The proof of the proposition and hence that of Theorem \ref{countingmeasures} is now over.
\end{proof}

\section{Examples of quasifuchsian representations}
\label{sfuchsianquasifuchsian}
In this section, we will be interested in discrete and faithful fuchsian and quasifuchsian representations of the fundamental group of a surface with negative curvature. These examples are by far the most simple ones, but we think it is worth explaining how we treat them.

\subsection{Quasifuchsian representations}

\paragraph{Boundary correspondence.} Let $\Sigma$ be a closed surface of genus $\geq 2$ and $g_0$ be a Riemannian metric with negative curvature (a priori variable). The Riemannian universal cover $\widetilde{\Sigma}_0$ is a disc which can be compactified by a $C^1$ circle $\Sigma_0(\infty)$ (the circle at infinity of $\Sigma$). This circle inherits a natural orientation. We also have a discrete and faithful representation $\rho_0:\pi_1(\Sigma)\to\Isom^+(\widetilde{\Sigma}_0)$: we have a natural action of $\pi_1(\Sigma)$ on $\Sigma_0(\infty)$ by orientation preserving diffeomorphisms.

Denote by $\Sigma_0^{(3)}(\infty)$ the set of oriented triples $(\xi_+,\xi_0,\xi_-)\in\left(\Sigma_0(\infty)\right)^3$. Since $\Gamma_0$ preserves orientation, the group $\pi_1(\Sigma)$ acts diagonally on $\Sigma_0^{(3)}(\infty)$. There is a natural identification $T^1\widetilde{\Sigma}_0\to\Sigma_0^{(3)}(\infty)$ which associates to any vector $v$ the triple $(pr_+(v),pr_0(v),pr_-(v))$, where:
\begin{itemize}
\item $pr_+(v)\in\Sigma_0(\infty)$ is the limit point of the geodesic ray determined by $-v$;
\item $pr_-(v)\in\Sigma_0(\infty)$ is the limit point of the geodesic ray determined by $v$;
\item $pr_0(v)\in\Sigma_0(\infty)$ is the limit point of the geodesic orthogonal to $v$ wich satisfies $pr_+(v)<pr_0(v)<pr_-(v)$ for the orientation.
\end{itemize}
This identification is an equivariance (recall that $\pi_1(\Sigma)$ acts on $T^1\widetilde{\Sigma}_0$ by differential of elements of $\Gamma_0$ and diagonally on $\Sigma^{(3)}_0(\infty)$). Moreover, the geodesic determined by $v$ can be parametrized by a point of the arc $[pr_+(v),pr_-(v)]$.

Let $g_1$ be a another negatively curved metric on $\Sigma$. Denote by $\widetilde{\Sigma}_1$ the universal cover and by $\Sigma_1(\infty)$ the sphere at infinity.

The \emph{boundary correspondence} is the bihölder homeomorphism $h:\Sigma_0(\infty)\to\Sigma_1(\infty)$ which conjugates the actions at infinity given by $\rho_0$ and $\rho_1$ (i.e. $h\circ\rho_0(\gamma)=\rho_1(\gamma)\circ h$ for $\gamma\in\pi_1(\Sigma)$). We refer to Section 5.9 of Thurston's notes \cite{T} for the existence and properties of this object.

\paragraph{Quasifuchsian representations.}

Take two hyperbolic metrics $g_-$ and $g_+$ on $\Sigma$: they can be uniformized by copies $\Gamma_-$ and $\Gamma_+$ of $\pi_1(\Sigma)$ inside $\PSL_2(\R)$. By a Bers' simultaneous uniformization \cite{Be}, there exist
\begin{itemize}
\item a discrete subgroup $\Gamma< \PSL_2(\C)$ which leaves invariant a Jordan curve $\Lambda$, as well as the two connected components $D_+$ and $D_-$ of $\C\PP^1\moins\Lambda$;
\item two analytic maps from the upper and lower half planes $H_{\pm}:\pm\Hyp^2\to D_{\pm}$ which conjugate $\Gamma_{\pm}$ and $\Gamma$. In particular, there exists an isomorphism $\rho:\pi_1(\Sigma)\to\Gamma$ and all elements of $\Gamma$ are hyperbolic;
\item $H_{\pm}$ extend as homeomorphisms $\pm\Hyp^2\cup\R\PP^1\to D_{\pm}\cup\Lambda$ which are equivariant. In particular, $(H_-)^{-1}\circ H_+:\R\PP^1\to\R\PP^1$ is a boundary correspondence between $(\Sigma,g_+)$ and $(\Sigma,g_-)$.
\end{itemize}
These objects are uniquely defined up to conjugacy by elements of $\PSL_2(\C)$. The group $\Gamma$ is called a \emph{quasifuchsian group} and the isomorphism $\rho:\pi_1(\Sigma)\to\Gamma<\PSL_2(\C)$ is called a \emph{quasifuchsian representation}.

Bowen showed in \cite{Bo3} that either $\Lambda$ is a geometric circle, or it has Hausdorff dimension $>1$. In the first case, $H_{\pm}$ are elements of $\PSL_2(\R)$, so $\Gamma_+$ and $\Gamma_-$ are conjugated by an isometry of $\Hyp^2$ and the two metrics $g_+$ and $g_-$ are in the same Teichmüller class. This corresponds to the \emph{Fuchsian case}.

Consider such a group $\Gamma$: it comes with objects $\Gamma_{\pm},\Lambda$ and $H_{\pm}$ described above. As explained in the previous paragraph, we have two boundary correspondences $h_{\pm}:\Sigma(\infty)\to\R\PP^1$ which conjugate the actions of $\Gamma_0=\rho_0(\pi_1(\Sigma))$ and $\Gamma_{\pm}$. Since these objects are defined up to conjugacy, we can assume that the two boundary correspondences $h_-\circ(h_+)^{-1}$ and $(H_-)^{-1}\circ H_+$ coincide (we then have $H_+\circ h_+=H_-\circ h_-$). 

\paragraph{The Lyapunov sections.} We consider a negatively curved Riemannian metric $g_0$ on $\Sigma$ and a quasifuchsian representation given by the isomorphism
$$\rho:\pi_1(\Sigma)\to\Gamma< \PSL_2(\C).$$
Consider the two associated foliated bundles $(\Pi,M,\Sigma,\C\PP^1,\F)$ and $(\Pi_{\ast},T^1\F,T^1\Sigma,\C\PP^1,\widehat{\F})$ with the lifted metric in the leaves. Since all elements of $\Gamma$ are hyperbolic isometries of the upper half space $\Hyp^3$ and since the group $\Gamma_0$ is not a copy of $\Z$, this representation has no invariant measure on $\C\PP^1$.

We can define three Hölder continuous sections $\sigma^{\star}:T^1\Sigma\to T^1\F$, $\star=+,-,0$. First use the boundary correspondence to identify $T^1\widetilde{\Sigma}_0$ with $\widetilde{\Sigma}_0^{(3)}$.

Now the three maps $\widetilde{s}^{\star}=H_+\circ h_+\circ pr_{\star}=H_-\circ h_-\circ pr_{\star}:T^1\widetilde{\Sigma}\to\C\PP^1$, for $\star=+,-,0$ are equivariant. The corresponding graphs $\widetilde{\sigma}^{\star}:T^1\widetilde{\Sigma}_0\to T^1\widetilde{\Sigma}_0\times\C\PP^1$ descend to the quotient and produce the desired sections. By construction, they commute with the geodesic flows, $\sigma^+$ commutes with the strong unstable foliations and $\sigma^-$ with the strong stable ones. These two sections are the two Lyapunov sections.

\subsection{Measures associated to fuchsian and quasifuchsian representations}

Let $F:T^1\Sigma\to\R$ be a Hölder continuous potential. Assume that 
$$\rho:\pi_1(\Sigma)\to\Gamma<\PSL_2(\C)$$
is a quasifuchsian representation. We obtain a foliated bundle $(\Pi,M,\Sigma,\C\PP^1,\F)$ by suspension of $\rho$. The unit tangent bundle is also a foliated bundle $(\Pi_{\ast},T^1\F,T^1\Sigma,\C\PP^1,\widehat{\F})$ with the same holonomy. By Theorem \ref{tripleuniquegibbs}, there exists a unique $F$-harmonic measure for $\F$.

We want to compare these measures for different potentials $F$. We know that each of them projects down onto a measure equivalent to Lebesgue: we have to compare the conditional measures on the fibers. We know by Theorem \ref{disintegrationfharmonic} that the conditional measures are equivalent to $s_p^+\,_\ast\omega^F_p$, where $(\omega^F_p)_{p\in\Sigma}$ is the Ledrappier family on the unit tangent fibers and $s^+_p$ is the map induced by the Lyapunov section: it is a bihölder homeomorphism on its image, which is identified to a Jordan curve $\Lambda$.

We show a special interest in the three classes on $\Sigma_0(\infty)$ which are (see \cite{Ka1}, \cite{Ka2}, \cite{L2}):
\begin{itemize}
\item the harmonic class which describes the distribution at infinity of Brownian paths;
\item the visibility class which describes the distribution of Lebesgue-almost all geodesics;
\item the Patterson-Sullivan class which describes the behaviour at infinity of the $\Gamma_0$-orbits ($\Gamma_0<\Isom^+(\widetilde{\Sigma}_0)$ being a copy of $\pi_1(\Sigma)$ which uniformizes the surface).
\end{itemize}

As a combination of the works of Katok and Ledrappier (see \cite{Ka1}, \cite{Ka2}, \cite{L3}), we get that:

\begin{theorem}[Katok, Ledrappier]
\label{differentclasses}
The curvature of $(\Sigma,g_0)$ is constant if and only if two of the three classes of measures (harmonic, visibility and Patterson-Sullivan) coincide. In this case, the three classes of measures are the same.
\end{theorem}

Thus, this theorem allows us to finish the proof of Theorem \ref{tripleuniqueergodicity}: if we suspend a quasifuchsian representation of the fundamental group of a Riemannian surface with \emph{variable} negative curvature, then the unique harmonic measure, the projection of the unique SRB measure and the unique limit of large discs are mutually singular.

\subsection{Two questions}
\label{somequestions}

We would like to finish this papers by addressing two questions.

\paragraph{Gibbs measures for foliated hyperbolicity.} In this article, we only define Gibbs measures in the case of foliated bundles. For the geodesic flow tangent to general foliations, we were only able to define Gibbs measures for two particular potentials: we studied the $H$-Gibbs measures (see \cite{Al2}), as well as the Gibbs $u$-states (see \cite{Al3}). The reason we were able to do it is that is that harmonic and Lebesgue classes on horospheres exist even if the leaves don't cover a compact manifold.

If one knew, for some potential $F$, how to define a family of measures in the unstable manifolds satisfying the right quasi-invariance properties under the action of the foliated geodesic flow, one could call $F$-Gibbs an invariant measure whose conditional measures in the unstable manifolds lie in the class prescribed by the family.

In the case of foliated bundles, we only have to lift the families of measures which are defined in the base by the fibration. In the general case, it is open.
\begin{ques}
Let $(M,\F)$ be a compact foliated manifold with negatively curved leaves, and $F:T^1\F\to\R$ a continuous function, uniformly Hölder continuous in the leaves of $\hcF$. Is it possible to construct a family of measures $(\lambda^F_{F,v})_{v\in T^1\F}$ defined in the unstable manifolds which is absolutely continuous under the action of the foliated geodesic flow with the following Radon-Nikodym cocycle
$$\frac{d\left[G_T\,_\ast\lambda^{u}_{F,g_{-T}(v)}\right]}{d\lambda^{u}_{F,v}}(v)=\exp\left[\int_0^T (F\circ G_{-t}(v)-P(F))dt\right]?$$
If the answer is affirmative, what is the dynamical meaning of the number $P(F)$? In the case of the null potential, is it related with the foliation entropy \cite{GLW}?
\end{ques}

\paragraph{The dimension of fibers.} The results of uniqueness stated here are in the context of a projective fiber of complex dimension $1$. It is natural to ask the following question. 

\begin{ques}
Do the uniqueness results stated here hold when the fiber is a $\C\PP^{d}$, $d\geq 2$?
\end{ques}
The case of $\C\PP^2$ could easily be derived from a trick due to Bonatti and G\'omez-Mont described in the last part of \cite{BG}. The question of the unique ergodicity of the foliated horocyclic flow when the base is a closed hyperbolic surface, and the fiber is $\C\PP^{d}$, $d\geq 3$ is a work in progress by Bonatti, Eskin and Wilkinson (private conversation with Bonatti).

\vspace{10pt}
\paragraph{Acknowledgments.} I would like to thank warmly my advisor Christian Bonatti who read the first drafts of this paper with care, dedication and a lot of patience. It is also a pleasure to thank François Ledrappier, Barbara Schapira, Rémi Langevin and Patrick Gabriel for many useful conversations and mail exchanges. Last but not least, the present version of the paper owes a lot to the referee's comments. I am grateful to him/her for that.

\vspace{10pt}

\noindent \textbf{Sébastien Alvarez (salvarez@impa.br)}\\
\noindent  IMPA, Est. D. Castorina 110, 22460-320, Rio de Janeiro, Brazil.\\


\begin{thebibliography}{widest-label}

\bibitem[AR]{AR} F. Alcalde Cuesta, A. Rechtman, Averaging sequences, \emph{Pac. J. of Math.}, \textbf{255}, (2012), 1-23.

\bibitem[Al1]{Al1} S. Alvarez, Discretization of harmonic measures for foliated bundles, \emph{C. R. Acad. Sci. Paris, Ser I.}, \textbf{350}, (2012), 621-626.

\bibitem[Al2]{Al2} S. Alvarez, Harmonic measures and the foliated geodesic flow for foliations with negatively curved leaves, \emph{Ergod. Th. \& Dynam. Sys.}, \textbf{36}, (2016), 355-374.

\bibitem[Al3]{Al3} S. Alvarez, Gibbs u-states for the foliated geodesic flow and transverse invariant measures, \emph{preprint} [arXiv:1311.7121].

\bibitem[Al4]{Al4} S. Alvarez, A Fatou theorem for $F$-harmonic functions, \emph{preprint} [arXiv:1512.01842].

\bibitem[Al5]{Al5} S. Alvarez, Mesures de Gibbs et mesures harmoniques pour les feuilletages aux feuilles courbées négativement, \emph{Thèse de l'Université de Bourgogne}, (2013).

\bibitem[AY]{AY}S. Alvarez, J. Yang, Physical measures for the geodesic flow tangent to a transversally conformal foliation, \emph{preprint} [arXiv:1512.01842].


\bibitem[An]{An} D. Anosov, Geodesic flows on closed Riemannian manifolds with negative curvature, \emph{Proc. Steklov Math. Inst .A.M.S Transl.}, (1969).

\bibitem[AS]{AS} M.T. Anderson, R. Schoen, Positive harmonic functions on complete manifolds of negative curvature, \emph{Ann. of Math.}, \textbf{121}, (1985), 429-461.

\bibitem[ABV]{ABV} J. Alves, C. Bonatti, M. Viana, SRB measures for partially hyperbolic systems whose
central direction is mostly expanding, \emph{Invent. Math.}, \textbf{140}, (2000), 351-398.

\bibitem[Ba]{Ba} M. Babillot, On the mixing property for hyperbolic systems, \emph{Israel J. of Math.}, \textbf{129}, (2002), 61-76.

\bibitem[BL]{BL} M. Babillot, F. Ledrappier, Geodesic paths and horocycle flows on abelian covers, \emph{Lie Groups and ergodic theory (Mumbai, 1996)}, in Tata Inst. Fund. Res. Stud. Math, \textbf{14}, (1998), 1-32.

\bibitem[BMar]{BMar} Y. Bakhtin, M. Mart\'inez, A characterization of harmonic measures on laminations by hyperbolic Riemann surfaces, \emph{Ann. I. H.Poincaré. Prob. Stat.}, \textbf{44}, (2008), 1078-1089.

\bibitem[BY]{BY} M. Benedicks, L.S. Young, SRB-measures for certain H\'enon maps, \emph{Invent. Math.}, \textbf{112}, (1993), 541-576.

\bibitem[Be]{Be} L. Bers, Simultaneous Uniformization, \emph{Bull. Amer. Math. Soc.}, \textbf{66}, (1960), 94-97.

\bibitem[BDV]{BDV} C. Bonatti, L. D\'iaz, M. Viana, \emph{Dynamics Beyond Uniform Hyperbolicity. A global geometric and probabilistic perspective}, Encyclopaedia of Mathematical Sciences, 102. Mathematical Physics, III. Springer Verlag, 2005.

\bibitem[BG]{BG} C. Bonatti, X. G\'omez-Mont, Sur le comportement statistique des feuilles de certains feuilletages holomorphes, \emph{Monogr. Enseign. Math.}, \textbf{38}, (2001), 15-41.

\bibitem[BGM]{BGM} C. Bonatti, X.G\'omez-Mont, M.Mart\'inez, Foliated hyperbolicity and foliations with hyperbolic leaves, \emph{preprint} [arXiv:1510.05026].

\bibitem[BGVil]{BGVil} C. Bonatti, X. G\'omez-Mont, R. Vila, Statistical behaviour of the leaves of Riccati foliations, \emph{Ergod. Th. \& Dynam. Sys.}, \textbf{30}, (2010), 67-96.

\bibitem[BGV]{BGV} C. Bonatti, X. G\'omez-Mont, M. Viana, G\'en\'ericit\'e d'exposants de Lyapunov non-nuls pour des produits d\'eterministes de matrices, \emph{Ann. I. H. Poincar\'e. Anal. Non Lin.}, \textbf{20}, (2003), 579-624.

\bibitem[BV]{BV} C. Bonatti, M. Viana, SRB measures for partially hyperbolic systems whose central direction is mostly contracting, \emph{Israel J. Math.}, \textbf{115}, (2000), 157-193.

\bibitem[Bo1]{Bo1} R. Bowen, Symbolic dynamics for hyperbolic flows, \emph{Amer. J. of Math.}, \textbf{95}, (1973), 429-459.

\bibitem[Bo2]{Bo2} R. Bowen, \emph{Equilibrium states and the ergodic theory of Anosov diffeomorphisms}, Lect. Notes Math, \textbf{470}, Springer-Verlag, 1975.

\bibitem[Bo3]{Bo3} R. Bowen, Hausdorff dimension of quasi-circles, \emph{Inst. Hautes Études Sci. Publ. Math.}, \textbf{50}, (1979), 259-273.

\bibitem[BM]{BM} R. Bowen, B. Marcus, Unique ergodicity for horocycle foliations, \emph{Israel J. of Math.}, \textbf{26}, (1977), 43-67.

\bibitem[BR]{BR} R. Bowen, D. Ruelle, The ergodic theory of Axiom A flows, \emph{Invent. Math.}, \textbf{29}, (1975), 181-202.

\bibitem[BH]{BH} M. Bridson, A. Haefliger, \emph{Metric spaces of non-positive curvature}, Grundlehren Math. Wiss., \textbf{319}, Springer-Verlag, 1999.

\bibitem[CL]{CL} C. Camacho, A. Lins Neto: \emph{Geometric Theory of Foliations}, Birkhäuser, Boston Inc., 1985.

\bibitem[Car]{Car} M. Carvalho, Sinai-Ruelle-Bowen measures for $n$-dimensional derived from Anosov diffeomorphisms, \emph{13}, (1993), 21-44.

\bibitem[DK]{DK} B. Deroin, V. Kleptsyn, Random conformal dynamical systems, \emph{Geom. Funct. Anal.}, \textbf{17}, (2007), 1043-1105.

\bibitem[Gar]{Gar} L. Garnett, Foliations, the ergodic theorem and Brownian motion, \emph{J. Funct. Anal.}, \textbf{51}, (1983), 285-311.

\bibitem[GLW]{GLW} E. Ghys, R. Langevin, P. Walczak, Entropie géométrique des feuilletages, \emph{Acta Math.}, \textbf{160}, (1988), 105-142.

\bibitem[GPl]{GPl} S. Goodman, J. Plante, Holonomy and averaging in foliated sets, \emph{J. Diff. Geo.}, \textbf{14}, (1979), 401-407.

\bibitem[HPS]{HPS} M. Hirsch, C. Pugh, M. Shub, \emph{Invariant manifolds}, in Lecture Notes in Math., \textbf{583}, Springer Verlag, 1977.

\bibitem[K1]{K1} V. Kaimanovich, Invariant measures of the geodesic flow and measures at infinity on negatively curved manifolds, \emph{Ann. I. H. Poincar\'e. Phys. Théor.}, \textbf{53}, (1990), 361-393.

\bibitem[K2]{K2} V. Kaimanovich, Amenability, hyperfiniteness, and isoperimetric inequalities, \emph{C. R. Acad. Sci. Paris, Ser I.}, \textbf{325}, (1997), 999-1004.

\bibitem[KL]{KL} V. Kaimanovich, M. Lyubich, Conformal and harmonic measures on laminations associated with rational maps, \emph{Mem. Amer. Math. Soc.}, \textbf{173}, (2005).

\bibitem[Ka1]{Ka1} A. Katok, Entropy and closed geodesics,  \emph{Ergod. Th. \& Dynam. Sys.}, \textbf{2}, (1982), 339-365.

\bibitem[Ka2]{Ka2} A. Katok, Four applications of conformal equivalence to geometry and dynamics, \emph{Ergod. Th. \& Dynam. Sys.}, \textbf{8}, (1988), 115-140.

\bibitem[Kn]{Kn} G. Knieper, Spherical means on compact manifolds of negative curvature, \emph{J. Diff. Geom. Appl.}, \textbf{4}, (1994), 361-390.

\bibitem[KP]{KP} G. Knieper, N. Peyerimhoff, Ergodic properties of isoperimetric domains in spheres, \emph{J. Mod. Dyn.}, \textbf{2}, (2008), 239-358.

\bibitem[L1]{L1} F. Ledrappier, Ergodic properties of Brownian motion on covers of compact negatively curved manifolds, \emph{Bol. Soc. Bras. Mat.}, \textbf{19}, (1988), 115-140.

\bibitem[L2]{L2} F. Ledrappier, Harmonic measures and Bowen-Margulis measures, \emph{Israel J. of Math.}, \textbf{71}, (1990), 275-287.

\bibitem[L3]{L3} F. Ledrappier, Structure au bord des variétés à courbure négative, \emph{Sémin. de Th. Spec. et Géom., Grenoble}, (1994-1995), 93-118.

\bibitem[L4]{L4} F. Ledrappier, A renewal theorem for the distance in negative curvature, \emph{Proc. Symp. Pure Math.}, \textbf{57}, (1995), 351-360.

\bibitem[Man]{Man} R. Ma\~né, \emph{Ergodic theory and differentiable dynamics}, Springer Verlag, 1987.

\bibitem[M]{M} G. Margulis, Certain measures associated with U-flows on compact manifolds, \emph{Func. Anal. Appl.}, \textbf{4}, (1970), 55-67.

\bibitem[Ma]{Ma} M. Mart\'inez, Measures on hyperbolic surface laminations, \emph{Ergod. Th. \& Dynam. Sys.}, \textbf{26}, (2006), 847-867.

\bibitem[Os]{Os} V. Oseledets, A multiplicative ergodic theorem: Lyapunov characteristic numbers for dynamical systems, \emph{Trans. Moscow. Math. Soc.}, \textbf{19}, (1968), 197-231.

\bibitem[Pat]{Pat} S.J. Patterson, The limit set of a Fuchsian group, \emph{Acta. Math.}, \textbf{136}, (1976), 241-273.

\bibitem[PPS]{PPS} F. Paulin, M. Pollicott, B. Schapira, \emph{Equilibrium states in negative curvature}, Astérisque, \textbf{373}, 2015.

\bibitem[PS]{PS} Ya. Pesin, Ya .Sinai, Gibbs measures for partially hyperbolic attractors, \emph{Ergod. Th. \& Dynam. Sys.}, \textbf{2}, (1982), 417-438.

\bibitem[Ra]{Ra} M. Ratner, Markov partitions for Anosov flows on $n$-dimensional manifolds, \emph{Israel J. of Math.}, \textbf{15}, (1973), 92-114.

\bibitem[Rob]{Rob} T. Roblin, \emph{Ergodicité et équidistribution en courbure négative}, \emph{Mem. Soc. Math. Fr.}, \textbf{95}, 2003. 

\bibitem[Ro]{Ro} V.A. Rokhlin, On the fundamental ideas of measure theory, \emph{Amer. Math. Soc. Transl.}, \textbf{10}, (1962), 1-52.

\bibitem[Sc1]{Sc1} B. Schapira, On quasi-invariant transverse measures for the horospherical foliation of a negatively curved manifold, \emph{Ergod. Th. \& Dynam. Sys.}, \textbf{24}, (2004), 227-257.

\bibitem[Sc2]{Sc2} B. Schapira, Mesures quasi-invariantes pour un feuilletage et limites de moyennes longitudinales, \emph{C. R. Acad. Sci. Paris, Ser I.}, \textbf{336}, (2003), 349-352.

\bibitem[Sh]{Sh} M. Shub, \emph{Stabilité globale des systèmes dynamiques}, in Astérisque, \textbf{56}, Soc. Math. France, 1978.

\bibitem[Si]{Si} Ya. Sinai, Gibbs measures in ergodic theory, \emph{Russ. Math. Surveys}, \textbf{27}, (1972), 21-69.

\bibitem[T]{T} W. Thurston, Geometry and topology of 3-manifolds, Princeton Lecture Notes, 1980, currently available at library.msri.org/books/gt3m/

\bibitem[Y]{Y} L.S. Young, Large deviations in dynamical systems, \emph{Trans. Amer. Math. Soc.}, \textbf{318}, (1990), 525-543.

\end{thebibliography}
\end{document}